\documentclass[12pt,twoside, epsf]{article}

\usepackage{color,graphicx,times}

\usepackage{amssymb}
\usepackage{amsmath}
\usepackage{amscd}
\usepackage{float}
\usepackage{graphicx}
\usepackage{amsfonts}
\usepackage{pb-diagram}
\usepackage{eufrak}

\newcommand{\CrossStarGlyph}{\raisebox{-0.25\height}{\includegraphics[width=0.5cm]{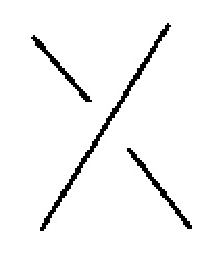}}}
 \newcommand{\CrossGlyph}{\raisebox{-0.25\height}{\includegraphics[width=0.5cm]{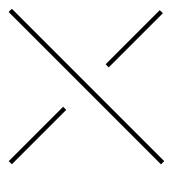}}}
\newcommand{\HSmoothGlyph}{\raisebox{-0.25\height}{\includegraphics[width=0.5cm]{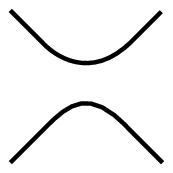}}}
\newcommand{\VSmoothGlyph}{\raisebox{-0.25\height}{\includegraphics[width=0.5cm]{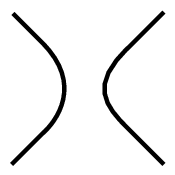}}}
\newcommand{\VirtualGlyph}{\raisebox{-0.25\height}{\includegraphics[width=0.5cm]{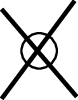}}}

 \newcommand{\CGlyph}{\raisebox{-0.25\height}{\includegraphics[width=0.5cm]{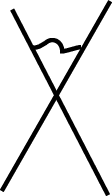}}}
\newcommand{\VGlyph}{\raisebox{-0.25\height}{\includegraphics[width=0.5cm]{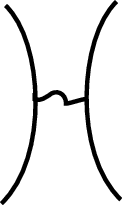}}}
\newcommand{\YGlyph}{\raisebox{-0.25\height}{\includegraphics[width=0.5cm]{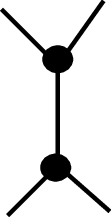}}}
\newcommand{\YMGlyph}{\raisebox{-0.25\height}{\includegraphics[width=0.5cm]{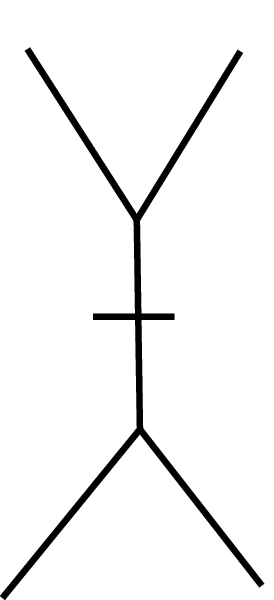}}}
\newcommand{\IMGlyph}{\raisebox{-0.25\height}{\includegraphics[width=0.5cm]{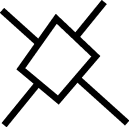}}}
\newcommand{\IMGlyphL}{\raisebox{-0.25\height}{\includegraphics[width=0.5cm]{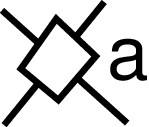}}}

\newcommand{\YMTens}{\raisebox{-0.25\height}{\includegraphics[width=0.5cm]{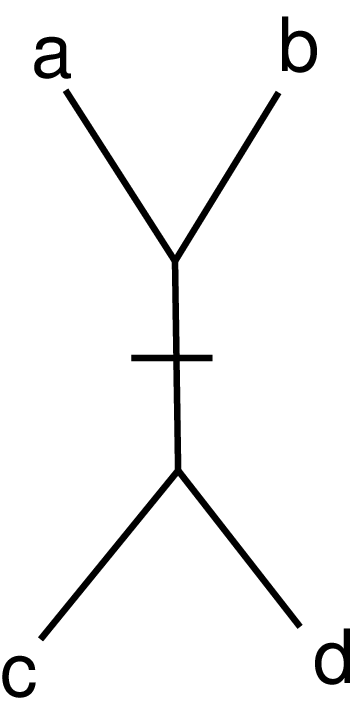}}}
\newcommand{\IMTens}{\raisebox{-0.25\height}{\includegraphics[width=0.5cm]{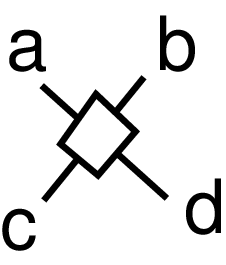}}}
\newcommand{\CTens}{\raisebox{-0.25\height}{\includegraphics[width=0.5cm]{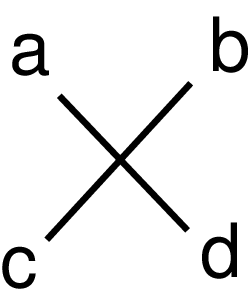}}}
\newcommand{\FTens}{\raisebox{-0.25\height}{\includegraphics[width=0.5cm]{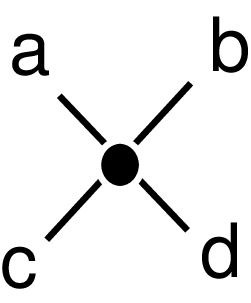}}}

\newcommand{\CDiag}{\raisebox{-0.25\height}{\includegraphics[width=0.5cm]{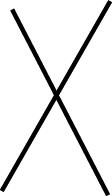}}}
\newcommand{\VDiag}{\raisebox{-0.25\height}{\includegraphics[width=0.5cm]{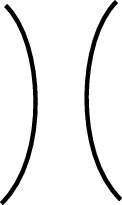}}}
\newcommand{\CDotDiag}{\raisebox{-0.25\height}{\includegraphics[width=0.5cm]{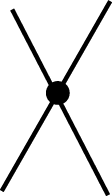}}}
\newcommand{\CCircleDiag}{\raisebox{-0.25\height}{\includegraphics[width=0.5cm]{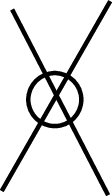}}}
\newcommand{\CBoxDiag}{\raisebox{-0.25\height}{\includegraphics[width=0.5cm]{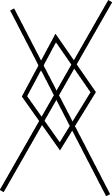}}}
\newcommand{\CCircleDiagL}{\raisebox{-0.25\height}{\includegraphics[width=0.5cm]{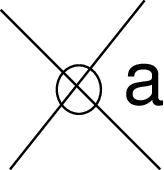}}}

\newcommand{\Horiz}{\raisebox{-0.25\height}{\includegraphics[width=0.5cm]{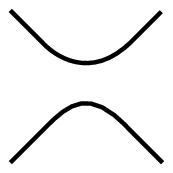}}}
\newcommand{\Virt}{\raisebox{-0.25\height}{\includegraphics[width=0.5cm]{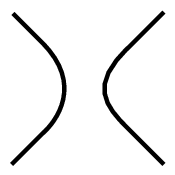}}}

\newcommand{\Wigg}{\raisebox{-0.25\height}{\includegraphics[width=0.9cm]{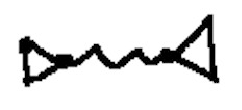}}}
\newcommand{\Contract}{\raisebox{-0.25\height}{\includegraphics[width=0.9cm]{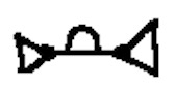}}}
\newcommand{\Edge}{\raisebox{-0.25\height}{\includegraphics[width=0.9cm]{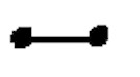}}}
\newcommand{\Delete}{\raisebox{-0.25\height}{\includegraphics[width=0.9cm]{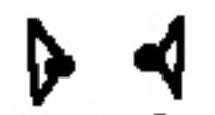}}}

\makeatother

\def\o{\overline}

\let \ttorg \tt \def \tt{\ttorg \obeyspaces}

\begin{document}

\date{}

\title{\bf Multi-Virtual Knot Theory}

\author{Louis H. Kauffman \\ 
Department of Mathematics, Statistics and Computer Science \\ 
University of Illinois at Chicago \\ 
851 South Morgan Street\\
Chicago, IL, 60607-7045}

\maketitle

\thispagestyle{empty}

\begin{abstract}
We generalize virtual knot theory to multi-virtual knot theory where there are a multiplicity of virtual crossings. Each virtual crossing type can detour over the other virtual crossing types, and over  
 classical or immersed crossings. New invariants of multi-virtual knots and links are introduced and new problems that arise are described. We show how the extensions of the Penrose coloring evaluation for trivalent plane graphs and our generalizations of this to non-planar graphs and arbitrary numbers of colors acts as a motivation for the construction of the multi-virtual theory.
 
 \end{abstract}

\noindent {\bf MSC2020-Mathematics Subject Classification System:} $05C15, 57K10, 57K12, 57K14.$\\

\section{Introduction}
This paper discusses   a generalization of virtual knot theory that we call
 multi-virtual knot theory. Multi-virtual knot theory uses a multiplicity of types of virtual crossings. As we will explain, this multiplicity is motivated by the way it arises first in
 a graph-theoretic setting in relation to generalizing the Penrose evaluation for colorings of planar trivalent graphs to all trivalent graphs, and later by its uses in a virtual knot theory. As a consequence, the paper begins with the graph theory as a basis for our constructions, and then proceeds to the topology of multi-virtual knots and links. The reader interested in seeing our generalizations of the original Penrose evaluation, can begin this paper at the beginning and see the graph theory.
 A reader primarily interested in multi-virtual knots and links can begin reading in section 4 with references to the earlier part of the paper.\\
 
 We begin Section 2 with a review of the formulation, due to Peter Guthrie Tait, of the the four color theorem for plane maps in terms of edge colorings of trivalent graphs. This reformulation asks to {\it properly color} a trivalent graph with three colors so that three distinct colors appear at every node of the graph. A graph is said to be $1$-connected if no removal of a single edge will disconnect it.
Tait's reformulation of the four color theorem states that every planar 1-connected trivalent graph can be properly colored. In \cite{P} Roger Penrose gave a recursive formula that calculates the number of proper colorings for planar trivalent graphs. His evaluation for a graph $G$ is denoted by the bracket $[G]$ and satisfies the formulas below.
$$[ \YGlyph ] = [ \VDiag ] -  [ \CCircleDiag ],$$
$$[ O \, G] = 3 [G],\, [O] = 3,$$
where $O$ denotes a graphical loop (possibly with self intersections) disjoint from the rest of the planar graph. See Sections 2,3 and 2.4 for detailed discussion of the Penrose evaluation and its 
generalizations. Here we indicate a circled virtual crossing in the Penrose expansion.  Sometimes the self-intersection  is written without the circle, but we explain (two paragraphs below) why we have to indicate it. One can prove that $[G]$ is the number of proper three-colorings of a planar graph $G.$ \\

Our approach to the Penrose bracket and its generalizations is based on the concept of a {\it formation} due to G. Spencer-Brown \cite{SB}. See Figure~\ref{Form1} for an illustration of a formation. In the initial case of three colors, we choose to illustrate directly the colors red (r) and blue (b) and to indicate the third color purple (p) as a superposition of red and blue. A formation is a planar drawing using red and blue Jordan curves so that the red curves are disjoint from one another, and the blue curves are disjoint from one another, but red and blue curves can share a finite set of disjoint intervals. The resulting drawing can be interpreted as a trivalent graph with a proper three coloring. The edges colored purple are in correspondence with the intervals in the drawing shared by red and. blue. An examination of the figure should clarify this statement. In the figure the corresponding trivalent graph is shown and the edges of the graph that correspond to shared red and blue are labeled with a slash. These purple edges in the coloring correspond to a perfect matching in the trivalent graph.  This illustrates how a proper coloring gives rise to an even perfect matching (See Section 2). Conversely, if one has a properly colored trivalent graph, then there is a formation that corresponds to the coloring. The red and blue curves are regarded as crossing one another when one curve starts on one side of the other, shares an interval, and continues on the other side of the other curve. \\

A particularly interesting property of the Penrose bracket is that it does not count correctly the number of proper colorings for non-planar graphs. By this, we mean that if $G$ is a non-planar graph presented as an immersion in the plane, and one applies the Penrose algorithm to this immersed graph, evaluating all immersed loops in the recursion by $3$ as indicated above, then the result can be an incorrect count of the colorings of the abstract graph. See Figure~\ref{penroseK33} for an example using the graph $K_{3,3}.$ In \cite{Kauffman} Kauffman solves the problem of extending the Penrose bracket to immersed graphs by distinguishing each initial immersion crossing (in this paper we use a box form at the immersion crossing) and adding a special
evaluation expansion at such crossings in the following form.
$$[ \CBoxDiag ] = 2[ \CDotDiag ] -  [ \CCircleDiag ].$$
The dotted crossing is regarded as a coloring vertex that demands that {\it all four lines incident to this vertex have the same color.}  When a boxed or circled crossing has colors assigned to its arcs it can be regarded as two straight arcs that cross
one another at the circled or boxed intersection. Each straight arc has one color. The two colors may be the same or different. These interpretations are important for understanding the relationships with colorings, but are not used in the formal expansions of the formulas themselves, except at the very end.  The straight arcs of the circled or boxed crossing can be colored with 
two colors that are different or the same. The technical reasons why this method adjusts the Penrose bracket count are reviewed in Section 2. Here it is important to point out that this second recursion is defined on the {\it states} $S$ of the Penrose bracket expansion. Each state consists in a collection of planar Jordan curves intersecting in round and square virtual crossings.
The above evaluation is equivalent to a graphical contraction deletion algorithm for a graph $G(S)$ where each state loop corresponds to  a node of $G(S)$ and two nodes are connected by an edge only when their 
loops touch at a boxed virtual crossing. Round virtual crossings do not affect the graph $G(S).$ Then the above relation becomes the contraction deletion relation shown below.
$$M_{\Wigg} = 2 M_{\Contract} - M_{\Delete}.$$ 
Here $\Contract$ denotes the contraction of an edge, which symbolically can be taken as an instruction that the ends must receive the same color. The expansion formula for the boxed crossing corresponds to this expansion formula, which in turn is a 
{\it dichromatic polynomial evaluation} $(-1)^{b(S)}Z_{G(S)}(\delta, v = -2)$ where $b(S)$ denotes the number of box virtual crossings in the state $S$, which is equal to the number of edges
in $G(S).$  The dichromatic polynomial (a version of the Tutte polynomial) is  defined by the 
recursion
$$Z_{G} = Z_{G'}  + v Z_{G''},$$
$$Z_{G \cup \bullet} = \delta  Z_{G}.$$
Here $G'$ is the result of deleting an edge from $G,$ and $G''$ is the result of contracting that same edge. For the proper three coloring, $\delta = 3.$ Later we will use this recursion for $v=2$ and arbitrary $\delta.$ The $\bullet$ denotes a node disjoint from the rest of the graph. 
This means that we evaluate a dichromatic polynomial \cite{LKStat,K,LKNew} for each state, and that the generalized Penrose evaluation is a sum with coefficients over these dichromatic evaluations.  This is a useful direction of reformulation and will be followed up in subsequent papers. We give further information about this dichromatic formulation and its relationships with chromatic polynomials at the end of Section 3. See Figures~\ref{EFF11},~\ref{Taut},~\ref{Taut1}.
See \cite{LKTutte,VCP,LK1,LK2,LK3,LKStat,Kauffman} for other points of view of this author about coloring graphs.\\

Contraction and deletion give rise directly to evaluation formulas. Let $\Wigg$ denote an edge of the graph associated with a boxed virtual crossing, as explained above. 
In Figure~\ref{Taut} we go back to the original definition of the box virtual crossing tensor. It is positive when the one color goes through the box, and it is negative when distinct colors go through the box. This means that for the corresponding edge in a graph $G = G(S)$ we have the evaluation formula
$$M_{\Wigg} = M_{\Contract} - M_{\Edge}$$
since a contracted edge corresponds to same color and a straight edge can stand for its end nodes having distinct colors.
For such a straight edge we also have the tautological color expansion for the chromatic polynomial
$$M_{\Edge} = M_{\Delete} - M_{\Contract}$$
since a deleted edge allows the end nodes to be the same or different, and this last expansion is the graphical counterpart of the logical statement that ``Different  $=$ All $-$ Same".
It was this tautology that led Hassler Whitney \cite{Whitney} to his recursive definition of the chromatic polynomial via deletion and contraction.
Combining these two identities we have
$$M_{\Wigg} = M_{\Contract} - (M_{\Delete} - M_{\Contract}) = 2 M_{\Contract} - M_{\Delete},$$
which is exactly our evaluation formula as explained above.\\

Now go back to the basic expansion corresponding to the box tensor $$M_{\Wigg} = M_{\Contract} - M_{\Edge}.$$ Given an arbitrary graph $G$ with wiggly edges (corresponding to box virtuals) repeated application of this identity will expand to a sum of copies of the graph with each edge either marked as a contraction, or left as a straight edge, multiplied by $(-1)^n$ where
$n$ is the number of straight edges in the marking. See Figure~\ref{Taut1} for an example. Each such marked graph will then be evaluated by the chromatic polynomial expansion
$$C_{\Edge} = C_{\Delete} - C_{\Contract}$$ applied to the graph contracted by the edges labeled as contractions. This description immediately reformulates to a subgraph $H$ of those edges
labeled straight, and the complementary subgraph $H^c$ of contraction edges. Let $G/H^c$ denote the result of contracting all the contraction edges in $G.$ we then have from this discussion that $$M_{G} = \sum_{H} (-1)^{|H|} C(G/H^c)$$  where $|H|$ denotes the number of edges in $H.$ Thus our special evaluation of box virtual states can be seen as an alternating sum of chromatic polynomials. Furthermore, we have (by comparison with the above) proved the formula 
$$M_{G}= (-1)^{|G|}Z_{G}(\delta, v = -2) = \sum_{H} (-1)^{|H|} C(G/H^c)(\delta),$$
where $G$ is any finite graph, $|G|$ is the number of edges of $G$ and $H$ runs over all sub-graphs of $G$ according to the above discussion.
For practical calculation of $M_{G}$ it may be best to use the dichromatic polynomial, but it is certainly of interest that this particular dichromatic evaluation can be expressed as an alternating sum of chromatic polynomials.\\

Returning to the adjustment of the Penrose evaluation, the key is to see that the Penrose definition is a summation over the possible colorings of the graph with products of $\pm i$ ($i^2=-1$)
at each node of the graph. We use a structure of a coloring called a {\it formation} (above and Section 2.1) to see that this product of square roots of negative unity is equal to the parity of the number of crossings of distinctly colored loops in the formation. Since these loops are Jordan curves in the plane, this parity is even when the graph is embedded in the plane. When there are extra immersion crossings
one needs to use a negative sign whenever the colors of crossing loops at an immersion point are distinct. This is the source of the contraction deletion relation discussed above. See Sections 2 and 3 for the details. The Penrose method attaches an epsilon tensor multiplied by $i$ to each node in the graph so that the Penrose bracket is the tensor contraction of this assignment to the graph. Our adjustment adds an extra tensor to the situation at each immersion crossing for a non-planar graph.\\

We then generalize from $\delta = 3$ colors to an arbitrary number $ \delta = n$ of colors by restricting the generalized bracket to trivalent graphs with a given perfect matching structure, denoted $(G,PM).$ A perfect matching ($PM$) for a graph $G$ is a choice of a collection of disjoint edges so that every node in the graph occurs on one of these edges. The generalized Penrose bracket is then well-defined for arbitrary $\delta = [O].$ We denote a perfect matching edge by an edge with a slash. The fully generalized Penrose bracket then takes the form below.
$$[\YMGlyph] = A[ \VDiag ]  + B[ \CCircleDiag ],$$
$$[ \CBoxDiag ] = 2[ \CDotDiag ] -  [ \CCircleDiag ].$$
$$[ O \, G] = \delta [G],\, [O] = \delta.$$
This is a perfect matching polynomial in three variables. With $A=1, B = -1, \delta = n$ this bracket counts the number of {\it proper $n$-colorings} of a perfect matching $(G,PM).$
By a proper $n$-coloring we mean a coloring of the non-PM edges of $G$, from a collection of $n$ distinct colors, so that two distinct colors occur at each node, and the {\it same} two distinct colors occur at the ends of any given perfect matching edge. This notion of $n$-coloring restricts to the classical notion of $3$-coloring with which our discussion began (by taking the remaining color of the three colors as a color for any perfect matching edge that meets two colors). In the generalization, no colors are assigned to perfect matching edges. In Section 3 we give a self-contained proof that the perfect matching polynomial counts $n$-colorings. In that section we discuss examples of such higher order colorings. One key example is the Petersen graph, which is seen to be uncolorable for all values of $n.$ We have used the perfect matching polynomial in relation to virtual knot theory and graph homology in \cite{BKR,BKM}. See also \cite{BMcC,BM-Vertex} for the remarkable work of Scott Baldridge and Ben McCarty. In this paper we give new relations to a wider version of virtual knot theory, and we expect further relations with link homology and graph homology in subsequent work.\\

A pivot in this discussion of coloring occurs as we point out (See Figure~\ref{pcode}) that by replacing the slashed perfect matching edge $\YMGlyph$ with a knot diagrammatic crossing $\CrossStarGlyph$ we can translate
the $n$-coloring count for perfect matching graphs to a corresponding coloring count for link diagrams (with boxed virtual crossings in the case of translation from non-planar graphs). One then has the formulas below for 
finding the number of $n$-colorings of link diagrams.
$$PK( \CrossStarGlyph ) = PK(\VDiag ) - PK( \CCircleDiag  ),$$
$$PK(O) = n,$$
$$PK(\CBoxDiag) = 2 PK(\CDotDiag) - PK(\CCircleDiag).$$
In this color count a link diagram is properly $n$-colored if the edges of the diagram (viewed as a $4$-valent graph) are colored from $n$ colors so that at each crossing, pairs of distinct colors appear in correspondence to our rule for perfect matching edges. See Figure~\ref{pcode}. This means that we have a knot diagrammatic interpretation of the coloring problem and a polynomial expansion for it in the category of link diagrams. The polynomial $PK$ is not a topological link invariant but it is a close relative of the bracket polynomial invariants that we shall now define for multiple virtual knot and link theory.\\

Section 4 begins the definition and exploration of multiple virtual knot theory in this paper. We use the graph theory described above, with its use of two distinct virtual crossings to motivate the definition and construction of multiple virtual knot theory where there can be an arbitrary number of virtual crossings. These virtual crossings all detour over one another, while the classical crossings do not detour over the virtual crossings. This paper can be read as an introduction to virtual knot theory, and the reader may wish to read  earlier introductions and research by the author \cite{VKT,DVK,vkt,rotvkt}. Multiple virtual knot theory is presented here in terms of knot and link diagrams that have
classical (and sometimes immersed) crossings and virtual crossings. The classical crossings interact with one another via the Reidemeister moves and the virtual crossings interact via detour moves. In a detour move a consecutive sequence of virtual crossings of the same type can be excised from the diagram and the endpoints of this excision can be reconnected by any diagrammatic arc that intersects the rest of the diagram in a sequence of virtual crossings of this same type. In sections 4.1, 4.2, 4.3 and 4.4 we discuss this definition of multiple virtual knot theory and we give a sketch of the topological interpretation of this definition in the case of single virtual knot theory.
Single virtual knot theory can be interpreted in terms of knots and links embedded in thickened surfaces taken up to knot theoretic equivalence in these surfaces and up to one-handle stabilization. In Sections 4.3 and 4.4 we describe a similar but more restricted topological interpretation for multiple virtual knot theory that is obtained by adding labeled handles at the site of each virtual crossing. Handles with the same label can merge, and handle stabilization arises in relation to these operations.\\

Section 4.5 defines a generalized bracket polynomial for multiple virtual knot theory. The basic bracket expansion is given by the formula
$$\langle \CrossGlyph \rangle = A \langle \HSmoothGlyph \rangle + B \langle \VSmoothGlyph \rangle.$$
This expansion can also be seen as a knot theoretic version of deletion contraction by way of the well-known checkerboard graph and medial graph constructions in 
knot theory and graph theory. See \cite{LKNew, LKStat,LKTutte,KauffJaeger,Kauff:KP}. By itself, the bracket polynomial is a breakthtough between low dimensional topology and the graph theoretical structures of the chromatic, dichromatic and Tutte polynomials. Here we extend the bracket to an invariant of multiple virtual knots. The most general version uses $B=A^{-1}$ and $\delta = -A^{2} - A^{-2}$ to give an invariant of regular isotopy of multiple virtual knots and links. In this formulation we have a state summation formula  in the form $$\langle K \rangle = \sum_{S} \langle K | S \rangle \langle S \rangle$$ where $ \langle K | S \rangle$ denotes a product of $A's$ amd $B's$ corresponding to a choice of smoothing of each crossing in the diagram $K.$ The state $S$ indicates such a choice of smoothings. $S$ consists in a collection of curves in the plane intersecting transversely in virtual crossings of different types. $\langle S \rangle$ denotes the equivalence class of $S$ under detour moves on these virtual crossings. In general these equivalence classes can be difficult to determine, but in many special cases the generalized invariant can yield much information.\\

We study a special case of this new invariant by combining it with the graphical deletion contraction algorithms explained in this introduction that are applicable to mutiple virtual diagrams. The main specialization of the bracket is called the {\it chromatic bracket} and it is specifically defined for multiple virtual systems with two distinct virtual crossings, just as in the graph theoretic Penrose evaluations described above. The equations then are as follows.
$$\langle \CrossGlyph \rangle = A \langle \HSmoothGlyph \rangle + B \langle \VSmoothGlyph \rangle$$
$$\langle O \rangle = \delta$$
$$\langle \IMGlyph \rangle  = 2 \langle \CDotDiag \rangle  - \langle \VirtualGlyph \rangle $$
With these conventions the chromatic bracket is a Laurent polynomial associated with the diagram $K.$ When $B=A^{-1}$ and $\delta = -A^{2} - A^{-2}$  we obtain an invariant of regular isotopy for vitual links of multiplicity two.
When $A=1$ and $B = -1$ we obtain a coloring polynomial written in the language of link diagrams as we have explained above. Section 4.5 gives many topological examples for these polynomials, including an example
(See Figure~\ref{EFF28}) that cannot be detected by the generalized bracket and is related to links that are not detectable by the Jones polynomial. The section finishes with an example of a knot $K$ where the non-triviality of its generalized bracket depends upon the non-triviality of a state diagram that is not detected by the chromatic calculation (see Figure~\ref{EFF30} and Figure~\ref{EFF31}). We develop quandle techniques for this detection in the next section of the paper.\\

Section 5 introduces generalizations of the quandle and the biquandle to multi-virtual theory and discusses a number of examples, including a generalized Alexander polynomial.\\

Section 6 discusses multiple virtual link cobordism and in particular shows that the non-trivial knot $K$ from the end of Section 5 is a multi-virtual slice knot.\\

Section 7 discusses the parity bracket polynomial in the multiple virtual context.\\

Section 8 discusses multiple virtual generalizations of the arrow polynomial. The arrow polynomial in standard virtual theory has many irreducible state loops. In the multiple generalization the classification of the arrow state
configurations presents many new open problems.\\

Section 9 discusses welded knot theory in the multiple context. Again there are many questions and new research problems.\\

Section 10 discusses rotational multiple virtuals. Here the bracket polynomial in its generalized form can be used and we give examples of distinguishing chirality and detecting differences for planar diagrams.\\

Section 11 discusses the definition of the multiple virtual braid group and states a Markov Theorem that will be proved and explored elsewhere.\\

Section 12 is an Epilogue about the general structure of knot theoretic generalizations.\\

\section {The Penrose Formula and Its Generalizations}
Before discussing the Penrose formula we give background on colorings of cubic graphs. The background follows closely our previous paper \cite{Kauffman} in order to provide the groundwork for the significant 
generalizations, graph theoretic and topological, that occur in the present paper. In section 3 we begin new graph theoretic results related to the Penrose evaluation.\\

\subsection{\bf Cubic Graphs and Formations}
A {\em cubic graph} is a graph in which every vertex is locally incident to three edges. The vertex either belongs to three
distinct edges, or there are two edges at the vertex with one of them a loop. A
{\em coloring} (proper coloring) of a cubic graph $G$ is an assignment of three 
labels $r$ (red), $b$ (blue), and $p$ (purple) to the edges of the graph so that
three distinct labels are incident to each vertex of the graph. This means that 
three distinct edges belong to each vertex and it is possible to
label the graph so that three distinct colors occur at each vertex. A
graph with a loop is not colorable. \\

The simplest uncolorable cubic graph is shown in Figure~\ref{peter}.  We refer to this graph as the {\em dumbell}.  The dumbell is a
planar graph. The figure shows a more complex dumbell and the Petersen graph, a non-planar uncolorable graph.\\

An edge in a connected plane graph is an {\em isthmus} if the deletion of that
edge results in a disconnected graph. A connected plane
cubic graph without isthmus is loop-free. \\

Petre Guthrie Tait (see \cite{Kempe,Heawood,Petersen}) reformulated the four-color conjecture (here called the {\em Map Theorem}) for plane maps to a statement about 
colorability of plane cubic graphs. In this form the theorem says \\

\noindent {\bf Map Theorem for Cubic Graphs.}  A connected plane cubic graph without
isthmus is properly edge-colorable with three colors. \\

 Let $G$ be a cubic graph and let $C(G)$ be a coloring of $G.$ Using the colors
$r$, $b$ and $p$ we write purple as a formal product of red and blue: $$p =
rb.$$

\noindent There are single colored paths on the coloring $C(G)$ in the
colors red and blue. Each red or blue path eventually returns to its starting
point, making a circuit in that color. The red circuits are disjoint from one
another. The blue circuits are disjoint from one another. Red circuits
and blue circuits can meet along edges in $G$ that are colored purple
($p=rb$).   For a plane graph $G$, a meeting of two circuits can take
the form of one circuit crossing the other in the plane, or one circuit can share
an edge with another circuit,  leaving on the same side of that other circuit.
These two
planar configurations are called a {\em cross} and a {\em bounce} respectively. \\

\begin{figure}[htb]
     \begin{center}
     \begin{tabular}{c}
     \includegraphics[width=5cm]{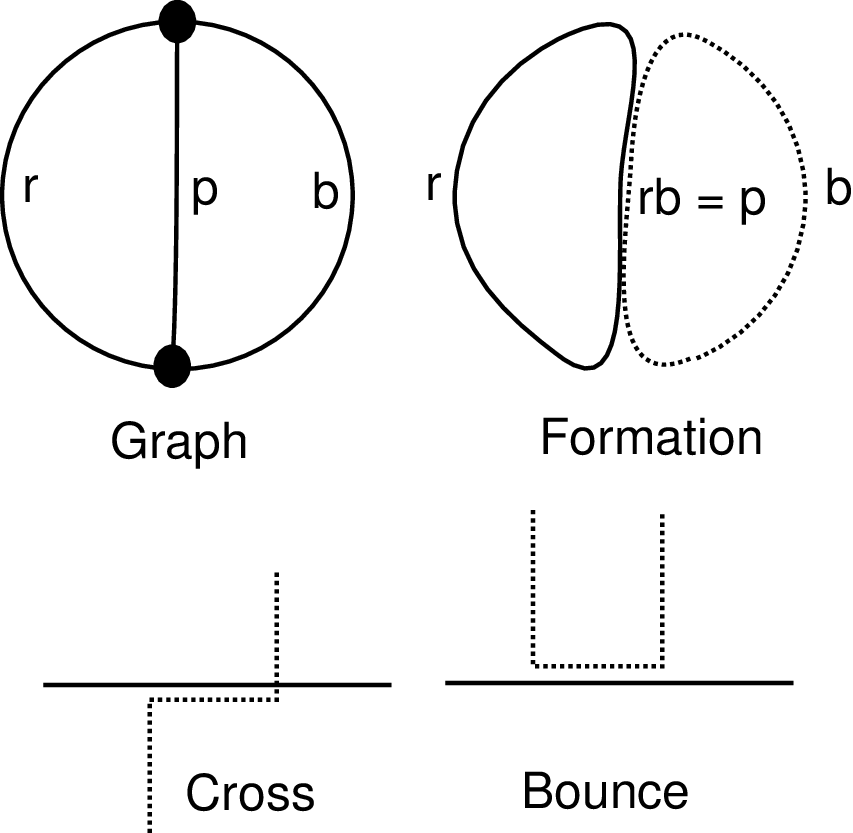}
     \end{tabular}
     \caption{\bf Coloring and Formation.}
     \label{theta}
\end{center}
\end{figure}

\noindent {\bf Definition.} A {\em formation} \cite{SB} is a (finite) collection of 
simple closed curves in the plane.  Each curve is colored either red or blue such that the
red curves are disjoint from one another, the blue curves are disjoint from one
another and red and blue curves can meet in a finite number of segments (as
described above for graphical circuits).\\

\noindent Associated with a formation $F$ there is a cubic graph
$G(F)$,  obtained by identifying the shared segments in the formation as edges in
the graph. The endpoints of these segments are the nodes of $G(F).$ The
unshared segments of each simple closed curve constitute the other edges of
$G(F).$  A formation $F$ is said to be a formation for a cubic graph $G$ if $G =
G(F).$ Then we say that $F$ {\em formates} $G.$ \\

\noindent A {\em plane formation} is a formation such that each simple closed
curve in the formation is a Jordan curve in the plane. For a plane formation,
each shared segment between two curves of different colors is either a bounce or
a crossing (see above), that condition being determined by the embedding of the
formation in the plane. \\

The notion of a formation is abstracted from the circuit decomposition of a
colored cubic graph.  We have the proposition \cite{Kauffman}: \\

\noindent {\bf Proposition.}  Let $G$ be a cubic graph.  Let $Col(G)$ be the set of
colorings of $G$.  Then $Col(G)$ is in one-to-one correspondence with the set of
formations for $G.$ 
\bigbreak

\noindent The Map Theorem is equivalent to the
\bigbreak

\noindent {\bf Formation Theorem.}  Every connected plane cubic graph with no
isthmus has a formation. \\

\begin{figure}[htb]
     \begin{center}
     \begin{tabular}{c}
     \includegraphics[width=5cm]{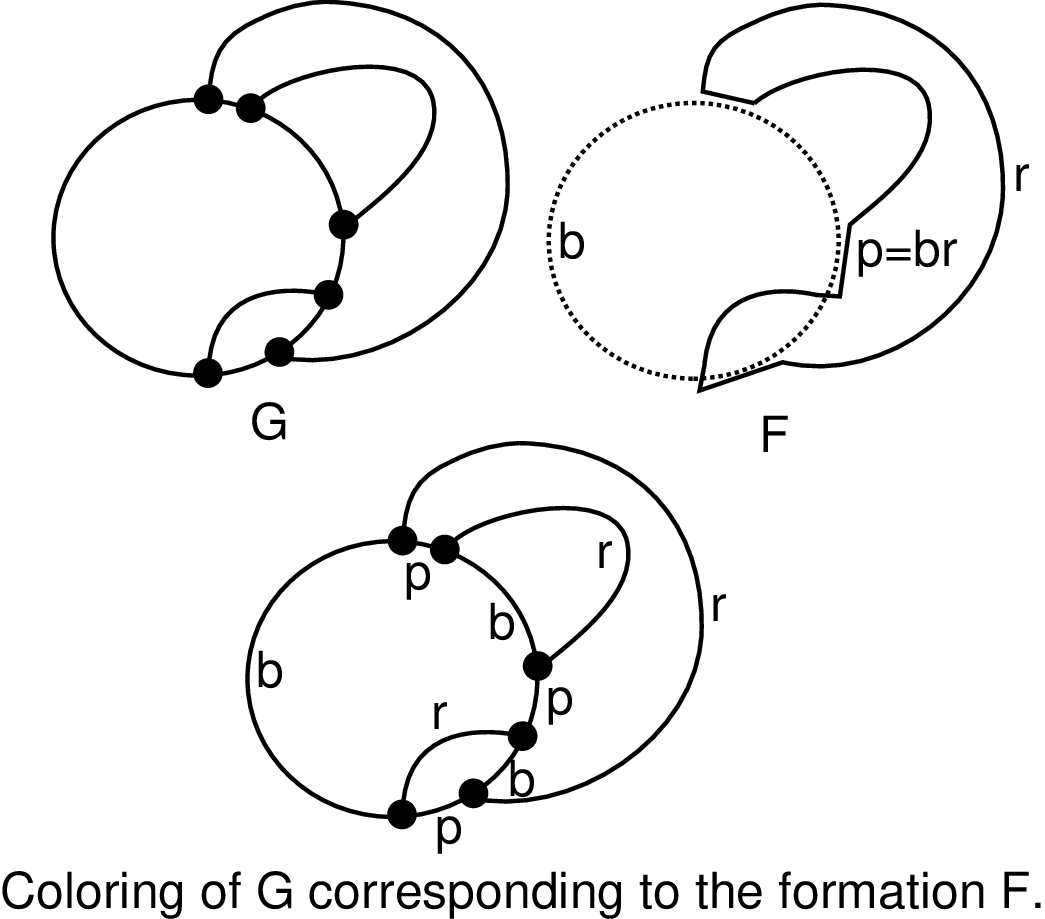}
     \end{tabular}
     \caption{\bf Second Example of Coloring and Formation.}
     \label{colform}
\end{center}
\end{figure}

This version of the Map Theorem is due to G. Spencer-Brown \cite{SB}. Just as one can enumerate graphs,
one can enumerate formations.  In particular, plane formations are generated by
drawing systems of Jordan curves in the plane that share segments according to
the rules explained above. This gives evidence for the Map
Theorem, since one can enumerate formations and observe that all plane cubic
graphs are occurring in the course of the enumeration.  See Figures~\ref{theta} and \ref{colform}. \\

\noindent {\bf Remark.}    In depicting formations, we have endeavored to keep the
shared segments separated for diagrammatic clarity. The separated
segments are amalgamated in the graph that corresponds to the formation.\\

\subsection{Cubic Graphs and Perfect Matchings}
A cubic map $G$ is said to be {\em properly colored} with $3$ colors if the edges of $G$
are colored so that the colors incident to any node of $G$ are distinct.\\

We give another equivalent version of the Four Color Theorem (FCT). 
Call a disjoint collection $E$ of edges of $G$ that includes all the vertices of $G$ a {\em perfect matching}
of $G.$ Then $C(E,G) = G - Interior(E)$ is a collection of {\em cycles}. $E$ is called a {\em even} perfect matching of $G$ if every cycle in $C(E)$ has
an even number of edges.\\

\noindent The following statement is equivalent to the Four Color Theorem: \\

\noindent {\bf Theorem.} Let $G$ be a plane cubic graph with
no isthmus. There there exists an even  perfect matching of $G.$
\bigbreak

\noindent {\bf Proof.} Take a proper coloring of $G.$ Choose a color and designate all the edges of $G$ that receive this color. This designation is an even perfect matching. For the rest of the proof, see \cite{Kauffman}.
 $\hfill\Box$ \\

\begin{figure}[htb]
     \begin{center}
     \begin{tabular}{c}
     \includegraphics[width=5cm]{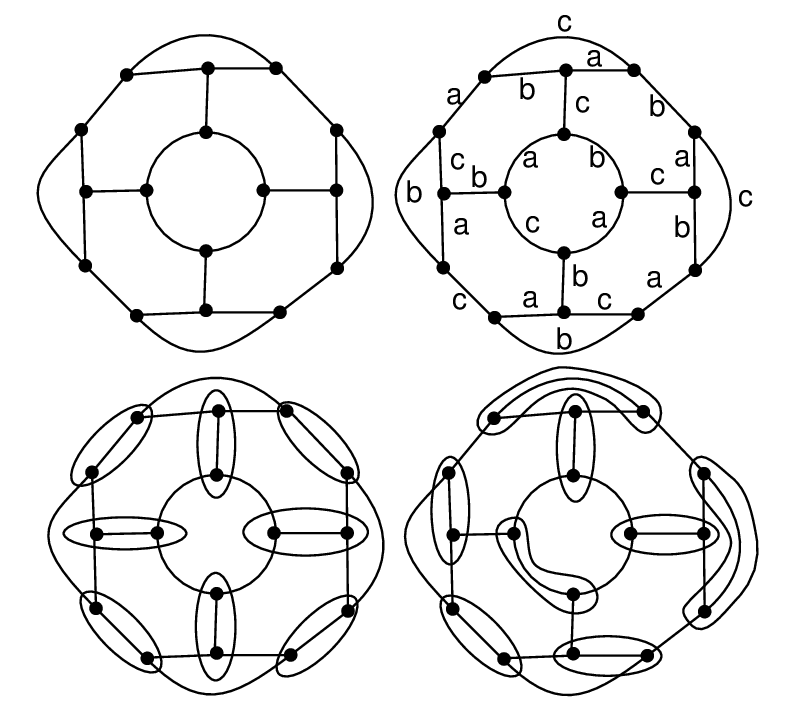}
     \end{tabular}
     \caption{\bf Perfect Matchings of a Cubic Plane Graph}
     \label{perfect}
\end{center}
\end{figure}

\noindent {\bf Remark.} See Figure~\ref{perfect} for an illustration of two perfect matchings of a graph $G.$  One perfect matching is even, and the  corresponding coloring is shown with each pefect matching edge marked with the letter $c$ and the remaining edges marked with $a$ and $b.$
The Theorem shows that one could
imagine dividing the proving of the FCT into two steps: First prove that every cubic plane isthmus-free graph has a perfect
matching. Then prove that it has an  even  perfect matching. In fact, the existence of a perfect matching is available
\cite{Petersen}, while the existence of an even perfect matching is apparently quite hard.\\

\noindent {\bf Proposition.} Every cubic graph with no isthmus has a  perfect matching.
\bigbreak

\noindent {\bf Proof.}  See  \cite{Petersen}, Chapter $4.$ $\hfill\Box$
\bigbreak

\begin{figure}[htb]
     \begin{center}
     \begin{tabular}{c}
     \includegraphics[width=5cm]{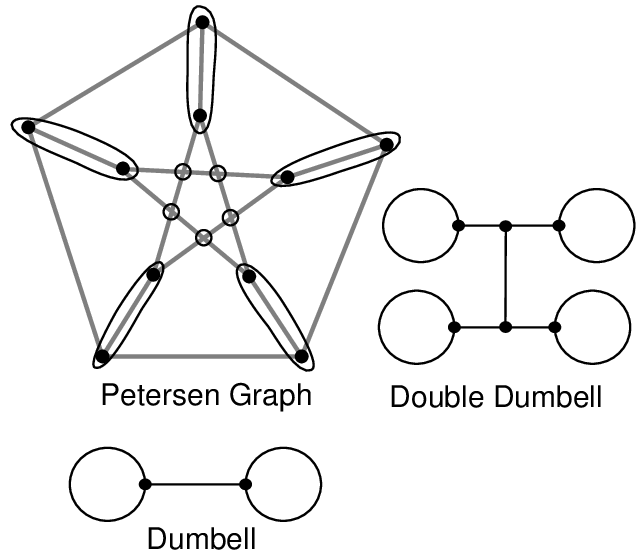}
     \end{tabular}
     \caption{\bf Petersen and Dumbbells}
     \label{peter}
\end{center}
\end{figure}

\noindent {\bf Remark.}  Two uncolorable graphs are indicated in Figure~\ref{peter}. These graphs have perfect matchings, but
no even perfect matching. The second example in  Figure~\ref{peter} is the ``dumbell graph".  The first
example is the non-planar Petersen Graph.  We show the Petersen with a perfect matching with two
five cycles. No  perfect matching of the Petersen is even. The third  "double dumbell" graph shown in Figure~\ref{peter} has no
perfect matching.
\bigbreak

\subsection{Penrose Evaluation Formula}
Roger Penrose \cite{P} gave a formula for finding the number of proper 
edge 3-colorings of a plane cubic graph $G$ with three distinct colors at each node of the graph.
The Penrose formula satisfies the identities shown below and in Figure~\ref{penrose}.
$$[ \YGlyph ] = [ \VDiag ] -  [ \CDiag ],$$
$$[ O \, G] = 3 [G],$$
$$[O] = 3,$$ where $O$ denotes a Jordan curve disjoint from $G$ in the plane.\\

\begin{figure}[htb]
     \begin{center}
     \begin{tabular}{c}
     \includegraphics[width=5cm]{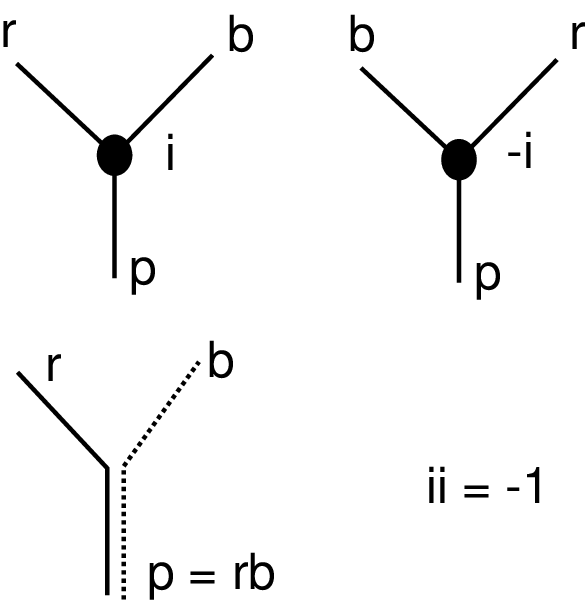}
     \end{tabular}
     \caption{\bf Node as Epsilon Tensor}
     \label{epsilon}
\end{center}
\end{figure}

This formula gives a recursive computation of the number of colorings for planar cubic graphs.
As Penrose explains in his paper \cite{P}, the formula is a consequence of associating to each node in the graph
an ``epsilon" tensor 
$$E_{ijk}=\sqrt{-1}\epsilon_{ijk}$$ 
\noindent as shown in Figure~\ref{epsilon}. One takes the colors from the set $\{1,2,3\}$ and the tensor $\epsilon_{ijk}$ takes value $1$
for $ijk = 123, 231,312$ and $-1$ for $ijk = 132, 321, 213.$  The tensor is $0$ when $ijk$ is not 
a permutation of $123.$  We use $r = 1, b=2, p = 3$ when we want letter designations for the colors.\\

One then evaluates the graph $G$ by taking the sum over all possible color
assignments to its edges, of the products of the $E_{ijk}$ associated with its nodes. Call this 
evaluation $[G].$ \\

Due to the epsilon tensor the Penrose formula satisfies the identities shown above, and so we have the result below. \\

\begin{figure}[htb]
     \begin{center}
     \begin{tabular}{c}
     \includegraphics[width=14cm]{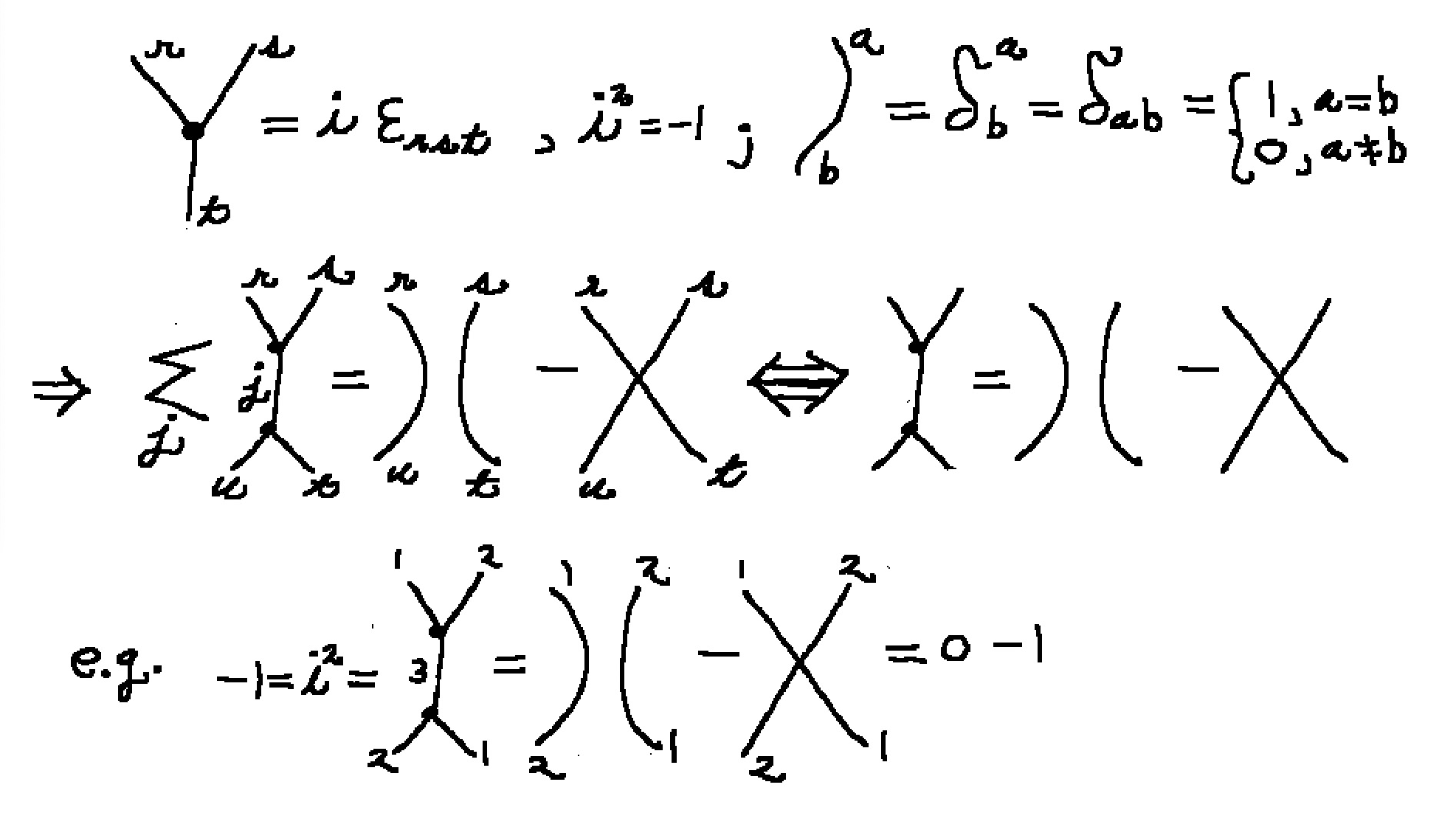}
     \end{tabular}
     \caption{\bf Epsilon Equation}
     \label{epsilonequation}
\end{center}
\end{figure}

\begin{figure}[htb]
     \begin{center}
     \begin{tabular}{c}
     \includegraphics[width=5cm]{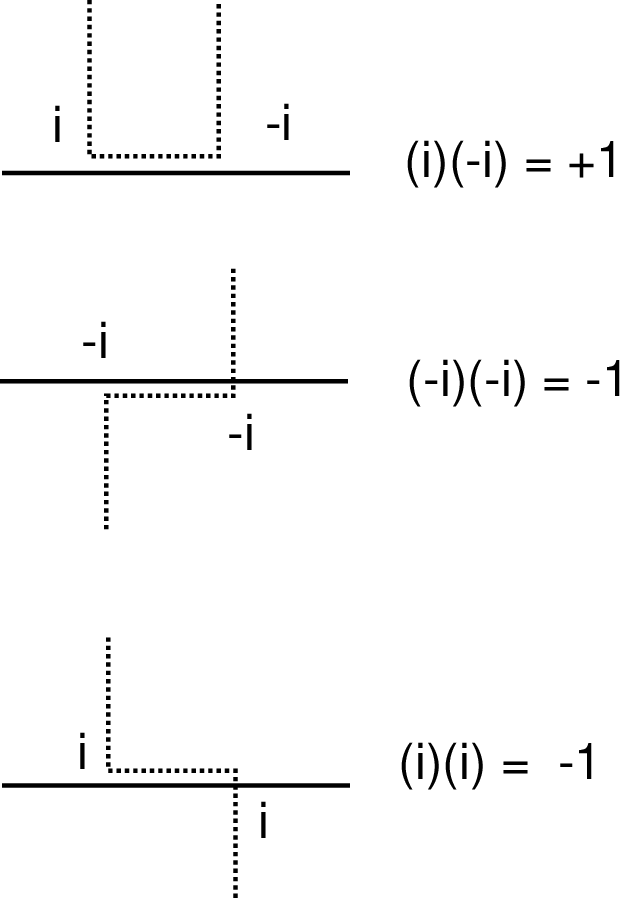}
     \end{tabular}
     \caption{\bf Bounce and Cross Contribute +1 and -1.}
     \label{bouncecross}
\end{center}
\end{figure}

\begin{figure}[htb]
     \begin{center}
     \begin{tabular}{c}
     \includegraphics[width=5cm]{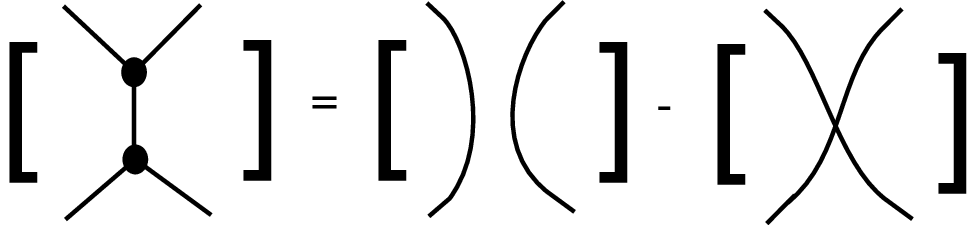}
     \end{tabular}
     \caption{\bf Penrose Formula.}
     \label{penrose}
\end{center}
\end{figure}

\begin{figure}[htb]
     \begin{center}
     \begin{tabular}{c}
     \includegraphics[width=5cm]{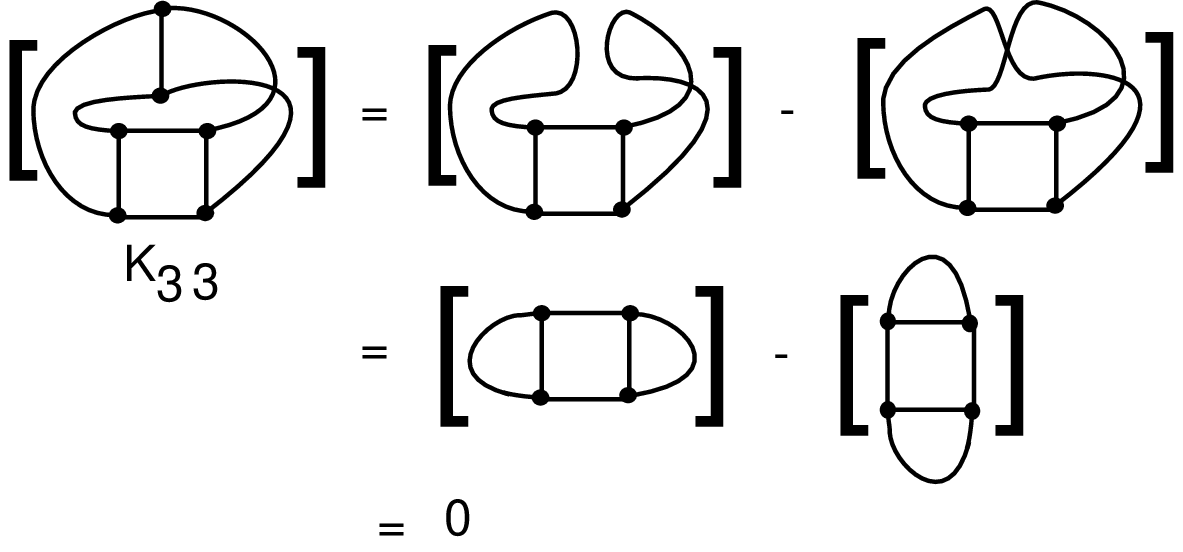}
     \end{tabular}
     \caption{\bf Penrose on K33 is Zero.}
     \label{penroseK33}
\end{center}
\end{figure}

\noindent {\bf Theorem (Penrose).} If $G$ is a planar cubic graph, then $[G]$, defined by the formulas
$$[ \YGlyph ] = [ \VDiag ] -  [ \CDiag ],$$
$$[ O \, G] = 3 [G],$$
$$[O] = 3,$$ as described above,
is equal to the number of
distinct proper colorings of the edges of $G$ with three colors (so that every node sees three colors at its edges).\\

\noindent {\bf Proof.} The formulas in the statement of the Theorem follow from properties of the epsilon tensor.
It is easy to see that the tensor satisfies the formula $$\sum_{j=1}^{3}\epsilon_{rsj}\epsilon_{jtu} = - \delta_{ru}\delta_{st} + \delta_{rt}\delta_{su}$$
and this becomes, after multiplication of each epsilon by $i$, $$\sum_{j=1}^{3}E_{rsj}E_{jtu} = \delta_{ru}\delta_{st} - \delta_{rt}\delta_{su},$$ which translates via the 
diagrammatic of Figure~\ref{epsilon} and Figure~\ref{epsilonequation}  and the identification of a line segment with labels $a$ and $b$ at its ends with a Kronecker delta, $\delta_{ab}$, to $$[ \YGlyph ] = [ \VDiag ] -  [ \CDiag ].$$ The fact that a loop evaluates to $3$ is because we use three colors.\\

It follows from the above description that only proper colorings of $G$ contribute to 
the summation $[G],$ and that each such coloring contributes a product of $\pm \sqrt{-1}$ from the tensor
evaluations at the nodes of the graph. In order to see that $[G]$ is equal to the number of 
colorings for a plane graph, one must see that each such contribution is equal to $+1.$
The proof of this assertion is given in Figure~\ref{bouncecross}  with reference to Figure~\ref{epsilon} where we see that in a formation for a coloring 
each bounce contributes $+1 = -\sqrt{-1}\sqrt{-1}$ while each crossing contributes $-1.$ Since, by the Jordan Curve Theorem, there
are an even number of crossings among the curves in the formation, it follows that the total 
product is equal to $+1$. This completes our proof of the Penrose Theorem. $\hfill\Box$
\bigbreak

The original Penrose formula only counts colors correctly for a plane graph. For example, view Figure~\ref{penroseK33} where we show that the Penrose bracket calculation of the $K_{3,3}$ graph is zero. This graph has $12$ colorings. \\

By using the properties of the epsilon tensor and the properties of formations we showed that each coloring contributed $+1$ to the state summation for a planar graph. This sign of $+1$ depended on our use of the Jordan curve theorem in counting colored Jordan curves that intersected one another an even number of times. When we use a planar diagram with immersion crossings the number of such intersections can be odd.\\

We modify the Penrose bracket so that it will calculate the number of colorings of any graph (given as a diagram in the plane with immersion crossings). These immersion crossings are artifacts of the way the diagram is drawn in the plane, and are necessary for non-planar graphs. Call these the {\em virtual crossings}. See Figure~\ref{penroseK33} for an illustration of an immersion of $K_{33}$ into the plane with one immersion (virtual) crossing.\\

\begin{figure}[htb]
     \begin{center}
     \begin{tabular}{c}
     \includegraphics[width=5cm]{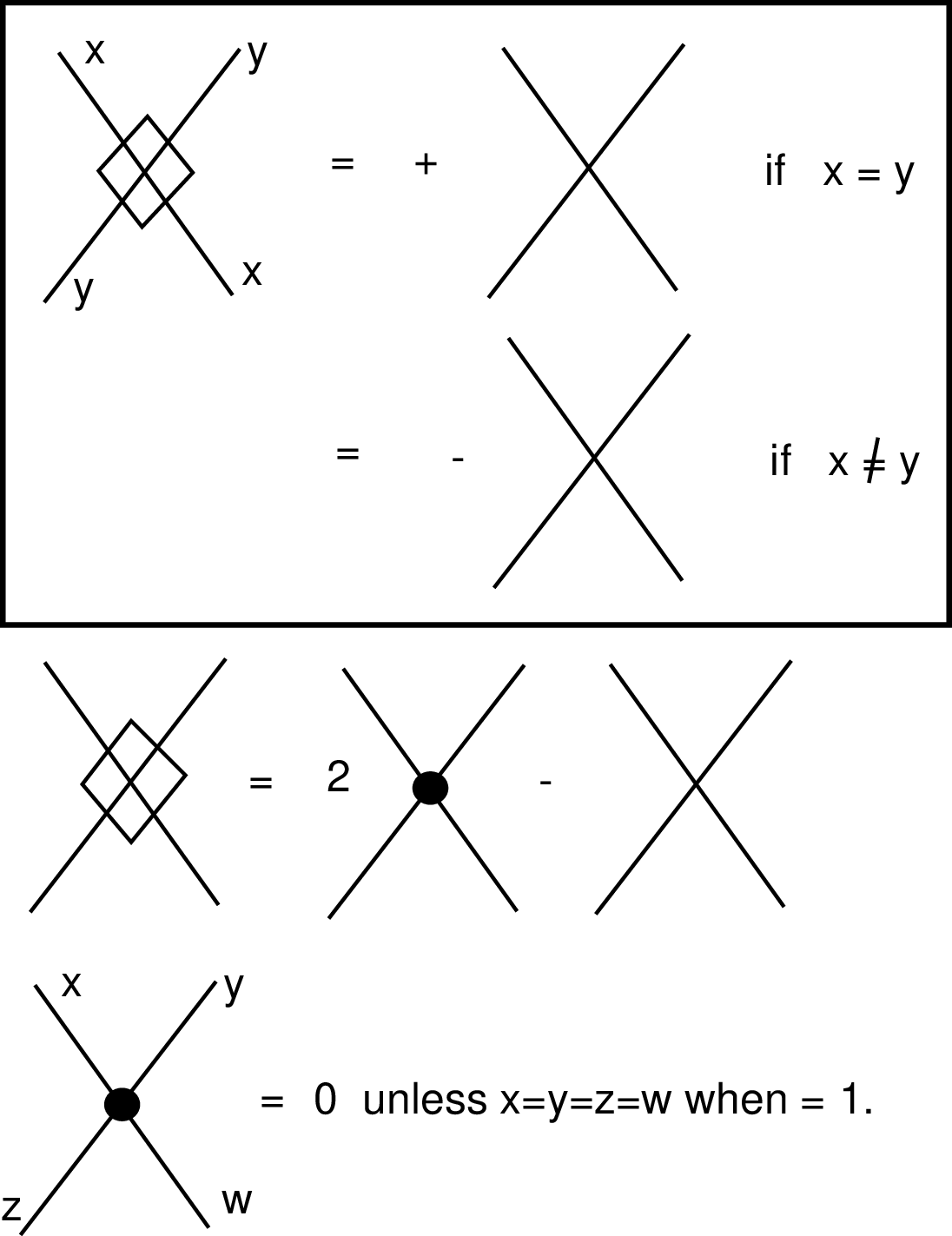}
     \end{tabular}
     \caption{\bf Crossing Tensor For Revised Penrose Bracket.}
     \label{crossingtens}
\end{center}
\end{figure}

\noindent{\bf The New Tensor.}
We add an extra tensor for each immersion crossing in the original diagram. The new tensor is indicated by a box around that crossing as in $\CBoxDiag$ (as we did before with a circle for virtual knot crossings but this is distinct from that method). The boxed crossing is a tensor depending on its four endpoints. It is zero unless the labels at the ends of a given straight segment in the crossed segments are the same. When these ends are the same we can label the crossed segments with two colors $x$ and $y,$ as in Figure~\ref{crossingtens}. Then the value of the tensor is $+1$ if $x=y$ and $-1$ if $x \neq y,$ as shown in the Figure~\ref{crossingtens}. The Penrose bracket is computed just as before and it satisfies the same formulas as before. The new crossings that occur in the expansion formula are standard crossings. {\it Only the initial immersion crossings in the diagram are boxed.} Examine the proof above and note that now every coloring will contribute $+1$ to the state sum. \\

\noindent {\bf Remark.}  We write formally the equation
$$[ \CBoxDiag ] = 2[ \CDotDiag ] -  [ \CDiag ],$$ or equivalently
$$[ \CBoxDiag ] = 2[ \CDotDiag ] -  [ \CCircleDiag ].$$
The dotted crossing is regarded as a coloring vertex that demands that {\it all four lines incident to this vertex have the same color.} The ordinary crossing is taken as usual to indicate that each crossing line has a color independent of the other line. In the second equation we have replaced the ordinary crossing with a circled crossing. This means the same thing, and will be used later to discriminate these two types of virtual crossing as boxed and circled. 
Note that evaluation via this difference formula gives the same results as the tensor in Figure~\ref{crossingtens}, as is easily checked in the two cases when $x$ and $y$ are equal or unequal.
\\

And thus we have the Theorem \cite{Kauffman}:\\

\noindent {\bf Theorem (Penrose Bracket for All Cubic Graphs).} Let $G$ be any cubic graph equipped with an immersion into the plane so that the transverse crossings of interiors of edges of $G$ are boxed as described above. Interpret the 
nodes of $G$ as epsilon tensors and the boxed crossings as sign tensors as described above. Let $[G]$ denote the tensor contraction of the immersed graph $G$ for the three color indices. Then 
$[G]$ is equal to the number of proper three-colorings of the cubic graph $G.$ The new version of the Penrose bracket satisfies the basic formulae
$$[ \YGlyph ] = [ \VDiag ] -  [ \CDiag ],$$
$$[ O \, G] = 3 [G],$$
$$[O] = 3.$$ 
The virtual crossings resulting from the immersion of $G$ are expanded by the formula
$$[ \CBoxDiag ] = 2[ \CDotDiag ] -  [ \CDiag ]$$
as described above so that the dotted crossing is regarded as a coloring vertex that demands that {\it all four lines incident to this vertex have the same color.} \\

\noindent {\bf Proof.} The proof of this result follows from the discussion preceding the statement of the Theorem. 
$\hfill\Box$
 \\

\begin{figure}[htb]
     \begin{center}
     \begin{tabular}{c}
     \includegraphics[width=8cm]{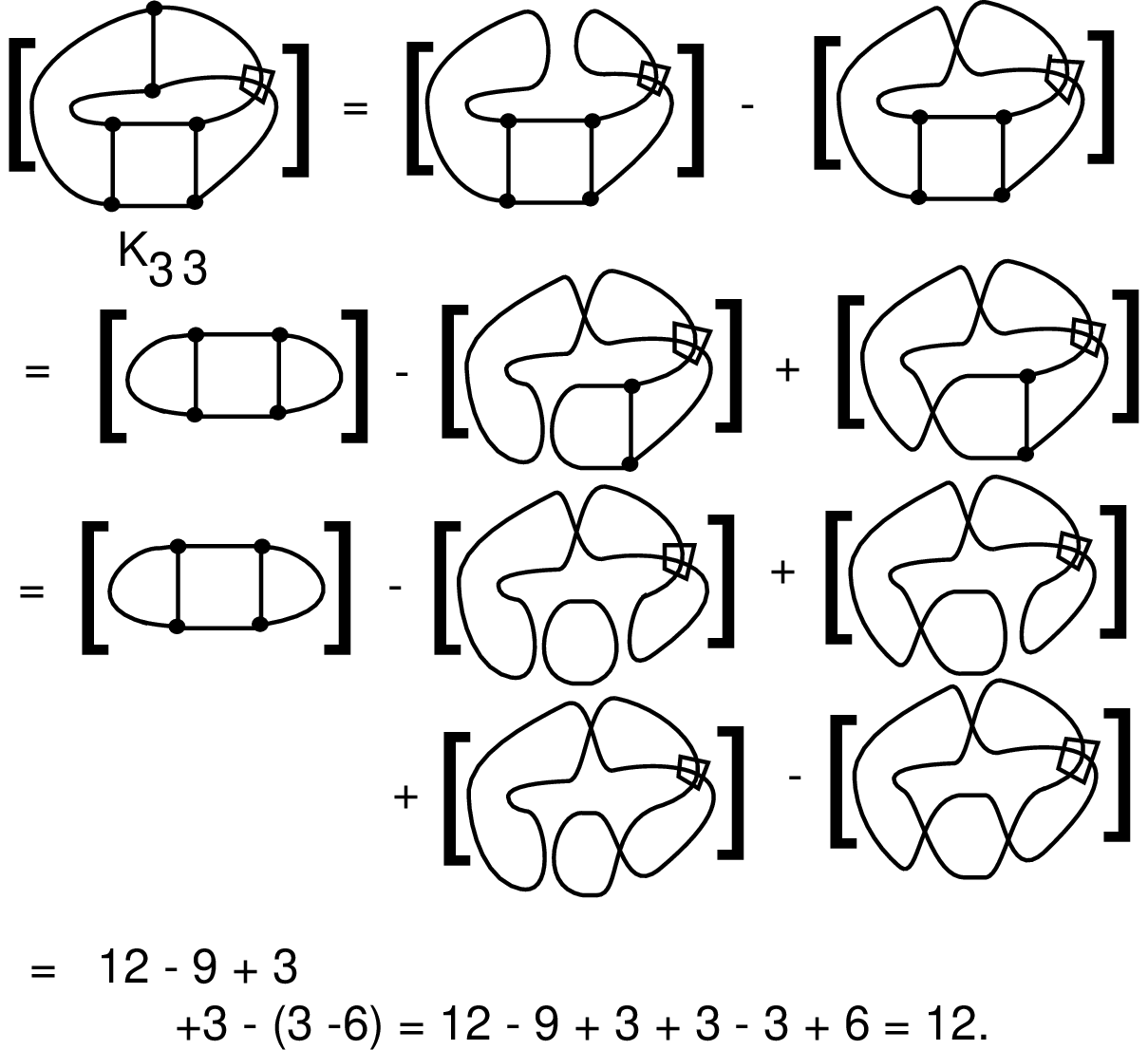}
     \end{tabular}
     \caption{\bf Revised Penrose evaluation on $K_{33}$ Counts Colorings.}
     \label{revisedpenroseK33}
\end{center}
\end{figure}

\noindent {\bf Remark.} See \cite{Kauffman} for detailed examples applying this result.
Examine Figure~\ref{revisedpenroseK33}. In this figure we have an immersion of the $G = K_{3,3}$ with one crossing and this crossing is boxed. We then expand the Penrose bracket until we have collections of circles for the second term of the first expansion. We keep the first term and rewrite it without a boxed crossing because this is a self-crossing and will always have equal colors and hence 
will contribute a $+1$ in the evaluation. In the other part of the expansion to collections of circles, each final diagram can be decided, using the rule for the sign of a boxed crossing. Note that in the very last diagram we have two curves that cross once with a boxed crossing and once with an unboxed crossing. When these two curves are colored the same the circle gives $+1$ and when they are colored unequally, the circle gives $-1.$ Thus these two circles contribute $3 -6 = -3.$ The reader will note that the initial second term adds up to zero and the total for $[G]$ is $12,$ the number of colorings of this graph.\\
\\

 We have an extension of the Penrose bracket that calculates the number of colorings of any cubic graph. Our more abstract point of view can be compared with \cite{EM} where a generalization of the Penrose bracket is made for graphs embedded in surfaces.\\

\subsection{Loops and Logical State Expansions}
\begin{figure}[htb]
     \begin{center}
     \begin{tabular}{c}
     \includegraphics[width=8cm]{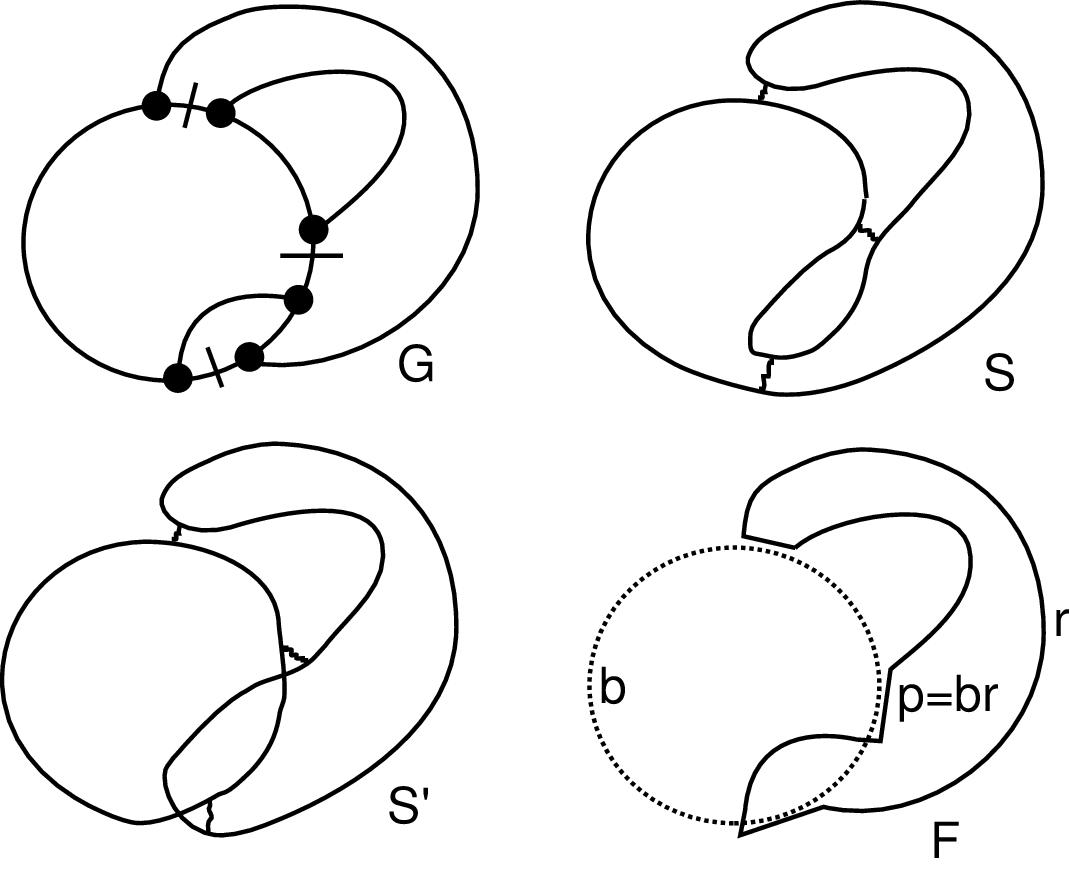}
     \end{tabular}
     \caption{\bf Graph, State, Switch and Formation.}
     \label{state}
\end{center}
\end{figure}

\begin{figure}[htb]
     \begin{center}
     \begin{tabular}{c}
     \includegraphics[width=8cm]{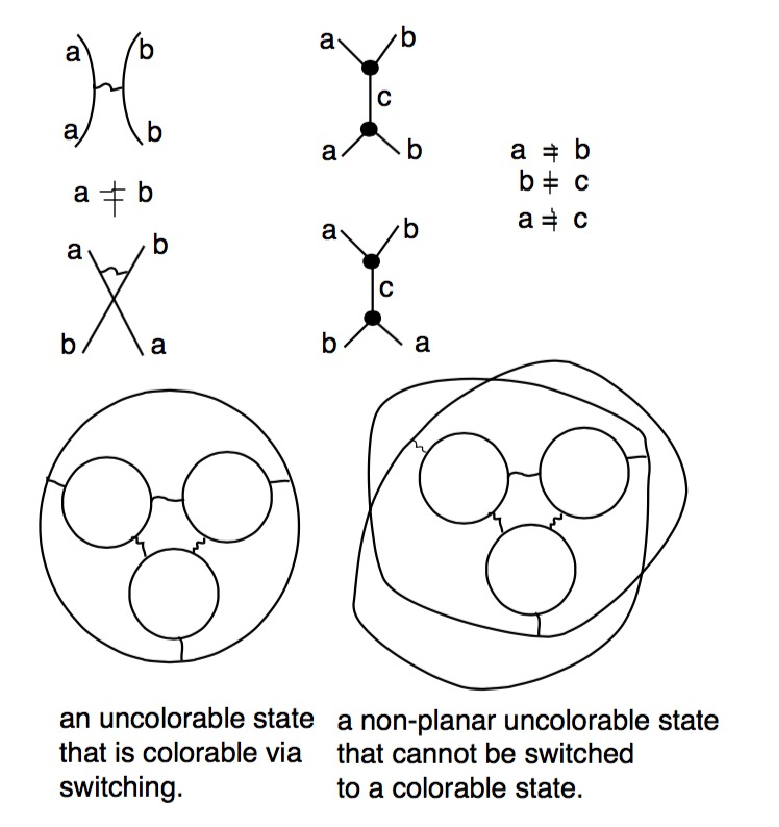}
     \end{tabular}
     \caption{\bf Crossed and Uncrossed Sites and Two Examples of States}
     \label{site}
\end{center}
\end{figure}

In this section we give a reformulation of the map theorem that is very close in spirit to the diagrammatics of low-dimensional knot theory.
Consider a collection of disjoint Jordan curves drawn in the plane. Let there be a collection of {\it sites} between distinct curves as indicated in Figures ~\ref{state} and \ref{site}.
A site is a region in the plane containing two arcs, each from one of the component curves, equipped with an indicator, marking and connecting the two curves.
The indicator is to be interpreted as an instruction that the two curves at the site are to be colored differently from a set of three distinct colors $\{r,b,p\}.$ In Figure~\ref{state} we illustrate a graph $G$ 
with three edges marked that form a perfect matching (as described in the previous section). Adjacent to $G$ on the right is a state $S$ with three sites that correspond to the edges in the perfect matching.
Below these two parts of the figure, on the left, is a diagram of a state $S'$ obtained by switching two sites of the state $S$ so that $S'$ is colorable. The coloring of $S$ is rewritten as a formation in the 
diagram on the lower right.\\

Each site corresponds, as shown in Figure~\ref{site}, to an edge in a cubic graph. Gven a cubic graph and a perfect matching of its nodes, we can convert each pair of nodes in the matching to a site. See Figure~\ref{state} for an illustration of this remark. In the figure we have marked two edges to be converted and, here we obtain a single Jordan curve $S$ from the replacements.  \\

\noindent {\bf Definition.} A collection of Jordan curves in the plane with a chosen set of sites is called a {\it planar state S}. Given a planar state S, one can attempt to color its curves according to the site indicators. Two curves connected by an indicator must be colored with different colors. A state may be uncolorable. The state S in Figure~\ref{state} is not colorable. We introduce a local operation of {\it site switching} as shown in Figure~\ref{site}. To switch a site we replace two
locally parallel arcs with two crossing arcs. Loops that cross are allowed to have the same or different colors unless a state indicator forces them to be colored differently. On switching a site, the site indicator remains in place as shown in Figure~\ref{site}. Note that in Figure~\ref{site} we switch two sites and arrive at a  state $S'$ that is colorable. This gives a coloring of the original graph $G,$ and a formation $F$ for that graph.  \\

\begin{figure}[htb]
     \begin{center}
     \begin{tabular}{c}
     \includegraphics[width=8cm]{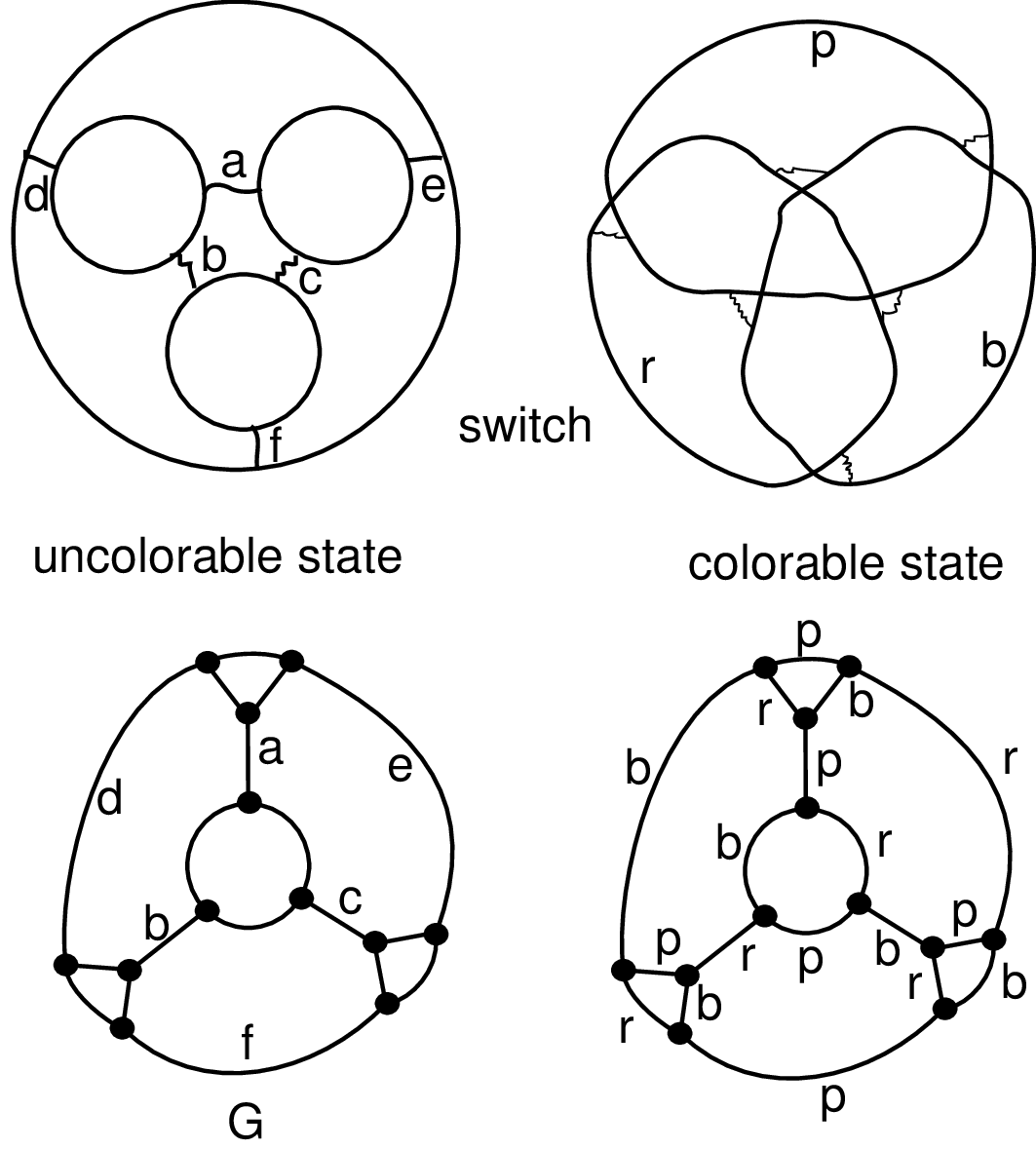}
     \end{tabular}
     \caption{\bf Coloring a Planar State and Its Graph.}
     \label{three}
\end{center}
\end{figure}

At the bottom left of Figure~\ref{site} we see a state that is uncolorable but can be converted to a colorable state by switching. The state itself is not directly colorable since it consists in four mutually touching planar loops. There is only one way to color this example. One must switch all of the sites as shown in Figure~\ref{three}. One obtains
a confguration of three loops and a coloring of the corresponding graph. \\

\begin{figure}[htb]
     \begin{center}
     \begin{tabular}{c}
     \includegraphics[width=8cm]{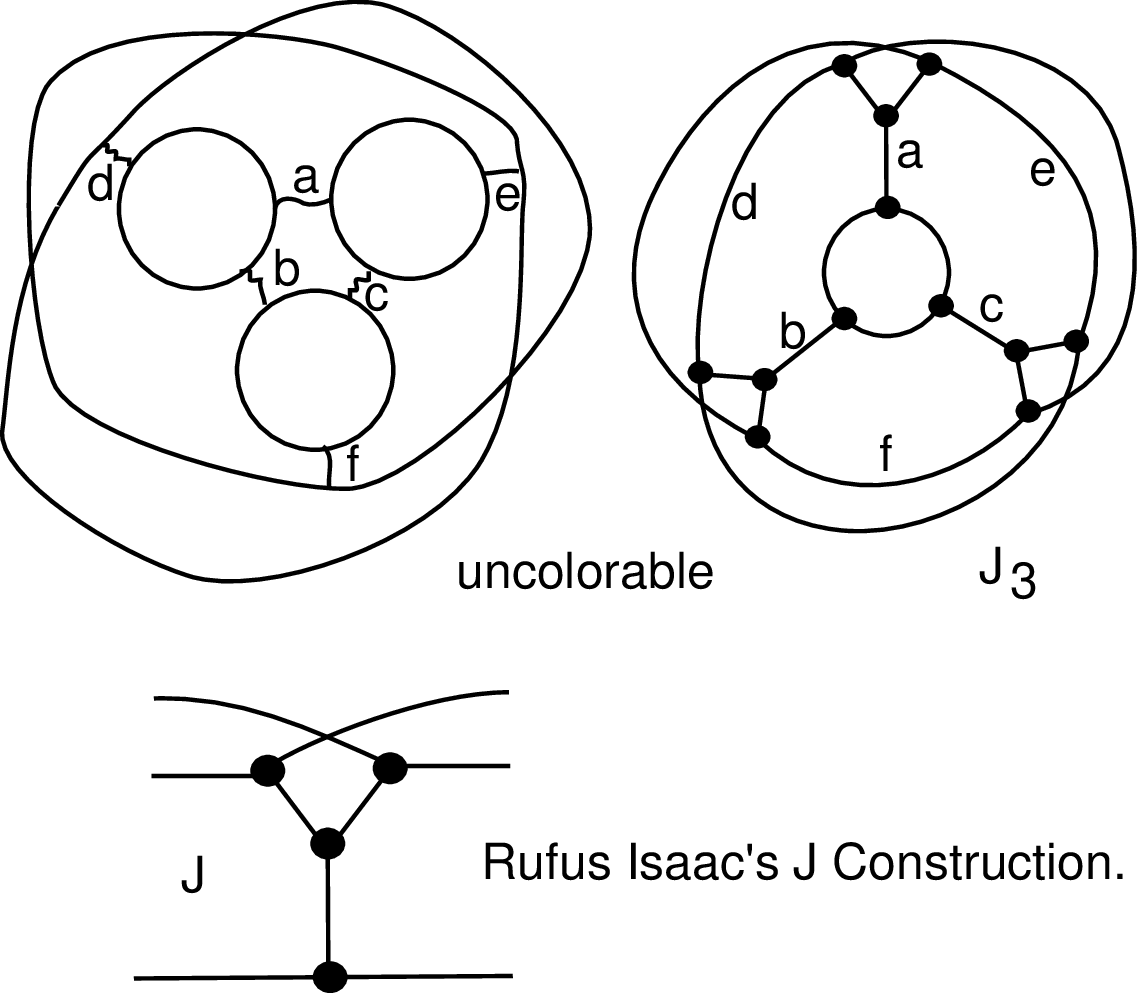}
     \end{tabular}
     \caption{\bf A NonPlanar State Corresponds to Isaac's $J_{3}.$ }
     \label{IsaacJ}
\end{center}
\end{figure}

The other example in Figure~\ref{site} is a non-planar variant of the example just discussed.  In Figure~\ref{IsaacJ}  we show the graph corresponding to this state. The graph is $J_{3},$ the well-known uncolorable graph constructed by Rufus Isaacs \cite{Isaacs} as a circular composition
of three copies of the special graphical element $J$ shown in this figure. An odd number of the $J$ elements placed in a circular configuration will always be uncolorable. These graphs are called $J_{2n+1}$ for $n = 1,2, \cdots .$ Uncolorable, isthmus-free, cubic graphs were named {\em snarks} by Martin Gardner \cite{M} inspired by  
Lewis Carroll's poem ``The Hunting of the Snark" \cite{LC}. Isaac's snarks are fundamental. The graph $J_{3}$ can be 
collapsed to the minimal non-planar uncolorable, the Petersen graph. We will return to this example and discuss the Petersen graph below. \\

 We say that a planar site is an {\it isthmus} if there is a pathway from one side of the site to the other that does not cross any of the Jordan curves.  
{\it We can now give a statement of our state-calculus reformulation of the 
map theorem.}\\

\noindent {\bf Planar State Theorem.} Let $S$ be a planar state with no isthmus. Then, after switching some subset of the sites of $S$ to form a new (possibly non-planar) state $S'$, the state $S'$ is colorable with three colors.\\

\noindent {\bf Proof.} See \cite{Kauffman} for the proof of this result.    $\hfill\Box$
\\

\begin{figure}[htb]
     \begin{center}
     \begin{tabular}{c}
     \includegraphics[width=8cm]{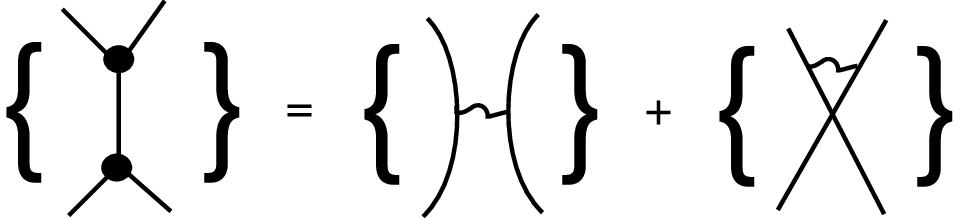}
     \end{tabular}
     \caption{\bf A Logical State Expansion.}
     \label{equation}
\end{center}
\end{figure}

Now we will use the state calculus to give a (tautological) expansion formula for colorings of arbitrary cubic maps. The formula is given below and in Figure~\ref{equation}.

$$\{ \YGlyph  \} = \{ \VGlyph \} + \{  \CGlyph  \}$$

The meaning of this equation is as follows: $\{ G \}$ denotes the number of proper three-colorings of the edges of a cubic graph $G.$ The graph $G$ is not assumed to be planar or represented in the plane.
The left-hand term of the formula represents a graph $G$ and one specific edge of $G.$ The two right hand terms of the formula represent the graph $G$ with the specific edge replaced either by a parallel connection of edges as illustrated, or by a cross-over connection of edges. In each case we have added the site-symbol, used in this section, to indicate that the two edges, in a coloring of the new graph, are to be colored with different colors. Thus the two graphs on the right-hand side of the formula have acquired a local state site of the type discussed in this section. The rest of the graphs may still have nodes. 
{\it One can define $\{ G \}$ directly as a summation over all states of the graph $G$, if we choose a perfect matching for $G$ and sum over all replacements of the edges in the perfect matching by parallel and crossed local sites.} Note that the crossed lines in the diagram for the formula are only an artifact of representation, and are not new nodes in these graphs. The truth of the formula follows from our previous discussion of the fact that a coloring of an edge will result in a double $Y$ form (as on the left-hand side of the formula) so that the two colors at the top of the $Y$ either match in parallel or match when crossed. If we are given the pattern of matched colors that are different, then we can reconstruct an edge and give it that color that is different from both of them. We shall refer to this formula as the 
{\it logical state expansion for three-coloring of cubic graphs.} This completes our review of the work in \cite{Kauffman}.\\

\section{New Penrose Evaluations and Polynomials}
\begin{figure}[htb]
     \begin{center}
     \begin{tabular}{c}
     \includegraphics[width=10cm]{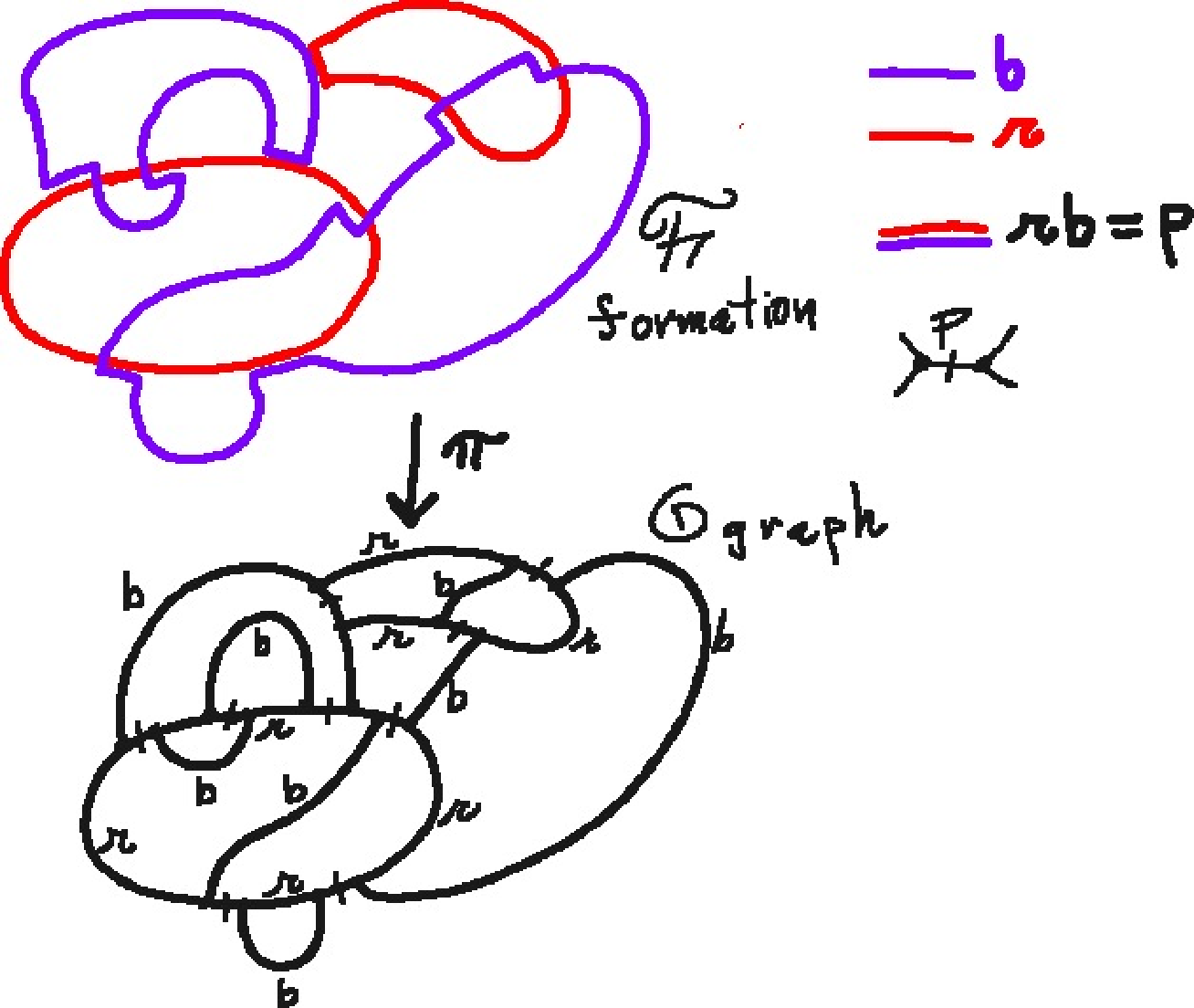}
     \end{tabular}
     \caption{\bf Standard Formation and Graph}
     \label{Form1}
\end{center}
\end{figure}

\begin{figure}[htb]
     \begin{center}
     \begin{tabular}{c}
     \includegraphics[width=10cm]{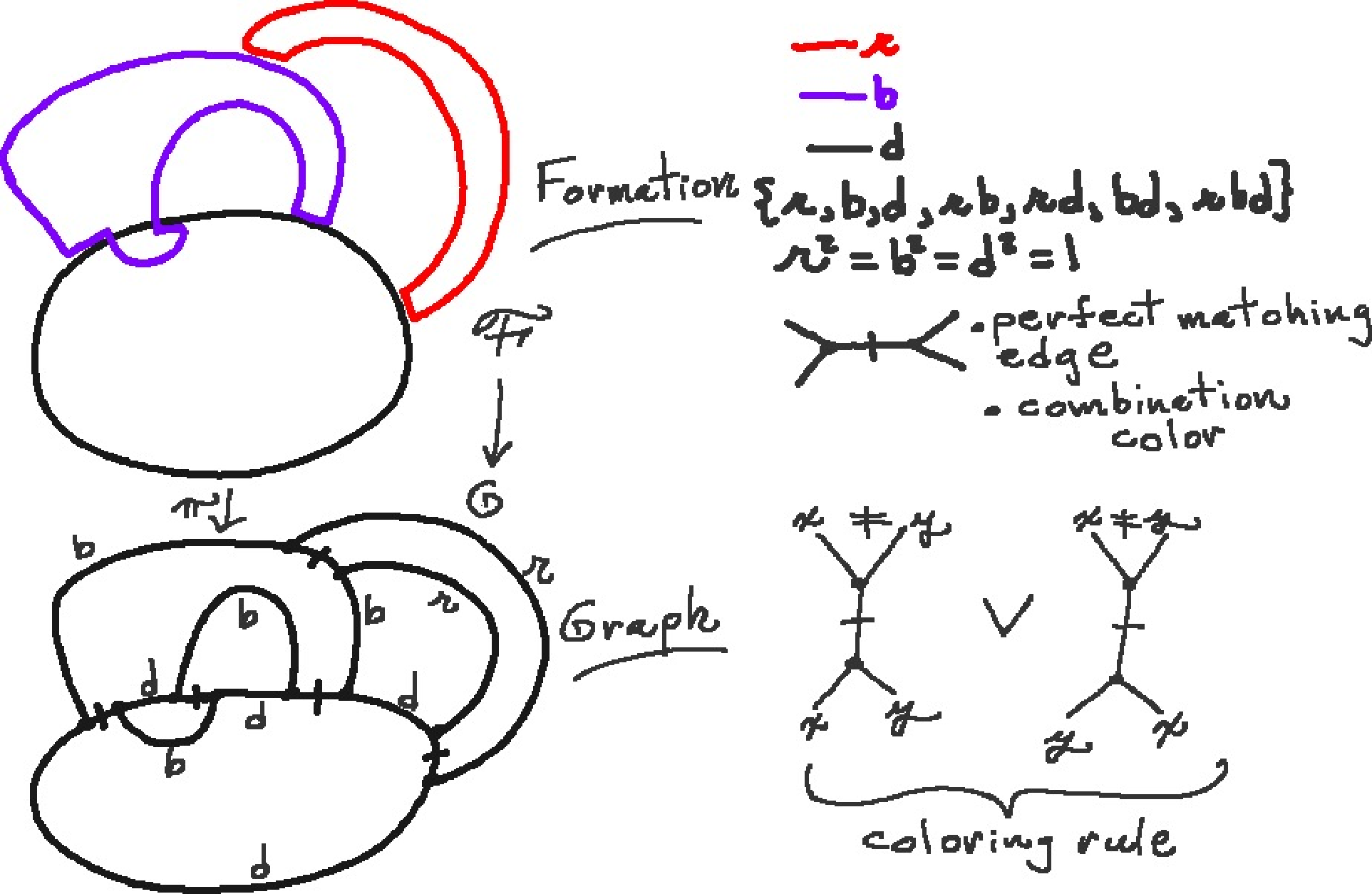}
     \end{tabular}
     \caption{\bf Generalized Formation and Graph}
     \label{Form2}
\end{center}
\end{figure}

In this section we generalize the concept of proper coloring of a trivalent graph with three colors to a notion of $n - coloring$ the graph using a palette of 
$n$ colors. See \cite{BKM} for a discussion of this concept from a different viewpoint. In order to explain this method, we direct the reader first to Figure~\ref{Form1} and Figure~\ref{Form2}. We point out a natural generalization of formations that uses $n$ colors.\\

In Figure~\ref{Form1} we show the formation version of three coloring using loop colors r (red) and b (blue) and the superimposed color p= rb (purple). There is a one-to-one correspondence between three colorings of a trivalent graph $G$ (with three distinct colors at each node) and {\it formations} of the graph, consisting in collections of loops in the colors red and blue (only) that are superimposed along arcs of the circles to form rb= red-blue  = purple = p assignments to these arcs. The formation is obtained from a colored graph by tracing loops in red and in blue, using the factorization of purple into red and blue. A colored graph is obtained from a formation by associating nodes to the ends of the arcs of two-color interaction in the formation.  The figure  illustrates a graph, a coloring of the graph and a corresponding formation for this coloring. Note that the set of purple edges of the graph forms a perfect matching structure for it. If we had fixed an even perfect matching for the graph (for example the one illustrated in the figure),then we could enumerate colorings of the remaining edges in just red and blue to obtain the subset of all colorings that used purple on these perfect matching edges. It would be of interest to make this two-color enumeration for a given choice of perfect matching.\\

We define a polynomial generalization of the Penrose evaluation by first choosing a perfect matching $M$ for a given graph $G,$ and then
expanding only on the perfect matching edges of $G$. By making this restriction, we obtain a polynomial that is well-defined for a given choice of perfect matching. The polynomial is defined for an immersion of $G$ in the plane, using  generalized Penrose recursion formulas described as follows
$$P(G,M)( \YMGlyph ) = A~P(G,M)(\VDiag ) + B~P(G,M)( \CCircleDiag  ),$$
$$P(G,M)(O) = \delta,$$
$$P(G,M)(\CBoxDiag) = 2 P(G,M)(\CDotDiag) - P(G,M)(\CCircleDiag).$$

The evaluation takes the form of a Penrose expansion except that we keep track of the original immersed crossings, denoted by $\CBoxDiag$ and the last relation expands further each original immersed crossing in terms of an ordinary (virtual) crossing in the expansion, $\CCircleDiag$, and a {\it fused} crossing in the form $\CDotDiag.$ If two loops are joined at a fused crossing
then they together contribute $\delta,$ the same as a single loop. In general, a complex of loops connected by fusions contributes only $\delta.$ This is the combinatorial definition of a three variable $(A,B, \delta)$ generalization of the Penrose evaluation. It is well defined for arbitrary trivalent graphs with chosen perfect matching. The result of the evaluation is independent of the immersion into the plane that is used for the computation.\\

In the rest of this section we specialize the polynomial $P(G,M)$ so that $A=1, B=-1, \delta = n$ where $n$ is a positive integer. Then the evaluation of $P(G,M)$ is interpreted as a color count to be described below.\\

{\it The concept of a formation generalizes to $n$ colors.} One has disjoint loops in the plane, each colored by one of the $n$ colors. Loops of different colors can share intervals and either cross or not cross one another just as before. In Figure~\ref{Form2} we show a formation where there are three colors of loops: red (r) , blue (b)  and dark (d). On making colored graph from such a formation, each node receives three distinct colors, one of which is a composite color just as purple was a composite of red and blue in the original formation. In this generalization, if x and y are the colors on two loops at an interaction vertex then the third edge at the vertex is labeled by xy.
In the figure we began with the basic set of colors $\{ r, b, d \}$ and this extends to include the superposition colors $\{ rb, rd , bd \}.$ Once again, the superposition colors occur on a collection of edges in the graph $G$ that form a perfect matching for $G$ (not necessarily an even perfect matching in this generalization). In the coloring of the graph shown in the figure we have indicated the perfect matching edges, but we have not colored them. They inherit the product of the other two colors that impinge on the nodes at the ends of the edge. Note that in colorings of the graph that are inherited from the $PM$ the two colors at either end of a perfect matching edge are the same. This fact is illustrated in the figure and we take it as the {\it coloring rule for associating $n$ colors to a given graph with perfect matching structure}. Thus $n$ colors are used for the edges outside the set of $PM$ edges, and for any $PM$ edge there are two distinct colors that occur at each of its end nodes. With this rule for coloring perfect matching graphs, we have a one to one correspondence between generalized formations for $G$ with $n$ colored loops and  colorings for $PM$ structures on $G.$\\

We now explain how our generalization of the Penrose formula can be used to count $n$ colorings of a given graph $G$ with a choice of perfect matching.\\

We define $PK(G,M)[n],$ a coloring polynomial in the variable $n$ (that can be taken to be a positive integer), defined for trivalent graphs $G$ with perfect
matching $M.$ This polynomial is a generalization of the $PK$ evaluation at $n=3$ studied in \cite{Kauffman}. The special case at $n=3$ counts the number of edge $3$-colorings of an arbitrary
trivalent graph $G$ (no perfect matchng required) via a generalization of the original Penrose evaluation \cite{Kauffman,P}. The key point about the evaluation  $PK(G)[3]$ is that it
gives the total number of colorings of the graph $G$ for an immersed representation of $G$ in the plane, and it follows the original Penrose expansion, with an extra caveat for the singularities
of the immersion. In our generalization to the perfect matching polynomial $PK(G,M)[n]$ we will follow the same procedure and we will obtain a count of special colorings of the perfect matching graph 
using $n$ colors. This coloring definition and the polynomials will be described below. The definition of the perfect matching polynomial is a combination of our method for extending the Penrose evaluation and (Baldridge, Kauffman, McCarty) methods for making perfect matching polynomials from Penrose-type evaluations \cite{BKM,BKR}.\\

\begin{figure}[htb]
     \begin{center}
     \begin{tabular}{c}
     \includegraphics[width=10cm]{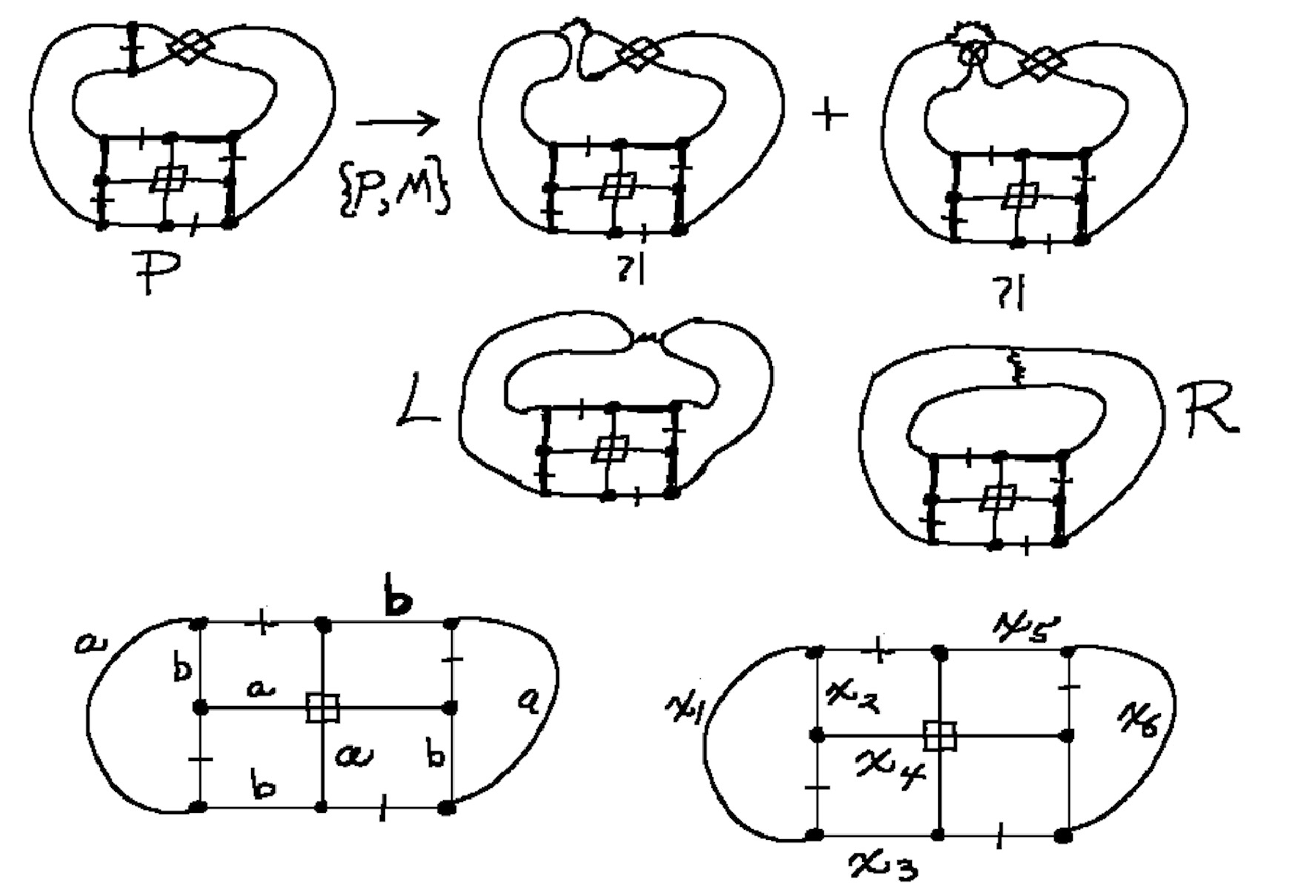}
     \end{tabular}
     \caption{\bf Petersen Graph is Strongly Uncolorable.}
     \label{Petersen}
\end{center}
\end{figure}

\begin{figure}[htb]
     \begin{center}
     \begin{tabular}{c}
     \includegraphics[width=12cm]{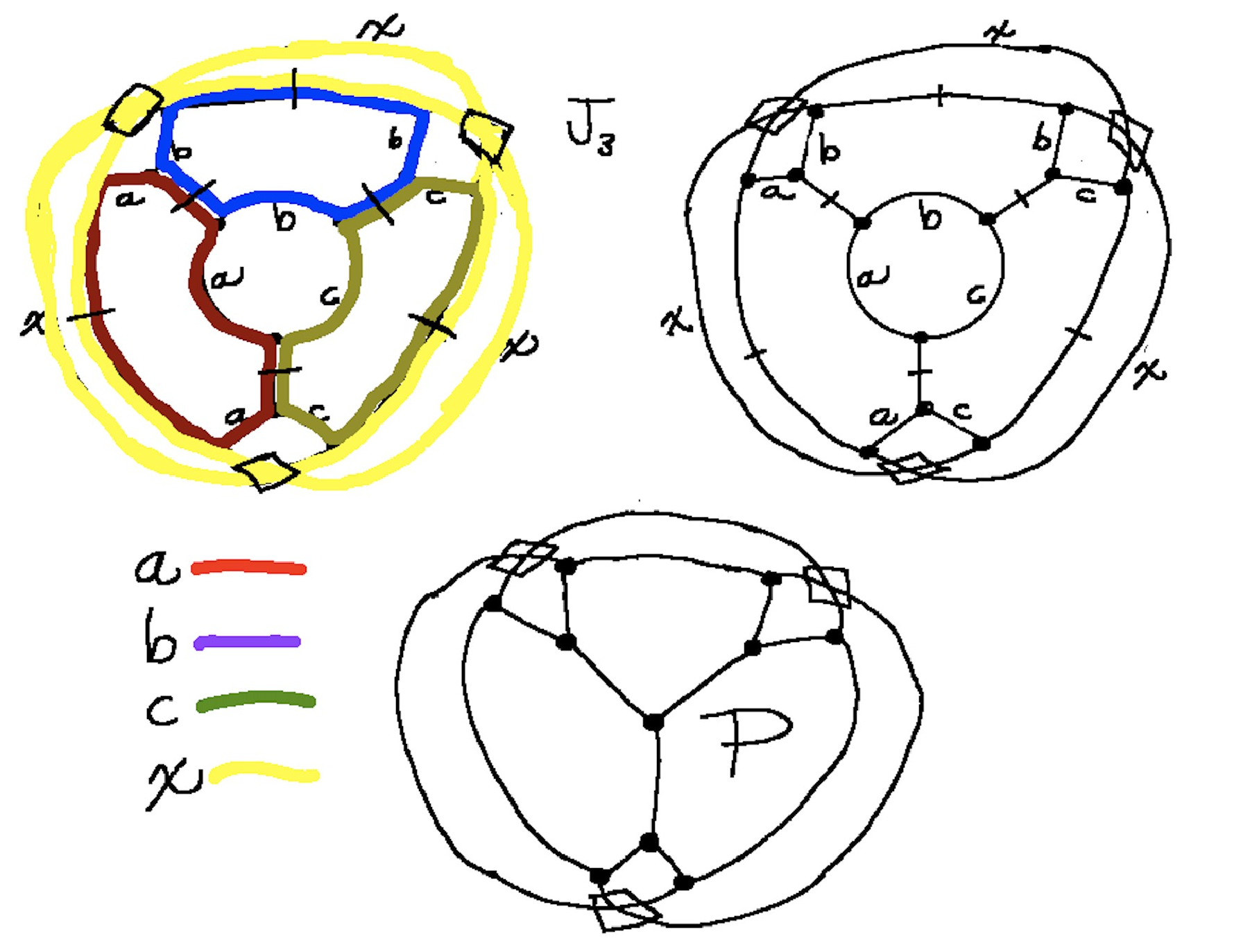}
     \end{tabular}
     \caption{\bf Isaacs $J_{3}$ can be PM-colored with four colors (but not with three colors).}
     \label{Isaacs}
\end{center}
\end{figure}

\begin{figure}[htb]
     \begin{center}
     \begin{tabular}{c}
     \includegraphics[width=11cm]{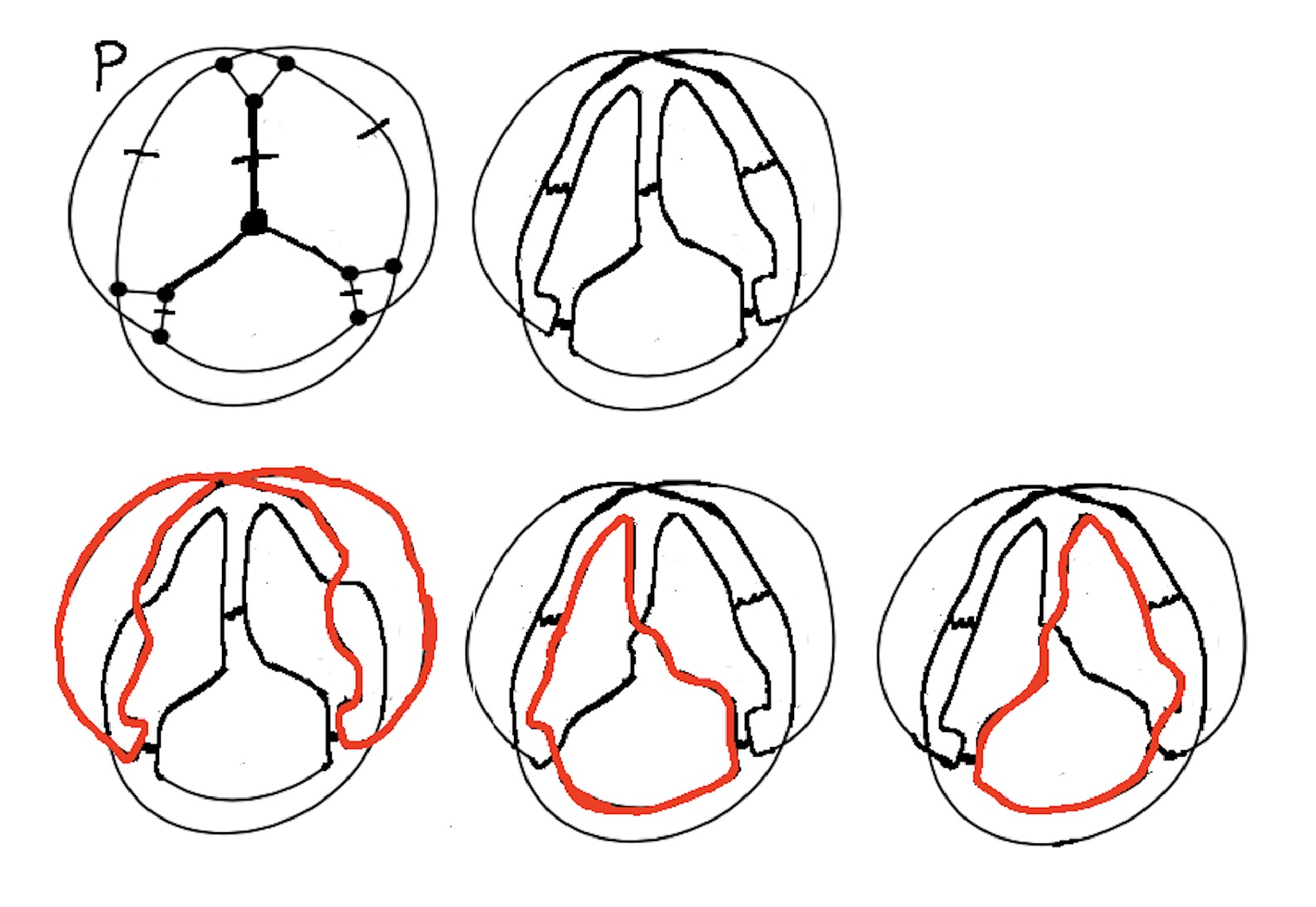}
     \end{tabular}
     \caption{\bf Strong Uncolorability of the Petersen Graph}
     \label{pgnarl}
\end{center}
\end{figure}

We begin by describing an {\it n-coloring} of a finite trivalent graph $G$ with perfect matching $M.$ A coloring consists in an assignment of colors to the {\it non-matching edges} of $G$ from the color set
$\{1,2,...,n\}$ so that only two distinct colors are touching  each matching edge and these  two colors both appear at each end of the matching edge. We call such a coloring of $(G,M)$ a perfect matching $n$-coloring.\\

We point out first a purely combinatorial interpretation of the coloring count for $(G,M).$ Define the {\it matching polynomial} $Match(G,M)$, also denoted by  $\{G,M \}$, by the recursion

$$\{ \YMGlyph  \} = \{ \VGlyph \} + \{  \CGlyph  \},$$
$$\{G O \} = n \{G\},$$
$$\{ O \} = n.$$

where it is understood that $$\{ G,M\} = \sum_{S} \{S\}$$ where each matching edge has been replaced by the glyphs in the recursion about to form a collection of 
state configurations consisting in loops connected by the wiggly glyphs in the form
$ \{ \VGlyph \} $ and $\{  \CGlyph  \}.$ The evaluation of a state $S$ is defined to be equal to the number of ways to color the loops in $S$ so that two distinct loops joined by 
a wiggly glyph are colored differently. Thus the state evaluation is a chromatic polynomial of the graph $G(S)$ associated with the state
by assigning a node to each loop in the state, and an edge to each pair of loops that share a wiggly line (indicating that the loops are colored differently). We use the color set of $n$ colors described above.\\

\noindent{\bf Remark.}. The matching polynomial can be reformulated as a Penrose-type expansion in the form shown below, where it is
given that loops that share a dark dotted intersection must be colored by the same color. Lines without wiggly lines or dotted intersections can be colored with the same color or with different colors. This formula is of interest, as it gives another solution to the problem to give a Penrose type formula that counts colors for all trivalent graphs. Consequences of this point of view will be 
studied in subsequent work. Note that for this formulation $G$ is any trivalent graph.
$$\{ \YMGlyph  \} = \{ \VDiag \} + \{  \CDiag  \} - 2\{ \CDotDiag \},$$
$$\{G O \} = n \{G\},$$
$$\{ O \} = n.$$\\

Since, by its definition, $\{G,M\}$ counts those colorings of the perfect matching graph so that exactly two distinct colors appear at each matching edge satisfying our conditions
for an $n$-coloring of $(G,M),$ it follows that $\{G,M\}$ is equal to the number of $n$-colorings of the perfect matching graph.\\

Note that $\{G,M\}$ is defined independent of any planar immersion of the graph $G.$  \\

Now we turn to the polynomial $PK(G,M)[n].$ This polynomial is defined by choosing an immersion of $G$ in the plane and then using the generalized Penrose recursion formulas described as follows
$$PK( \YMGlyph ) = PK(\VDiag ) - PK( \CDiag  ),$$
$$PK(O) = n,$$
$$PK(\CBoxDiag) = 2 PK(\CDotDiag) - PK(\CDiag).$$
These recursions mean that the evaluation takes the form of a Penrose expansion except that we keep track of the original immersed crossings, denoted by $\CBoxDiag$ and then the last relation
expands further each original immersed crossing in terms of an ordinary crossing in the expansion, $\CDiag$, and a {\it fused} crossing in the form $\CDotDiag.$ If two loops are joined at a fused crossing
then they together contribute $n,$ the same as a single loop. In general, a complex of loops connected by fusions contributes only $n.$ This is the combinatorial definition.
We now prove that $PK(G,M)[n]$ counts $n$-colorings of $(G,M).$\\

\noindent {\bf Theorem.}  $PK(G,M)[n] = \{ G,M\}.$\\

\noindent {\bf Proof.} Define two diagrammatic tensors as shown below.
The indices run in the set $\{1,2,\cdots n\}$ for $n$ colors.\\

$ \YMTens = 1$ if $a=c, b=d, a \ne b,$

$\YMTens =-1$ if $a=d, b=c$ and $a \ne b,$ 

$\YMTens = 0$ otherwise.

$\IMTens  = 1$ if $a=d, b=c, a=b,$

$\IMTens  = -1$ if $a=d,b=c, a\ne b,$

$\IMTens  = 0$ otherwise.\\

Note that with $\delta^{a}_{b}$ denoting a Kronecker delta, we have the formula
$$\YMTens = \delta^{a}_{c}\delta^{b}_{d} - \delta^{a}_{d}\delta^{b}_{c}$$
and that 
$$\IMTens = 2 \FTens - \CTens = 2 \FTens - \delta^{a}_{d}\delta^{b}_{c},$$
where
$\FTens$ is equal to $1$ only when $a=b=c=d$ and is $0$ otherwise.\\

Define $[G,M]$ as the tensor contraction of $(G,M)$ with respect to these tensors in the sense of Penrose \cite{P}. That is,
$[G,M]$ equals the sum over all possible index assignments to the {\it non-matching edges} of $G$ where we take the product of tensor values for each assignment of the
indices. It follows from the tensor definitions that in order for an index assignment to contribute to the summation, it must be a coloring of a state $S$ of the matching polynomial for $(G,M).$  
The contribution is (by the above assignments) equal to $(-1)^{A + B}$ where $A$ is the number of crossed glyph, $\CGlyph$, contributions, and $B$ is the number of immersion tensors, $\IMTens$, with 
$a\ne b.$ By the Jordan Curve Theorem (since the graphs are in the plane), $A+ B$ is even, and hence each state contributes $+1$ to the summation. This proves that 
$[G,M] = \{ G,M\}.$ On the other hand, it follows from the tensor definitions that $[\YMGlyph] = [\VDiag] - [\CDiag]$ and $[\CBoxDiag] = 2[\CDotDiag] - [\CDiag].$ Thus
$[G,M] = PK(G,M)[n].$ This completes the proof. $\hfill\Box$ \\

\begin{figure}
     \begin{center}
     \begin{tabular}{c}
     \includegraphics[width=11cm]{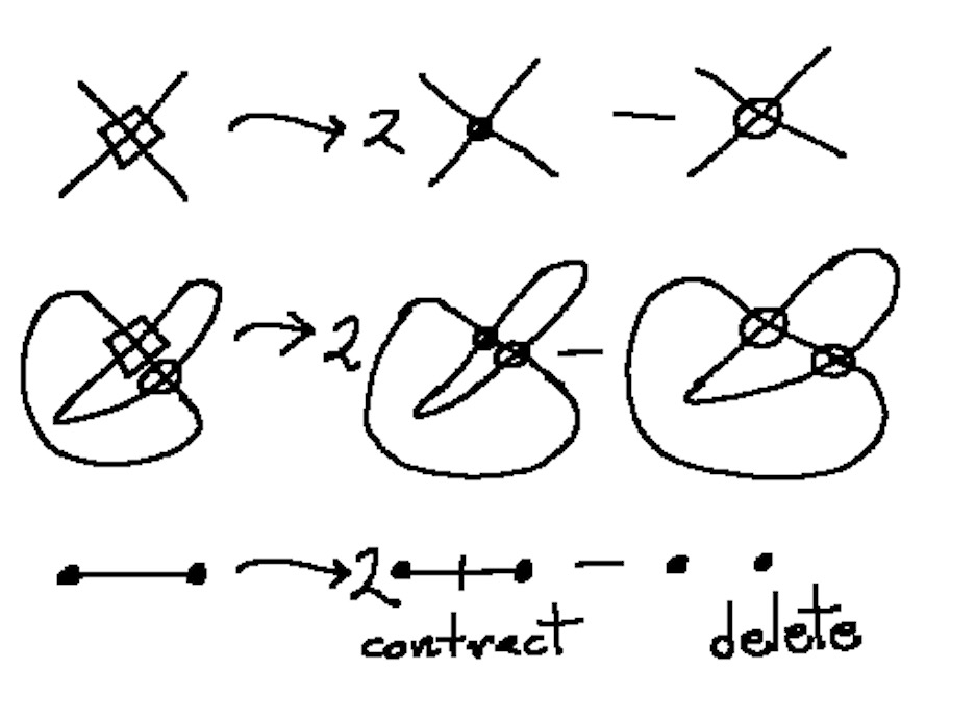}
     \end{tabular}
     \caption{\bf Contraction Deletion for Graphical Evaluation}
     \label{EFF111}
\end{center}
\end{figure}

\begin{figure}
     \begin{center}
     \begin{tabular}{c}
     \includegraphics[width=11cm]{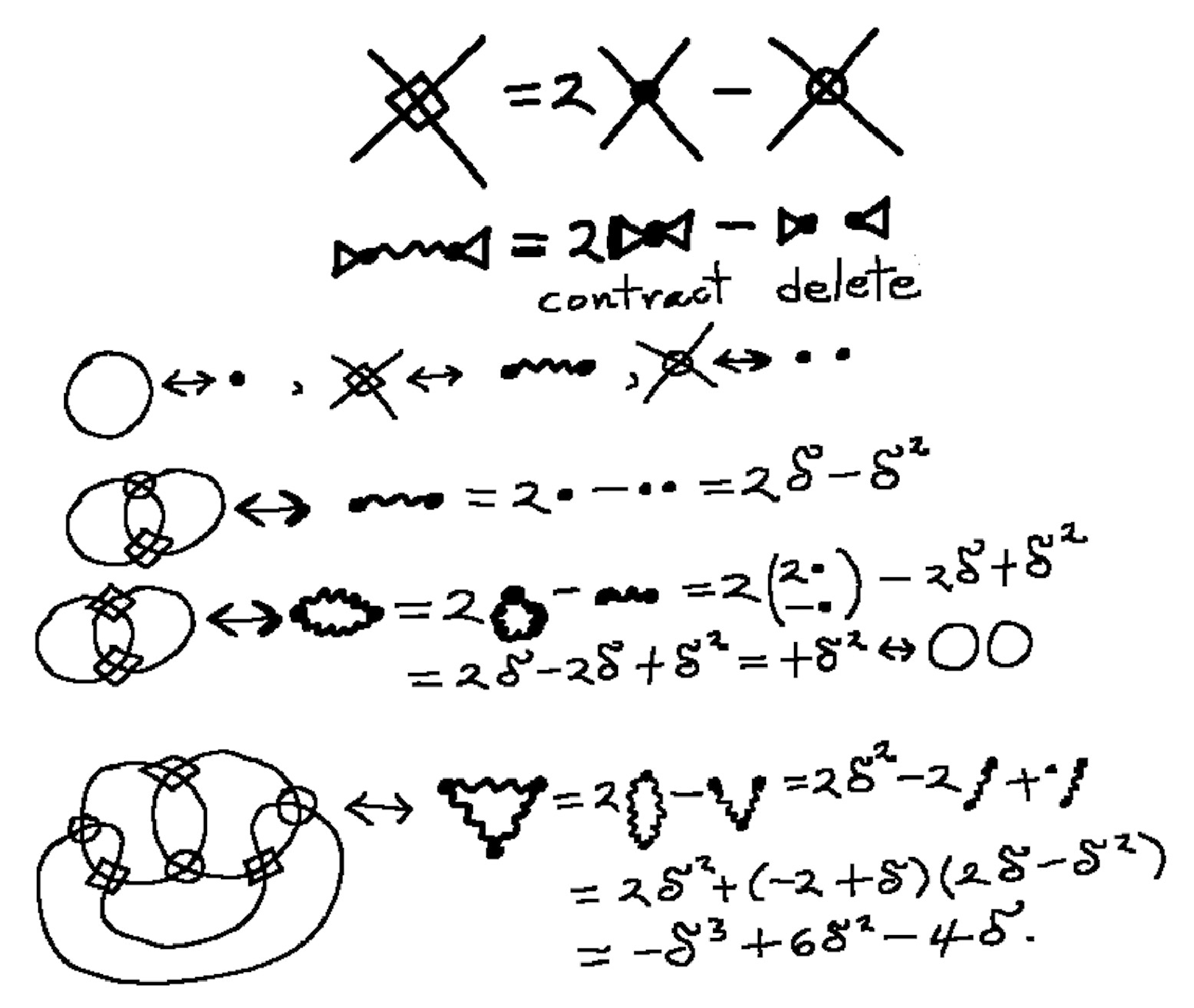}
     \end{tabular}
     \caption{\bf Contraction Deletion for Graphical Evaluation}
     \label{EFF11}
\end{center}
\end{figure}

\begin{figure}
     \begin{center}
     \begin{tabular}{c}
     \includegraphics[width=11cm]{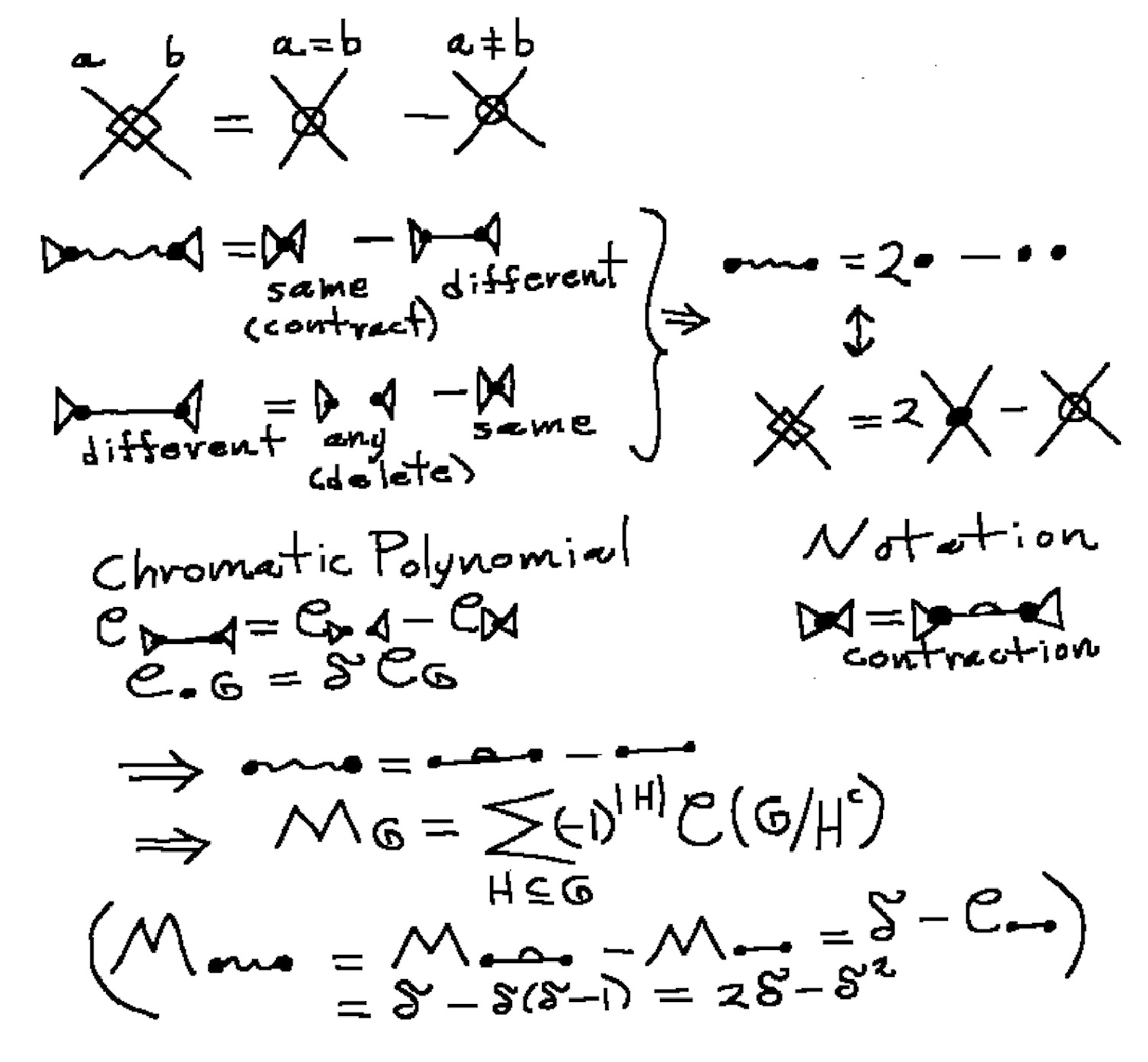}
     \end{tabular}
     \caption{\bf Tautologies of Expansion, Contraction and Deletion}
     \label{Taut}
\end{center}
\end{figure}

\begin{figure}
     \begin{center}
     \begin{tabular}{c}
     \includegraphics[width=11cm]{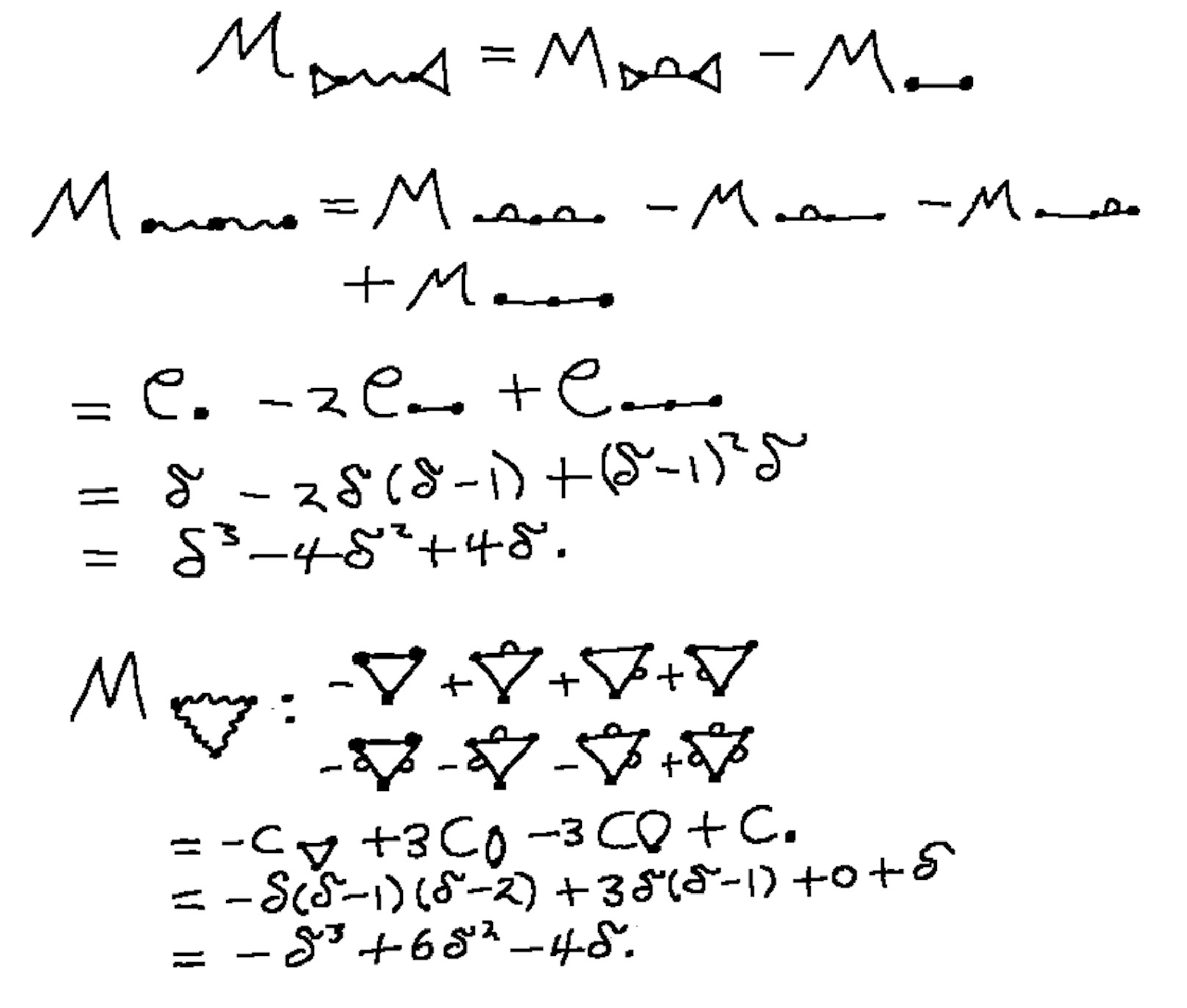}
     \end{tabular}
     \caption{\bf Example Chromatic Evaluation}
     \label{Taut1}
\end{center}
\end{figure}

\begin{figure}[htb]
     \begin{center}
     \begin{tabular}{c}
     \includegraphics[width=8cm]{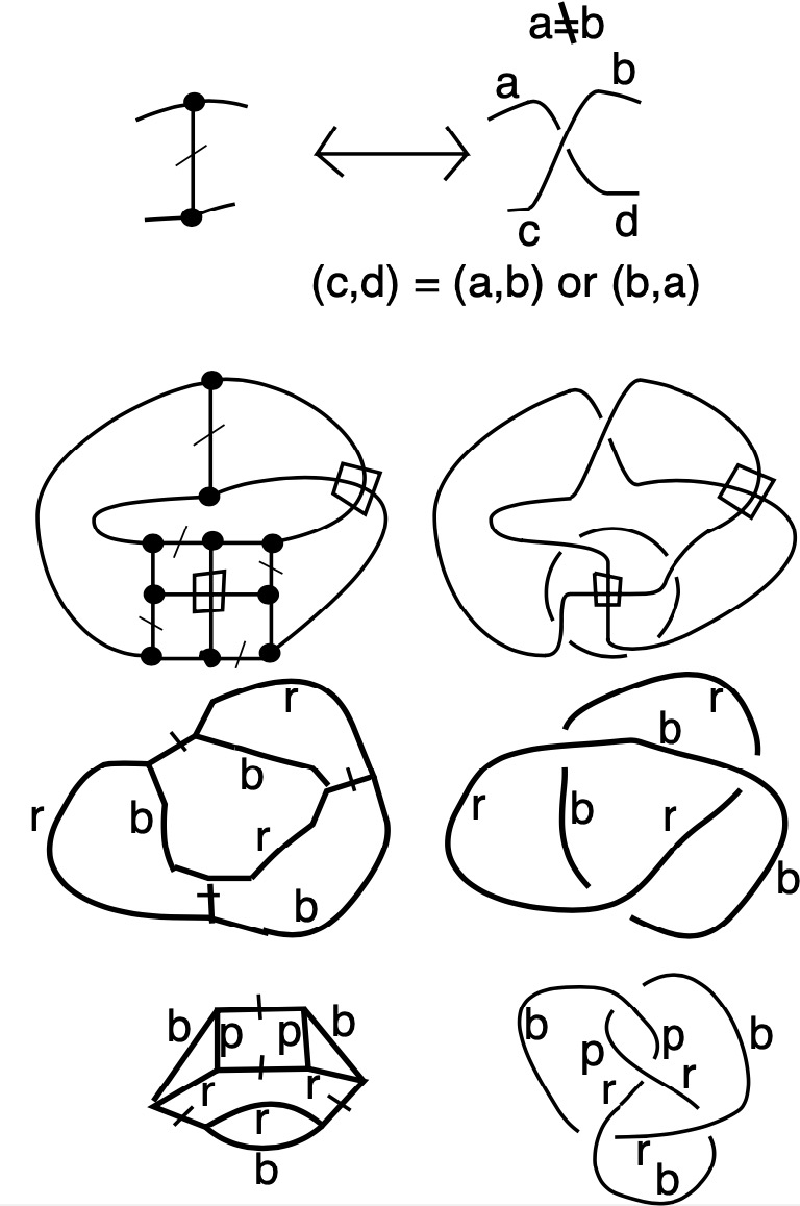}
     \end{tabular}
     \caption{\bf Knot Code for the Petersen Graph, an Unknot and a Figure Eight Knot}
     \label{pcode}
\end{center}
\end{figure}

\noindent {\bf Examples.} We end this part of the paper with some examples and remarks about evaluations. In particular, we point out how the extra evaluations we have discussed for
immersions of non-planar graphs can be seen as evaluations of the Tutte dichromatic polynomial for graphs associated with the diagrams, and that this dichromatic evaluation can be seen as an alternating sum of chromatic polynomials. Furthermore, the coloring rules, written in graphical language, can be translated into coloring rules for knot and link diagrams. This provides the beginning of a transition to knot and link theory that constitutes the remainder of the paper.\\

\begin{enumerate}
\item Examine Figure~\ref{Petersen}. This figure shows how the Petersen graph is strongly uncolorable. We have chose a specific perfect matching $M,$ and we show the tautological coloring expansion for $(P,M)$  along one perfect matching edge. The two terms of the expansion are labeled $L$ and $R.$  $L$ and $R$ are combinatorially identical and we must show that the arcs joined by a
wiggly arc are forced to have the same color (since a coloring of the Petersen graph would entail these arcs having different colors). The next part of the figure shows the verification of this fact.
If arc $x_1$ is labeled $a$ and arc $x_2$ is labeled $b$ with $a \ne b,$ then it follows from the coloring rules that $x_4 = a.$ Note that this is the case using the generalized $PM$-coloring rules, since $a$ and $b$ appear at one end of a perfect matching arc and $a$ appears at the other end, we must have an appearance of $b$ at the other end.  By the same reasoning, $x_3 = a, 
x_4 = a, x_5 =b, x_6 = a,$ and we find that the colors are forced as claimed. While we have shown here that the Petersen is strongly uncolorable for a given perfect matching $M,$ the reader will find that this same phenomenon occurs for any choice of perfect matching for the Petersen Graph.\\

\item In Figure~\ref{Isaacs} we show a perfect matching for the Isaacs $J_3$ graph and a $PM$- coloring using four colors. The coloring is illustrated as a formation with overlap coloring along the perfect matching edges.The figure also shows how $J_3$ collapses to the Petersen graph. It is of interest to point out that both $J_3$ and Petersen are uncolorable at n = 3, but $J_3$ can be colored using more colors, while the Petersen is strongly uncolorable and cannot be colored with $n$ colors for any $n.$\\

\item Figure~\ref{pgnarl} shows a perfect matching for the form of Petersen graph in Figure~\ref{Isaacs} and also shows representative states in the tautological expansion for this perfect matching. As the reader can see each state contains self-interaction sites. It is the case that every state in the tautological expansion of the Petersen contains self-interactions, and it is this property that prevents coloring for any $n.$\\

\item We note that the virtual box expansion can be regarded as a contraction-deletion expansion as shown in Figure~\ref{EFF111} and Figure~\ref{EFF11}. Specifically, to each loop state $S$ of the chromatic expansion
we associate a graph $G(S)$ whose nodes are in $1-1$ correspondence with the loops in the state. We assign an edge between nodes if they meet in a boxed virtual crossing. No edge occurs for the round virtual crossings. Note that we denote the edges in this fashion $\Wigg.$ This will allow us to use the $\Wigg$ for a contraction/deletion expansion corresponding to our use in the generalized Penrose evaluation, and to use a straight version of an edge for the usual coloring count expansion, as we explain in the next item. Then the $\Wigg$ contraction-deletion expansion is as shown in the Figure~\ref{EFF11}. We will also use $\Contract$ to denote a contracted edge. Then, if $M_{G}$ denotes the corresponding graphical evaluation, we have the formula $$M_{\Wigg} = 2 M_{\Contract} - M_{\Delete}.$$ The expansion formula for the boxed crossing corresponds to this expansion formula, which in turn is a 
{\it dichromatic polynomial evaluation} $(-1)^{b(S)}Z_{G(S)}(\delta, v = -2)$ where $b(S)$ denotes the number of box virtual crossings in the state $S$, which is equal to the number of edges
in $G(S).$  The dichromatic polynomial is  defined by the 
recursion
$$Z_{G} = Z_{G'}  + v Z_{G''},$$
$$Z_{G \cup \bullet} = \delta  Z_{G}.$$
Here $G'$ is the result of deleting an edge from $G,$ and $G''$ is the result of contracting that same edge. The $\bullet$ denotes a node disjoint from the rest of the graph. 
This means that we evaluate a dichromatic polynomial \cite{LKStat} for each state, and that the generalized Penrose evaluation is a sum with coefficients over these dichromatic evaluations. This is a useful direction of reformulation and will be followed up in subsequent papers.

It is interesting to note that if one uses the dichromatic polynomial to represent the partition function for the 
Potts model in statistical mechanics \cite{LKStat}, then $v = e^{\frac{-1}{kT}} - 1$ where $k$ is Boltzmann's constant and $T$ is the temperature variable in the model. When $v=-2$ we have
$-2 = v = e^{\frac{-1}{kT}} - 1$ so that $e^{\frac{-1}{kT}} = -1.$ Thus we can set $i \pi = \frac{-1}{kT} $ and conclude that this special point of evaluation of the dichromatic polynomial corresponds to imaginary temperature in the Potts model $ T = i \frac{1}{k \pi}.$ We will investigate this observation in a subsequent paper.\\

\item Contraction and deletion give rise directly to evaluation formulas. Let $\Wigg$ denote an edge the graph associated with boxed virtual crossing, as explained above. 
In Figure~\ref{Taut} we go back to the original definition of the box virtual crossing tensor. It is positive when the one color goes through the box, and it is negative when distinct colors go through the box (no more than two colors can appear at the box as described earlier in the paper). This means that for the corresponding edge in a graph $G = G(S)$ we have the evaluation formula
$$M_{\Wigg} = M_{\Contract} - M_{\Edge}$$
since a contracted edge corresponds to same color and a straight edge can stand for its end nodes having distinct colors.
For such a straight edge we also have the tautological color expansion
$$M_{\Edge} = M_{\Delete} - M_{\Contract}$$
since a deleted edge allows the end nodes to be the same or different, and this last expansion is the graphical counterpart of the logical statement that ``Different  $=$ All $-$ Same".
Combining these two identities we have
$$M_{\Wigg} = M_{\Contract} - (M_{\Delete} - M_{\Contract}) = 2 M_{\Contract} - M_{\Delete},$$
which is exactly our evaluation formula as explained above.\\

Now go back to the basic expansion corresponding to the box tensor $$M_{\Wigg} = M_{\Contract} - M_{\Edge}.$$ Given an arbitrary graph $G$ with wiggly edges (corresponding to box virtuals) repeated application of this identity will expand to a sum of copies of the graph with each edge either marked as a contraction, or left as a straight edge, multiplied by $(-1)^n$ where
$n$ is the number of straight edges in the marking. See Figure~\ref{Taut1} for an example. Each such marked graph will then be evaluated by the chromatic polynomial expansion
$$C_{\Edge} = C_{\Delete} - C_{\Contract}$$ applied to the graph contracted by the edges labeled as contractions. This description immediately reformulates to a subgraph $H$ of those edges
labeled straight, and the complementary subgraph $H^c$ of contraction edges. Let $G/H^c$ denote the result of contracting all the contraction edges in $G.$ we then have from this discussion that $$M_{G} = \sum_{H} (-1)^{|H|} C(G/H^c)$$  where $|H|$ denotes the number of edges in $H.$ Thus our special evaluation of box virtual states can be seen as an alternating sum of chromatic polynomials. Furthermore, we have (by comparison with the above) proved the formula 
$$M_{G}= (-1)^{|G|}Z_{G}(\delta, v = -2) = \sum_{H} (-1)^{|H|} C(G/H^c)(\delta),$$
where $G$ is any finite graph, $|G|$ is the number of edges of $G$ and $H$ runs over all sub-graphs of $G$ according to the above discussion.
For practical calculation of $M_{G}$ it may be best to use the dichromatic polynomial, but it is certainly of interest that this particular dichromatic evaluation can be expressed as an alternating sum of chromatic polynomials.\\

\item Figure~\ref{pcode} shows the simplest translation of perfect matching graphs to knot diagrams with virtual crossings. In this translation, we replace each perfect matching edge with a knot diagrammatic crossing as illustrated at the top of the figure. In reference to the figure, note that the crossing edges $a$ and $b$ are paired and declared unequal, while the crossing edges $c$ and $d$ are paired and it is delcared that either $a=c$ and $b=d$, or $ b=c$ and $a=d.$  This matches the corresponding condition for the perfect matching edge. With this convention we can articulate a coloring rule for semi-arcs of the knot or link diagram that corresponds directly to the rules for coloring the perfect matching. Here the semi-arcs of diagram are defined by the convention shown in Figure~\ref{pcode}. Color can change at an over-crossing as indicated in that figure.  \\

As the figure shows, the link diagram corresponding to a perfect matching graph may have virtual crossings if the original graph had such crossings, and the link diagram or knot diagram may be topologically knotted or unknotted depending on the distribution of the perfect matching edges. Virtual detour moves (See section 4 for the definition.) do not affect coloring of the knot or link diagram, but topological equivalence of knots (Reidemeister moves) can change the coloring properties. In the figure eight knot diagram at the bottom of the figure, we see that this diagram is the translation of a graph with a perfect matching that is not even. Consequently, all three colors appear on the arcs in the the three-coloring of the figure eight knot diagram. In the translation of the Petersen graph (with perfect matching) to a virtual link diagram, we see that the link consists in a small ring and a larger diagram with virtual crossings. Removal of the ring gives the larger diagram alone, which is uncolorable, but contains an isthmus (up to virtual detour moves on the surface of the two-sphere). The inclusion of the ring stops the isthmus and retains uncolorability. With this translation of perfect matching graphs to link diagrams, we see that coloring theorems can be studied in new terms. Classical link diagrams correspond to planar perfect matching graphs and so {\it classical link diagrams with no isthmus are colorable with three colors (via the above rules)} is an equivalent statement of the four color theorem.\\

\item We can translate our generalized Penrose formula for counting colorings of perfect matchings to the formulas shown below for counting link diagrammatic colorings.
In these formulas two virtual crossings appear (circled and boxed), and the formulas apply to a class of link diagrams with multiple virtual crossings. The rest of this paper is 
an exploration of the topology of multiple virtual crossing diagrams, motivated by these considerations from the coloring of graphs.

$$PK( \CrossStarGlyph ) = PK(\VDiag ) - PK( \CCircleDiag  ),$$
$$PK(O) = n,$$
$$PK(\CBoxDiag) = 2 PK(\CDotDiag) - PK(\CCircleDiag).$$

In the next sections we work with analogous expressions that give topological invariants of link diagrams with multiple virtual crossings.\\

\end{enumerate}

\section{Virtual and Multiple Virtual Knot Theory}
In this section we will work with a generalization of virtual knot theory that has a multiplicity of types of virtual crossing.
We begin with a description of standard virtual knot theory and how it will be modified. Then we use the ideas from previous sections, transposing them so that the two graphical virtual crossing types can be used in a knot theoretic context. As we shall see, there is much more to the multiple virtual knot theory than what initially is seen from the graph theoretic point of view. The paper will develop a number of structures in this topological context, and at the end we will return to the graph theory with new ideas. There are many variants of knot theory that occur in relation to virtual knot theory. The reader can see a survey of some of these variants in \cite{VKT,vkt,DVK,DKT,rotvkt,VKC,SL} and in the books by Manturov and Fenn \cite{VSA,FennBook}.\\

\subsection{Virtual Knot Theory}
Knot theory
studies the embeddings of curves in three-dimensional space.  Standard virtual knot theory studies the  embeddings of curves in thickened surfaces of arbitrary
genus, up to the addition and removal of empty handles from the surface. Virtual knots have a special diagrammatic theory, described below,
that makes handling them
very similar to the handling of classical knot diagrams. Many structures in classical knot
theory generalize to the virtual domain.
\bigbreak  

\begin{figure}[htb]
     \begin{center}
     \begin{tabular}{c}
     \includegraphics[width=10cm]{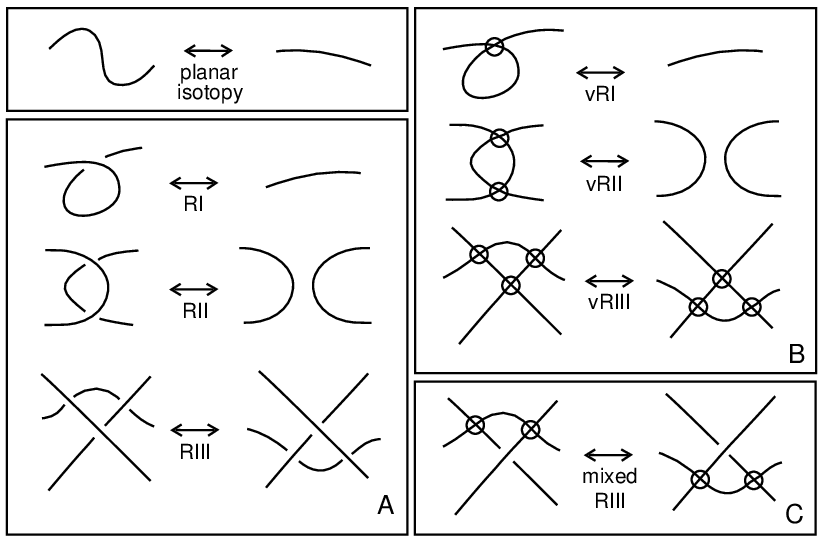}
     \end{tabular}
     \caption{\bf Moves}
     \label{Figure 1}
\end{center}
\end{figure}

In the diagrammatic theory of virtual knots one adds 
a single type of {\em virtual crossing} (see Figure~\ref{Figure 1} that is neither an over-crossing
nor an under-crossing.  A virtual crossing is represented by two crossing segments with a small circle
placed around the crossing point. 
\bigbreak

\begin{figure}[htb]
     \begin{center}
     \begin{tabular}{c}
     \includegraphics[width=10cm]{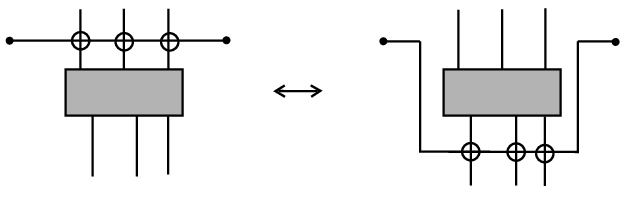}
     \end{tabular}
     \caption{\bf Detour Move}
     \label{Figure 2}
\end{center}
\end{figure}

Moves on virtual diagrams generalize the Reidemeister moves for classical knot and link
diagrams.  See Figure~\ref{Figure 1}.  One can summarize the moves on virtual diagrams by saying that the classical crossings interact with
one another according to the usual Reidemeister moves while virtual crossings are artifacts of the attempt to draw the virtual structure in the plane. 
A segment of diagram consisting of a sequence of consecutive virtual crossings can be excised and a new connection made between the resulting
free ends. If the new connecting segment intersects the remaining diagram (transversally) then each new intersection is taken to be virtual.
Such an excision and reconnection is called a {\it detour move}.
Adding the global detour move to the Reidemeister moves completes the description of moves on virtual diagrams. In Figure~\ref{Figure 1} we illustrate a set of local
moves involving virtual crossings. The global detour move is
a consequence of  moves of type (B) and (C) in Figure~\ref{Figure 1}. The detour move is illustrated in Figure~\ref{Figure 2}.  Virtual knot and link diagrams that can be connected by a finite 
sequence of these moves are said to be {\it equivalent} or {\it virtually isotopic}.
\bigbreak

\begin{figure}[htb]
     \begin{center}
     \begin{tabular}{c}
     \includegraphics[width=10cm]{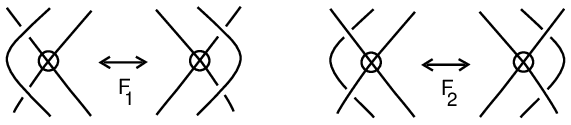}
     \end{tabular}
     \caption{\bf Forbidden Moves}
     \label{Figure 3}
\end{center}
\end{figure}

Another way to understand virtual diagrams is to regard them as representatives for oriented Gauss codes \cite{GPV}, \cite{VKT,DVK} 
(Gauss diagrams). Such codes do not always have planar realizations. An attempt to embed such a code in the plane
leads to the production of the virtual crossings. The detour move makes the particular choice of virtual crossings 
irrelevant. {\it Virtual isotopy is the same as the equivalence relation generated on the collection
of oriented Gauss codes by abstract Reidemeister moves on these codes.}  
\bigbreak

Figure~\ref{Figure 3} illustrates  two {\it forbidden moves}. Neither of these follows from Reidmeister moves plus detour move, and 
indeed it is not hard to construct examples of virtual knots that are non-trivial, but will become unknotted on the application of 
one or both of the forbidden moves. The forbidden moves change the structure of the Gauss code and, if desired, must be 
considered separately from the virtual knot theory proper. 
\bigbreak

\subsection{Interpretation of Virtuals Links as Stable Classes of Links in  Thickened Surfaces}
There is a useful topological interpretation \cite{VKT,DVK,Carter,Kamada3} for this virtual theory in terms of embeddings of links
in thickened surfaces.  Regard each 
virtual crossing as a shorthand for a detour of one of the arcs in the crossing through a 1-handle
that has been attached to the 2-sphere of the original diagram.  
By interpreting each virtual crossing in this way, we
obtain an embedding of a collection of circles into a thickened surface  $S_{g} \times R$ where $g$ is the 
number of virtual crossings in the original diagram $L$, $S_{g}$ is a compact oriented surface of genus $g$
and $R$ denotes the real line.  We say that two such surface embeddings are
{\em stably equivalent} if one can be obtained from another by isotopy in the thickened surfaces, 
homeomorphisms of the surfaces and the addition or subtraction of empty handles (i.e. the knot does not go through the handle).

\begin{figure}
     \begin{center}
     \begin{tabular}{c}
     \includegraphics[width=10cm]{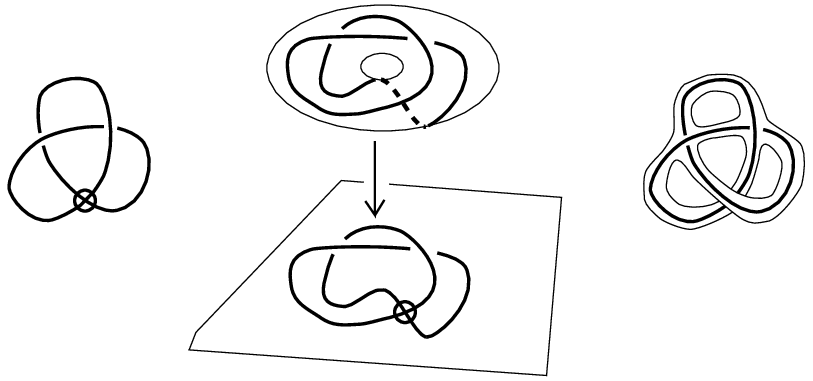}
     \end{tabular}
     \caption{\bf Surfaces and Virtuals}
     \label{Figure 4}
\end{center}
\end{figure}


\noindent We have the
\smallbreak
\noindent
{\bf Theorem 1 \cite{VKT,DKT,DVK,Carter}.} {\em Two virtual link diagrams are isotopic if and only if their corresponding 
surface embeddings are stably equivalent.}  
\smallbreak
\noindent
\bigbreak  

\noindent In Figure~\ref{Figure 4} we illustrate some points about this association of virtual diagrams and knot and link diagrams on surfaces.
Note the projection of the knot diagram on the torus to a diagram in the plane (in the center of the figure) has a virtual crossing in the 
planar diagram where two arcs that do not form a crossing in the thickened surface project to the same point in the plane. In this way, virtual 
crossings can be regarded as artifacts of projection. The same figure shows a virtual diagram on the left and an ``abstract knot diagram" \cite{Kamada3,Carter} on the right.
The abstract knot diagram is a realization of the knot on the left in a thickened surface with boundary and it is obtained by making a neighborhood of the 
virtual diagram that resolves the virtual crossing into arcs that travel on separate bands. The virtual crossing appears as an artifact of the
projection of this surface to the plane. The reader will find more information about this correspondence \cite{VKT,DKT} in other papers by the author and in
the literature of virtual knot theory.
\bigbreak
 
\subsection{Multiple Virtual Knot Theory}
We use the following notations:
ordinary virtual crossings denoted by $\CCircleDiag,$ and a second type of virtual crossing
denoted by $\CBoxDiag.$  
In principle, we shall have as many virtual crossings as we like and they can be labeled by glyphs like the above with 
a letter attached and the letter running over some chosen index set as in $\CCircleDiagL,$ or $\IMGlyphL.$\\

\noindent {\bf Detour Axiom for Multiple Virtual Crossings.} Each virtual crossing type detours over all the types of virtual crossings including itself, and each virtual crossing type detours over the classical crossings.
Figures  \ref{Figure 2}, \ref{EFF3} and  \ref{EFF4} illustrate the detour moves.\\ 

\noindent {\bf Definition.} {\it Multiple Virtual Knot Theory (MVKT)} is the diagrammatic theory generated by diagrams with multiple virtual crossings, using the Detour Axiom stated above in conjunction with the 
classical Reidemeister moves. Note that the virtual Reidemeister moves as shown in Figure~\ref{Figure 1} are consequences of the generalized detour move. The figure illustrates virtual Reidemeister moves for the case of one type of virtual and in the generalized theory there are such moves for each type of virtual crossing. For a given instantiation of MVKT, a set of virtual crossing types is chosen with a method for indicating them. Note that in one application of the detour move, only one virtual crossing type will be used.
Thus if one excises a sequence of virtual crossings of a given type, then the new connection will create only virtual crossings of that type. If the excision does not remove any virtual crossings, then the reconnection can use
any single virtual crossing type that is available.\\

Note the two component ``link" with two virtual crossings in Figure~\ref{EFF2}. This link is not equivalent to a disjoint union of two circles. We will see a proof of this below. The figure also illustrates the generalized bracket expansion of a multi-virtual knot with one classical crossing and two virtual crossings. We will define the generalized bracket below. Figure~\ref{onecomp} shows a two-virtual diagram of one component and a series of detour moves that reduce it to a circle.{\it We conjecture that any single component diagram, decorated only with two virtual crossings, reduces by detour moves to a circle.}\\

The author has previously considered multiple virtual crossings with hierarchies of detour moves among them and there is a paper \cite{MFlat} taking steps with this idea.
In multiple virtual theory, all the virtual crossings are on the same footing. They can detour over one another, but classical crossings cannot detour over virtual crossings.\\

 Many questions arise about the relationship of classical knot theory, virtual knot theory and multiple virtual knot theory. It is known that classical knot theory embeds in virtual knot theory in the sense that if two classcial knot or link diagrams are equivalent as virtual diagrams, then they are equivalent as classical diagrams \cite{VKT}. The same result holds for multiple virtual knot theory.\\
 
 \noindent {\bf Theorem.} Let $K$ and $L$ be classical link diagrams in classical knot theory (KT), and let MVKT denote a multiple virtual knot theory (with some given cardinality of virtual crossing types). Since $K$ and $L$ can be viewed as diagrams in MVKT we have a tautological mapping $F: KT \longrightarrow MVKT $ where F(K) is the same diagram seen in the larger category of virtual moves. Then if $F(K)$ and $F(L)$ are equivalent in MVKT, they are equivalent in KT. Thus classical knot theory embeds in multiple virtual knot theory.\\

\noindent {\bf Proof.} The proof has the same structure as the proof for virtual knot theory in \cite{VKT}. It is known from the work of Waldhausen that the fundamental group and peripheral subgroup as a pair $(G,P)$ classify classical knots and links. Thus if $K$ and $L$ are classical links with peripheral subgroup pairs $(G(K), P(K))$ and  $(G(L), P(L))$ then $K$ is classically ambient istopic to $L$ if and only if the corresponding group pairs are isomorphic. The group pair of a classical link or a virtual link diagram $K$ can be given by a Wirtinger presentation read from the diagram in conjunction with longitudes $\lambda_{i}(K)$ (one for each link component) that is obtained by walking along a component and writing the product of the Wirtinger generators 
(corrsponding to oriented arcs of the diagram) that are underpassed in the course of the walk. These longitudes are not affected by detour moves. From this it follows that if $L$ is obtained from $K$ by a sequence of Reidemeister moves and virtual detour moves (in the sense of MVKT), then  $(G(L), P(L))$  is isomorphic as a pair with $(G(K), P(K)).$ If $K$ and $L$ are classical diagrams, it follows that they are classically equivalent. This completes the proof of the Theorem.  $\hfill\Box$

\begin{figure}
     \begin{center}
     \begin{tabular}{c}
     \includegraphics[width=10cm]{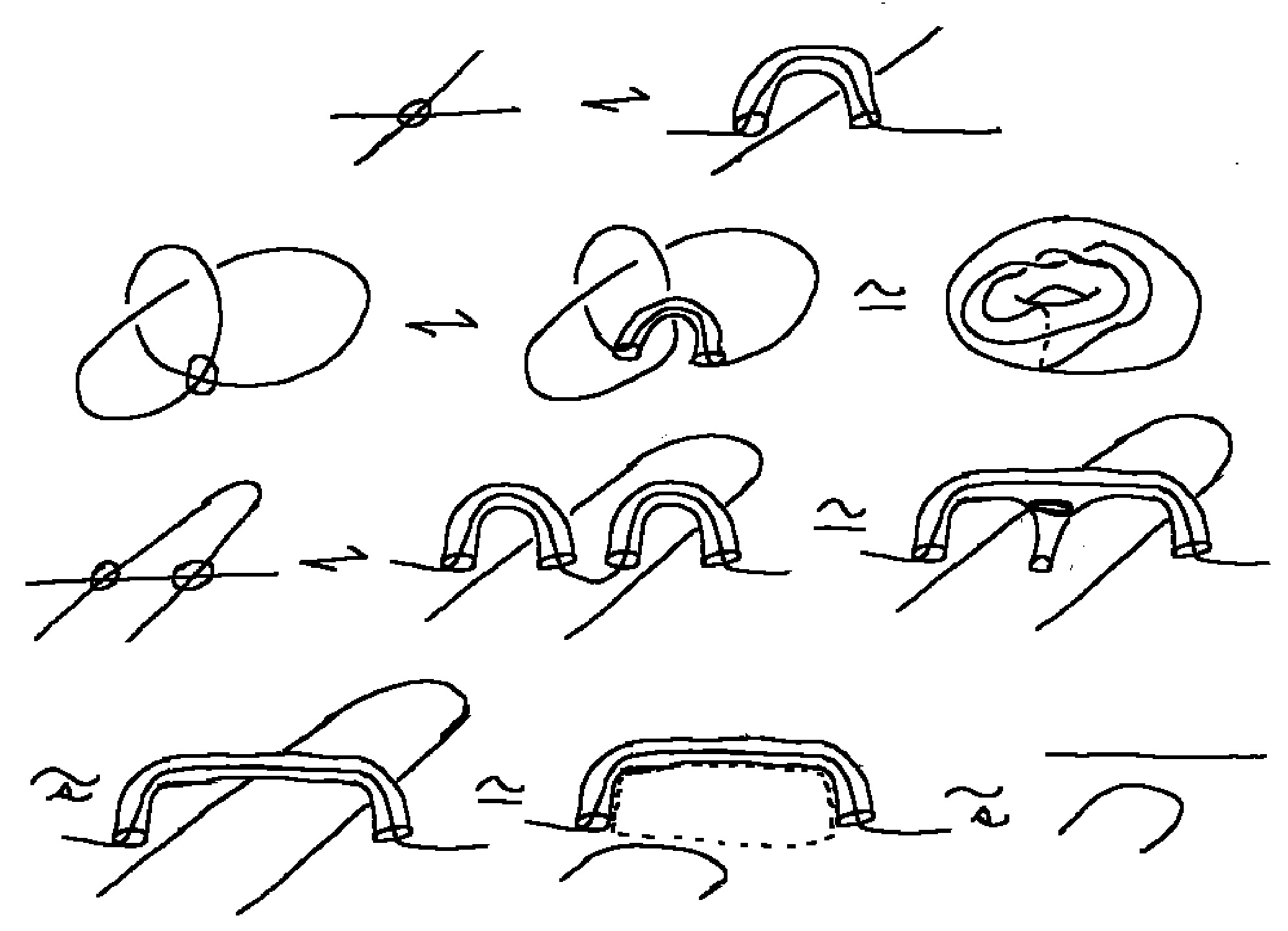}
     \end{tabular}
     \caption{\bf Replacing Virtual Crossings by Handle Detours}
     \label{Handles}
\end{center}
\end{figure}

\subsection{Surface Interpretations for Multiple Virtual Knot Theory}
As we mentioned above, single virtual knot theory can be interpreted in terms of knots in thickened surfaces. We can obtain a knot or link in a surface by using an abstract link diagram, obtained by producing a ribbon neighborhood of the virtual diagram and using two local ribbons for each virtual crossing as shown in Figure~\ref{Figure 4}. Another method that associates a surface to a virtual diagram is illustrated in Figure~\ref{Handles}. In that figure we add a handle to the two-sphere or plane for the virtual link at each virtual crossing. The virtual crossing is replaced by two arcs, one going through the handle, and one going underneath the handle. This construction can be generalized to multiple virtual knot theory by labelling the resulting handles according to the label of the virtual crossing. The handles then interact as shown in Figure ~\ref{Handles} only when they have the same label. In the figure we have illustrated how a detour move (virtual Reidemeister move two) is seen in terms of handle stabilization. The two handles interact and undergo a 1-handle stabilization that fuses them into a single handle. This allows the arc passage for the move and then we destabilize the resulting handle, taking us back to the two sphere or the plane. Handles that interact come from the same type of virtual crossing.  This is, of course, a mixed combinatorial and topological interpretation.\\

\subsection{Generalized Bracket Polynomial}
\begin{figure}
     \begin{center}
     \begin{tabular}{c}
     \includegraphics[width=8cm]{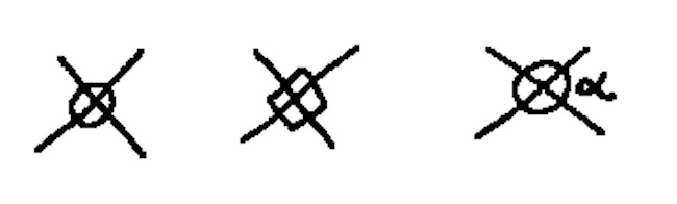}
     \end{tabular}
     \caption{\bf Virtual Crossing Notations}
     \label{EFF1}
\end{center}
\end{figure}

We now start the proper discussion of multiple virtual knots and links. In Figure~\ref{EFF1} we illustrate graphical notations for virtual crossings.
The traditional circle notation is shown along with a local box notation for a distinct virtual crossing, and a circle notation with an index $\alpha.$
The $\alpha$ can run over any convenient index set, allowing any chosen cardinality of virtual crossings. \\

\begin{figure}
     \begin{center}
     \begin{tabular}{c}
     \includegraphics[width=10cm]{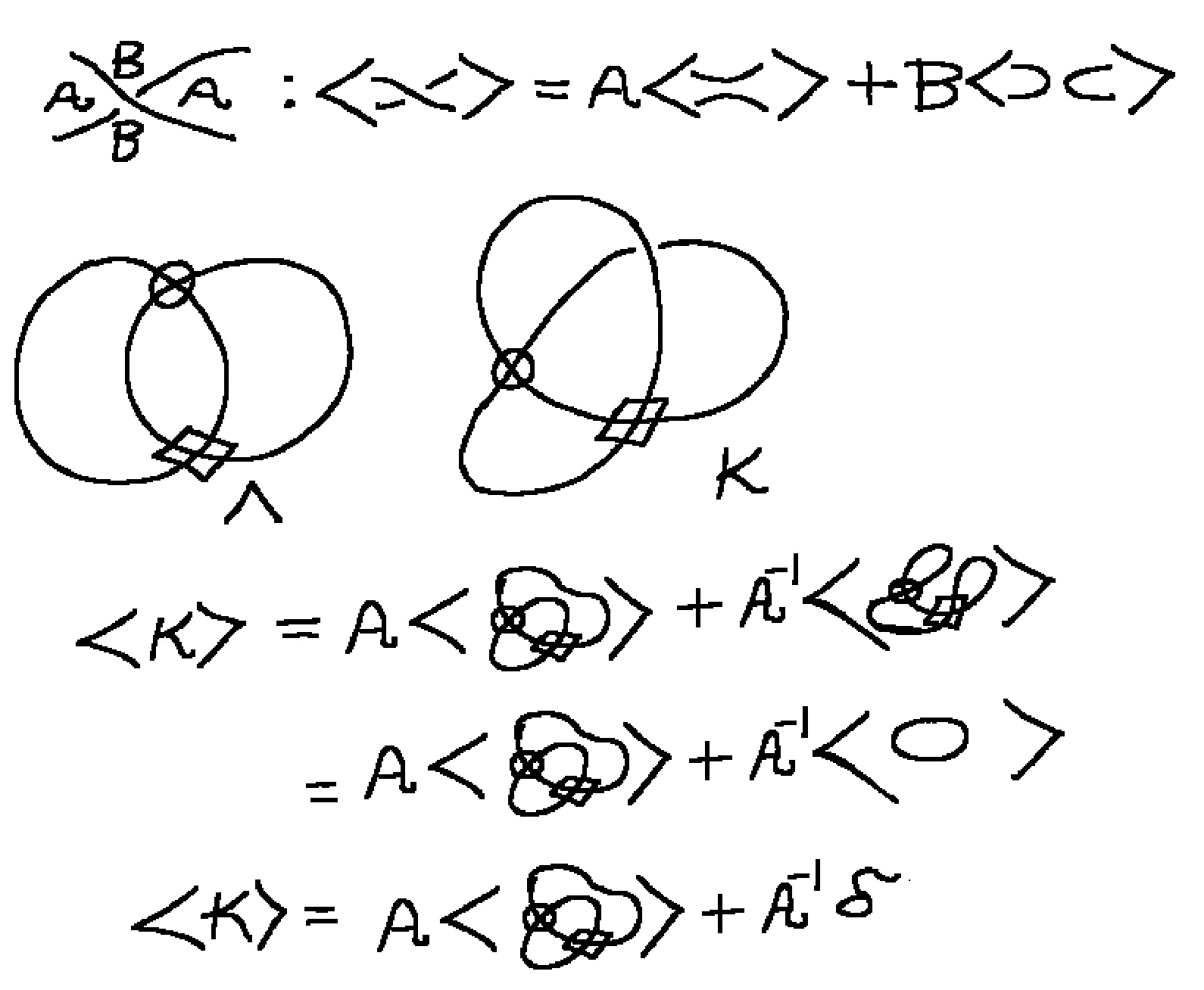}
     \end{tabular}
     \caption{\bf Double Virtual Link and Double Virtual Knot}
     \label{EFF2}
\end{center}
\end{figure}

\begin{figure}
     \begin{center}
     \begin{tabular}{c}
     \includegraphics[width=10cm]{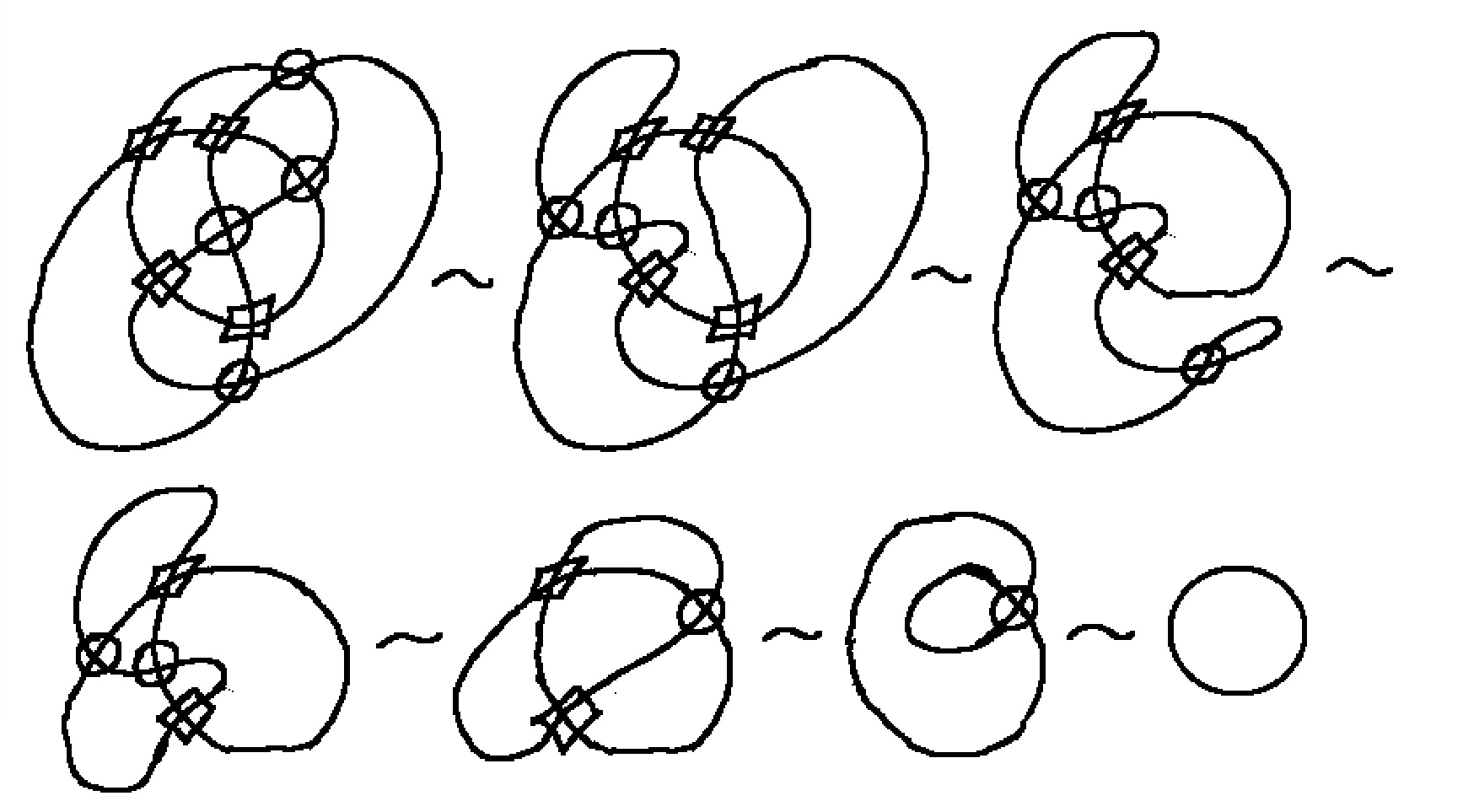}
     \end{tabular}
     \caption{\bf Single Component with Two Virtual Crossings}
     \label{onecomp}
\end{center}
\end{figure}

In Figure~\ref{EFF2} we show the simplest multiple virtual link, a link with two virtual crossings, one of circular type and one of box type.
Our axioms for virtual crossings do not permit two different types to cancel, and  this link can be proved to be non-trivial. We will shortly explain an algebraic invariant that verifies this fact.\\

\begin{figure}
     \begin{center}
     \begin{tabular}{c}
     \includegraphics[width=8cm]{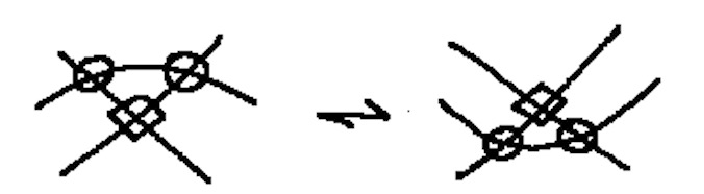}
     \end{tabular}
     \caption{\bf Detour Move for Virtuals across a Virtual}
     \label{EFF3}
\end{center}
\end{figure}

\begin{figure}
     \begin{center}
     \begin{tabular}{c}
     \includegraphics[width=8cm]{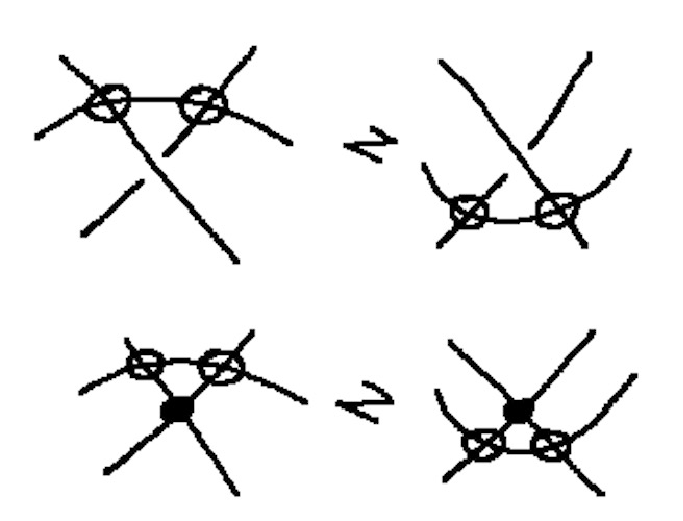}
     \end{tabular}
     \caption{\bf Detour Moves}
     \label{EFF4}
\end{center}
\end{figure}

In Figure~\ref{EFF3} and Figure~\ref{EFF4} we illustrate a simple example of a detour move of two virtuals of one type, across a virtual of another type. In general the detour move for any given virtual type allows a consecutive sequence of virtuals of a given type to be excised, and a new arc can be drawn using these resulting endpoints. This arc, where it crosses other parts of the diagram will meet these parts in the same virtual type that was excised. In Figure~\ref{EFF3}  we show a detour move in relation to another virtual crossing. In Figure~\ref{EFF4}  we show a detour move in relation to a classical crossing and in relation to a 4-valent graphical node. We will use such nodes in our constructions for invariants.\\

\begin{figure}
     \begin{center}
     \begin{tabular}{c}
     \includegraphics[width=8cm]{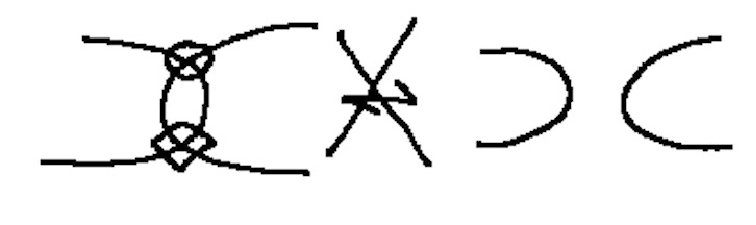}
     \end{tabular}
     \caption{\bf Non Cancellation of Distinct Virtual Crossings.}
     \label{EFF5}
\end{center}
\end{figure}

Figure~\ref{EFF5} illustrates the non-cancellation of distinct virtual crossings.\\

We first define a generalized bracket polynomial for any multiple virtual theory. This bracket invariant has the usual skein relation and loop evaluation.
$$\langle \CrossGlyph \rangle = A \langle \HSmoothGlyph \rangle + A^{-1} \langle \VSmoothGlyph \rangle$$
$$\langle O \rangle = -A^2 -A^{-2}= \delta$$\\

Applying the expansion formula recursively will yield collections of loops with multiple virtual crossings. The generalized bracket has state expansion 
$$<K> = \Sigma_{S} A^{a(S)-b(S)} \langle S \rangle$$ 
where $S$ is a state of the diagram $K$ obtained by choosing a smoothing for each 
crossing of $K$ and $a(S)$ is the number of smoothings of type $A,$ while $b(S)$ is the number of smoothings of type $B.$ Smoothing types correspond to the convention shown in 
Figure~\ref{EFF2} where the local regions of a crossing are labeled $A$ and $B$ so that the $A$ regions are swept by a counterclockwise turn of the overcrossing arc, and a smoothing is of type $A$ if the two $A$ regions are joined by the smoothing. Each state $S$ in the virtual bracket expansion is a diagram with multiple virtual crossings. For the generalized bracket, we take $\langle S \rangle$ to be the multiple virtual class (i.e. the equivalence class of the state up to virtual detour moves)  of the state $S$ with the caveat that any free disjoint circle in $S$ (or in a diagram equivalent to $S$) is evaluated as $\delta.$ Thus we have that $\langle S O \rangle = \delta \langle S \rangle.$\\

Since, in general, we do not know full information about the class of a state $S$ up to detour moves, this generalized bracket is a {\it picture valued} invariant where the pictures (diagrams) of the states represent their 
equivalence classes. There are different ways to create invariants for these states, and so it is important to leave the generalized bracket in this open-ended form.\\

In Figure~\ref{EFF2} we illustrate the bracket expansion of a multiple virtual knot $K.$ In this expansion we obtain one picture $\langle \Lambda \rangle$ where $\Lambda$ is the state indicated in the figure, consisting of two circles intersecting in two distinct virtual crossings.\\

Determining the virtual detour class of these states in the multi-virtual setting is an additional mathematical problem. To this end, we make a special evaluation in the case of {\it double multiple virtual theory} where we use two virtual crossing types, one indicated by a circle and the other by a box as in Figure~\ref{EFF1}.
Then we take $\langle S \rangle$ to be the graphical evaluation previously discussed in Section 3 and using the expansion in Figure~\ref{EFF7}.
We call this specialized bracket the {\it chromatic bracket} since its extra rule is motivated by the coloring problems for graphs. The chromatic bracket is then governed by the rules shown below.\\

$$\langle \CrossGlyph \rangle = A \langle \HSmoothGlyph \rangle + B \langle \VSmoothGlyph \rangle$$
$$\langle O \rangle = \delta$$
$$\langle \IMGlyph \rangle  = 2 \langle \CDotDiag \rangle  - \langle \VirtualGlyph \rangle $$\\

If we take $A=1$ and $B=-1$ and let $\delta = n$, a natural number, then the chromatic bracket counts colorings of a knot or link diagram as we have already described at the end of 
Section 3, items 6 and 7. That is, we have from the previous definition that $$\langle L \rangle(A=1, B=-1, \delta = n] = PK[L].$$ This, in turn is a translation of the generalized coloring evaluation of a graph with a given perfect matching. If the graph is non-planar and represented with box-type virtual crossings in the plane, then the associated virtual link diagram is the input for the Chromatic Bracket with 
$A=1,B=-1,\delta = n.$\\

If we take $B=A^{-1}$ and $\delta = -A^2 - A^{-2},$ then the Chromatic Bracket is invariant under the second and third Reidemeister moves and, by dint of the box virtual evaluation formula
(used exactly as we did in the coloring evaluations) it is a Laurent polynomial valued invariant. Each picture has been translated into an invariant polynomial.\\

We will use the same notation $\langle K \rangle$ for both the general bracket and for the chromatic bracket evaluated at $B = A^{-1}$ and $\delta = -A^2 - A^{-2}$. The general bracket takes its values in sums of Laurent polynomials in $A$ and $A^{-1}$ as coefficients of graphical equivalence classes of purely virtual diagrams. There will be cases where the chromatic method is not sufficient and we will examine these equivalence classes directly or by using another way to evaluate them. This will be detailed below.\\

\begin{figure}
     \begin{center}
     \begin{tabular}{c}
     \includegraphics[width=8cm]{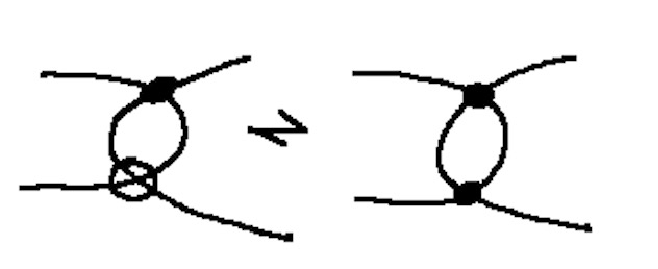}
     \end{tabular}
     \caption{\bf Nodal Equivalence}
     \label{EFF6}
\end{center}
\end{figure}

\begin{figure}
     \begin{center}
     \begin{tabular}{c}
     \includegraphics[width=8cm]{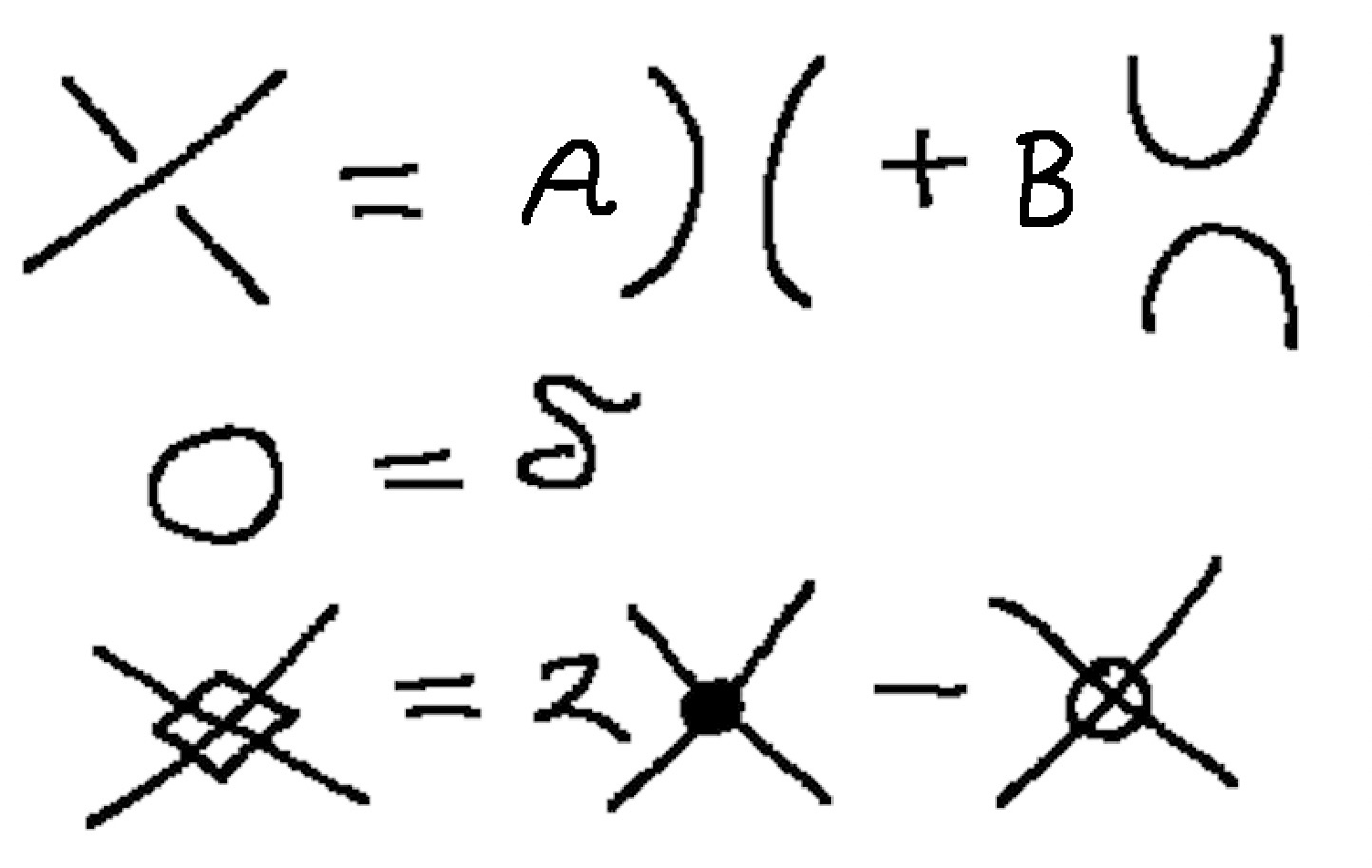}
     \end{tabular}
     \caption{\bf Chromatic Bracket}
     \label{EFF7}
\end{center}
\end{figure}

Figure~\ref{EFF6} illustrates our rule for the nodal crossing. This is in accord with the usage for the nodal crossing that was explained in Section 2.3 in relation to evaluating the 
generalized Penrose Polynomial. We will reenter this formalism now in a topological context. For this purpose, the reader should view Figure~\ref{EFF7}. In that Figure we show the expansion rules for the chromatic bracket polynomial. Here we use coefficients $A$ and $B$ for the two smoothings and a value of $\delta$ for the loop value.\\

The chromatic bracket polynomial is defined for a theory with two virtual crossings that we will indicate with round and box indicators at the 
corresponding crossing. The polynomial has the usual bracket expansion at a crossing coupled with our expansion for the box virtual crossings as we have used them for the graph theory earlier in the paper. In the formulas below we use the topological coefficients for the bracket with 
$B = A^{-1}$ and $\delta = -A^2 - A^{-2}.$ The resulting evaluation is invariant under all moves except the first Reidemeister move where it changes by a multiple of $(-A^3)$ or $(-A^{-3})$ for positive and negative curls respectively.\\

\begin{figure}
     \begin{center}
     \begin{tabular}{c}
     \includegraphics[width=8cm]{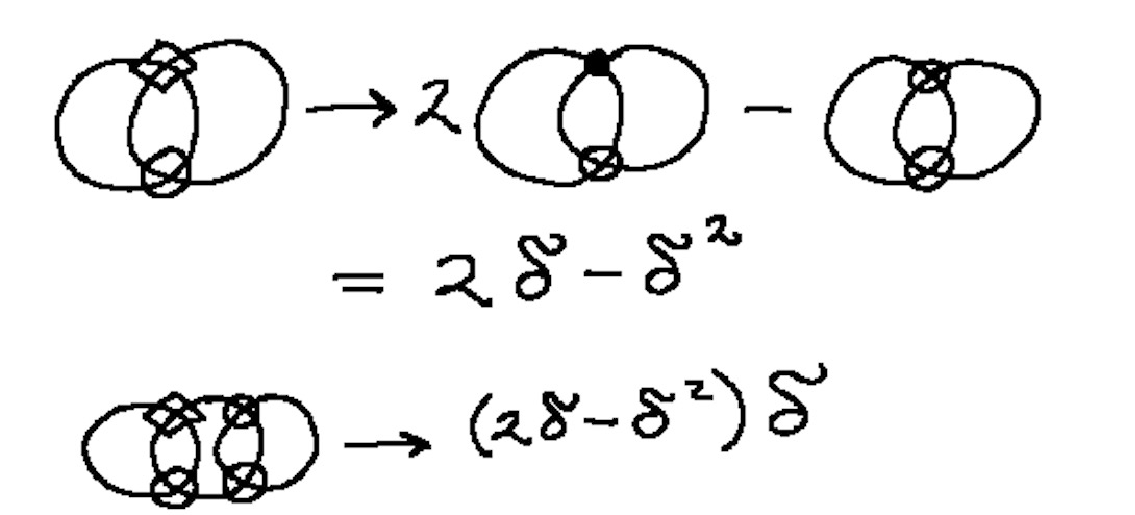}
     \end{tabular}
     \caption{\bf Loop Evaluations}
     \label{EFF8}
\end{center}
\end{figure}

The context for the chromatic  bracket is multi-virtuals with two crossing types: circle and box and the rule that the box virtual crossing is expanded into twice the nodal crossing minus the circle virtual. This is the identical formalism to our coloring polynomial from the Section 3. Hence one expands the bracket until there are only collections of virtual loops.
Then the box virtuals are expanded, and one has virtual loops with nodes. A collection of loops that is connected by nodes is evaluated as a single $\delta$. Thus if $S$ is a state in the expansion of this bracket then the evaluation of $S$, denoted $[S] = \delta^{c(S)}$ where $c(S)$ denotes the number of connected components of $S$ relative to the 
nodes. We illustrate this calculation on the virtual diagram in Figure~\ref{EFF8}.\\

Note that these state evaluations $\langle S \rangle$ can be reformulated in terms of dichromatic polynomial evaluations and alternating sums of chromatic polynomial evaluations of graphs $G(S),$ as we have done in the end of Section 3, items 4 and 5 and the figures: Figure~\ref{EFF11}, Figure~\ref{Taut}, Figure~\ref{Taut1}. Translation of the states $S$ to the graphs $G(S)$ is useful for doing these calculations.\\

We have the basic theorem:\\

\noindent {\bf Theorem.} The generalized bracket is invariant under detour moves for any choice of commuting variables $A,B, \delta.$
The generalized and chromatic brackets $<K>$ are invariant under the Reidemeister moves $R2$ and $R3$ when we take $B= A^{-1}, \delta = -A^2 - A^{-2}.$
We call this the {\it topological specialization} of the bracket. For arbitrary multiple virtuals the topological specialization of the generalized bracket is invariant under the R2 and R3 moves and takes values in the ring generated by the virtual classes of the states of the diagram with coefficients in the Laurent polynomials in $A.$\\

\noindent {\bf Proof.} It is only necessary to check that the expansion of the box virtual is compatible with the detour moves. This is verified diagrammatically in Figure~\ref{EFF9} and Figure~\ref{EFF10}. This completes the proof. $\hfill\Box$


\clearpage

In Figure~\ref{EFF8} we illustrate double virtual calculations for the simple link using two virtual crossings, and another example with three components.\\

\begin{figure}
     \begin{center}
     \begin{tabular}{c}
     \includegraphics[width=10cm]{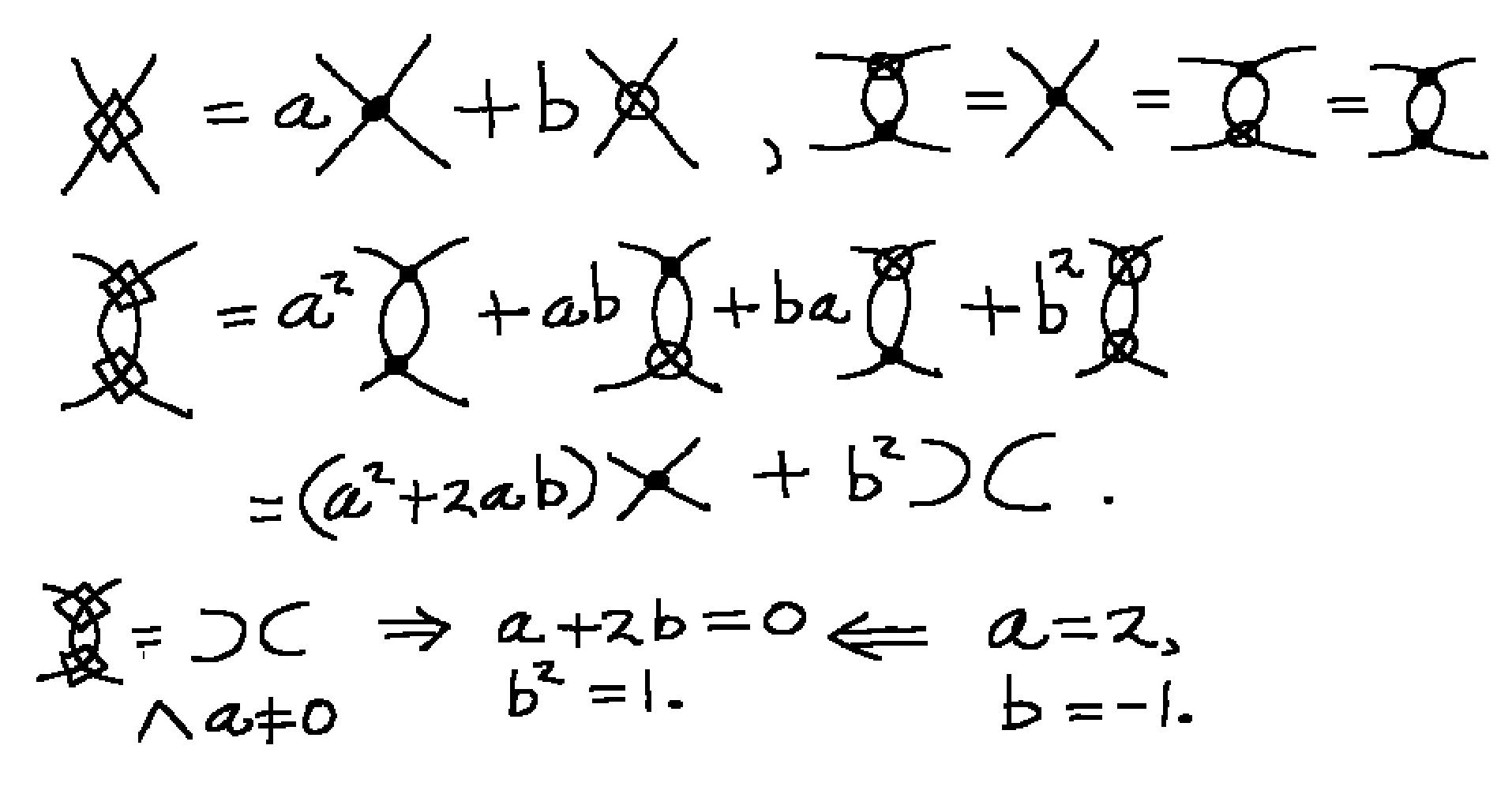}
     \end{tabular}
     \caption{\bf Graphical Expansion Verification}
     \label{Coeffs}
\end{center}
\end{figure}

\begin{figure}
     \begin{center}
     \begin{tabular}{c}
     \includegraphics[width=10cm]{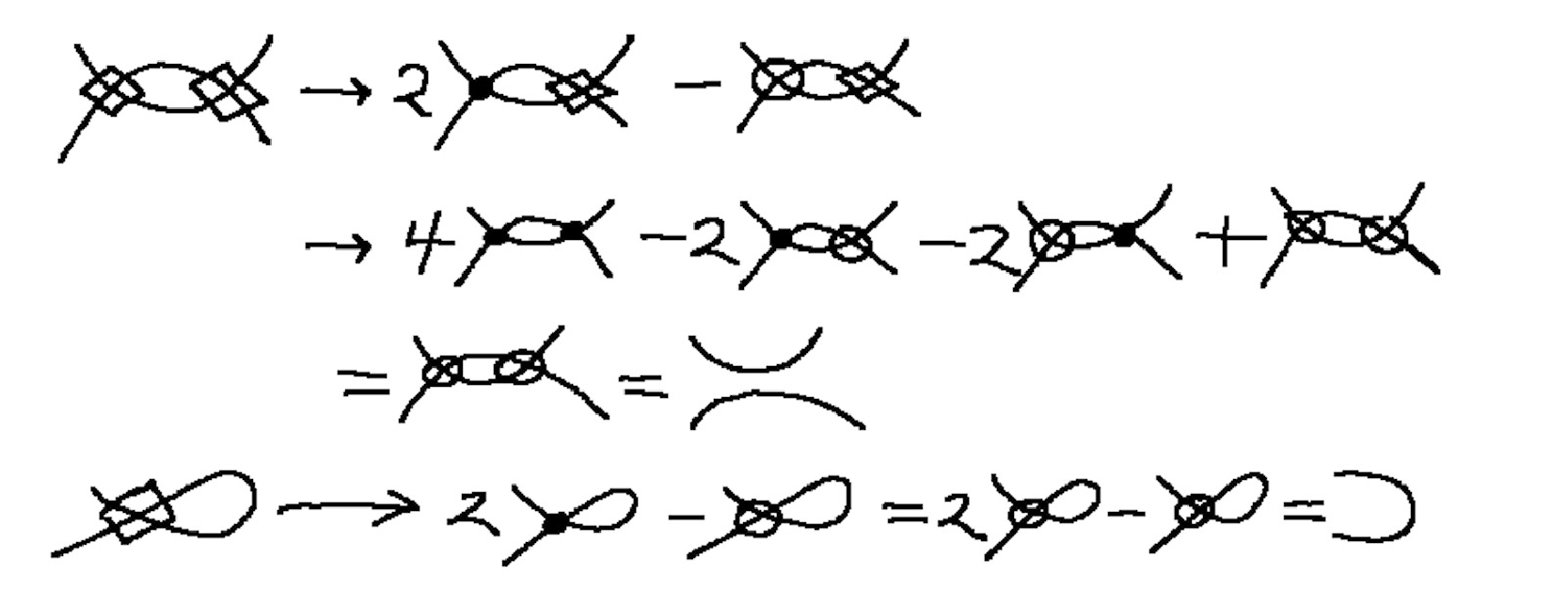}
     \end{tabular}
     \caption{\bf Graphical Expansion Verification}
     \label{EFF9}
\end{center}
\end{figure}

\begin{figure}
     \begin{center}
     \begin{tabular}{c}
     \includegraphics[width=10cm]{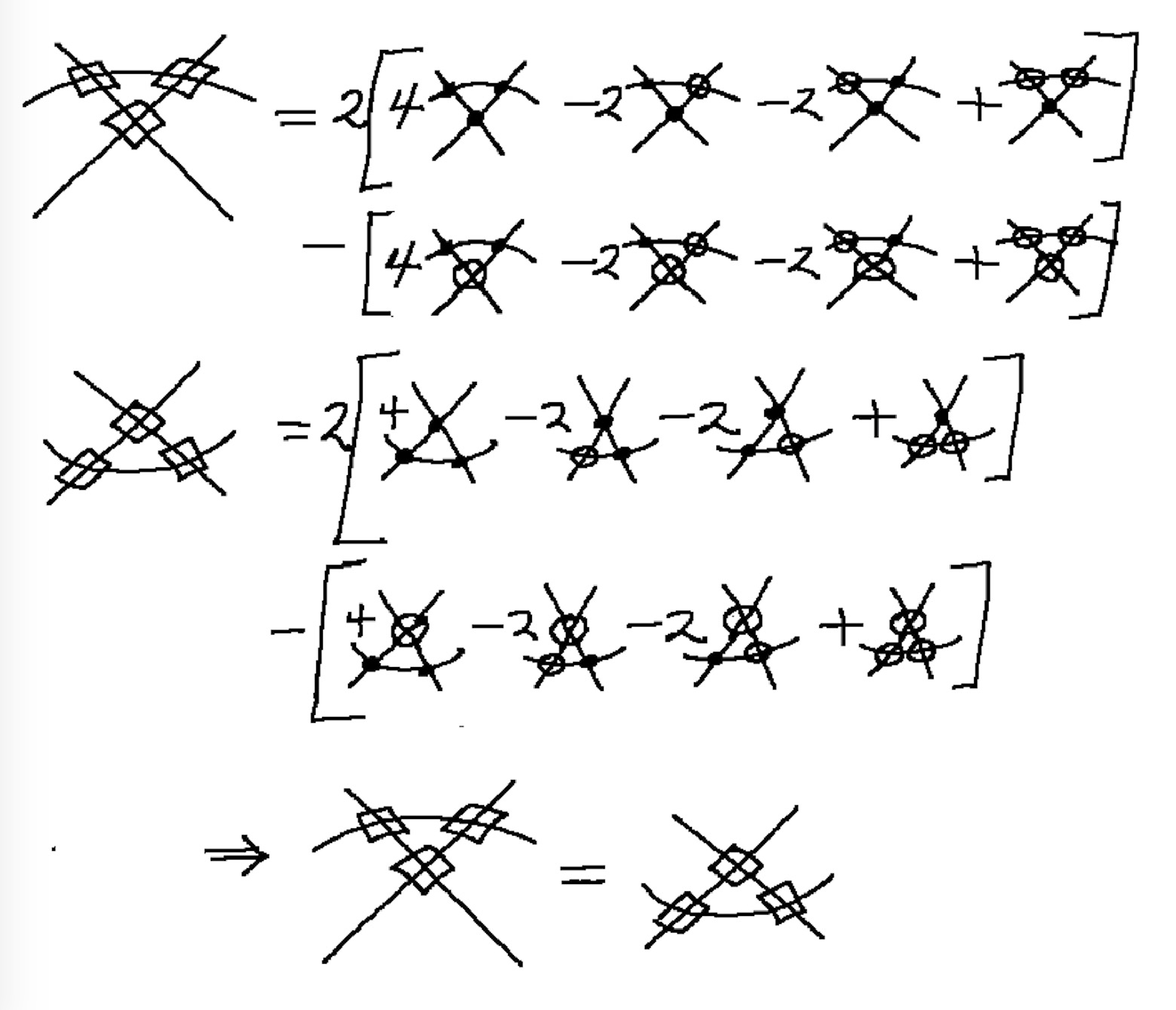}
     \end{tabular}
     \caption{\bf Graphical Expansion Verification}
     \label{EFF10}
\end{center}
\end{figure}

Figure~\ref{Coeffs} shows that the expansion for the box virtual crossing necessarily has the coefficients $2$ and $-1$ if it is to satisfy the simplest virtual detour move.
In Figure~\ref{EFF9} and Figure~\ref{EFF10} we illustrate the graphical calculations that verify that the expansion formula for the box virtual crossing is compatible with 
the detour moves and virtual Reidemeister three moves (also detour moves) for the doubled virtual crossings.  These figures constitute the crucial lemma that this method of evaluation leads to a regular isotopy invariant for doubled virtual links.\\

\clearpage

\begin{figure}
     \begin{center}
     \begin{tabular}{c}
     \includegraphics[width=10cm]{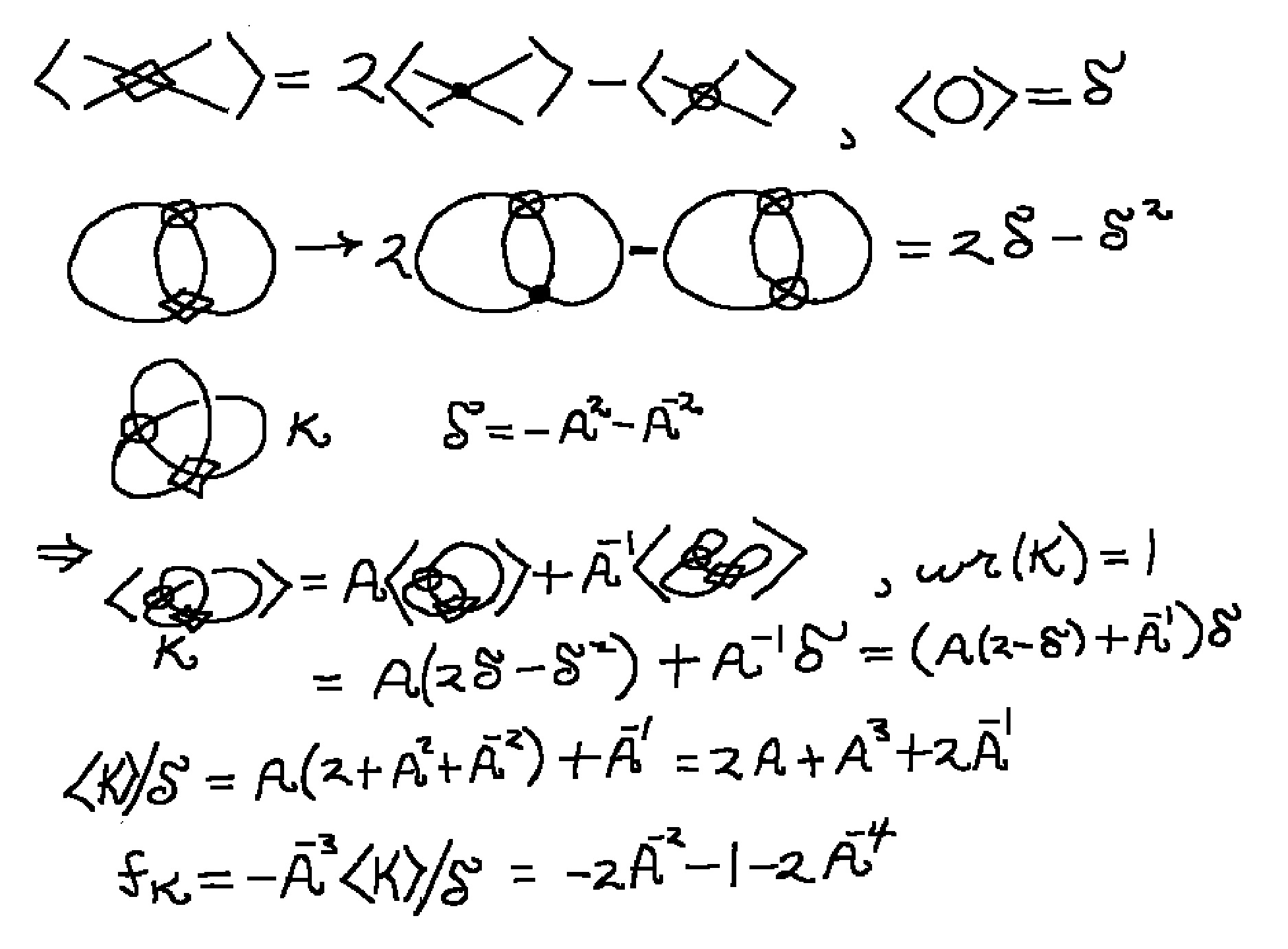}
     \end{tabular}
     \caption{\bf A Chromatic Bracket Evaluation}
     \label{EFF12}
\end{center}
\end{figure}

\begin{figure}
     \begin{center}
     \begin{tabular}{c}
     \includegraphics[width=8cm]{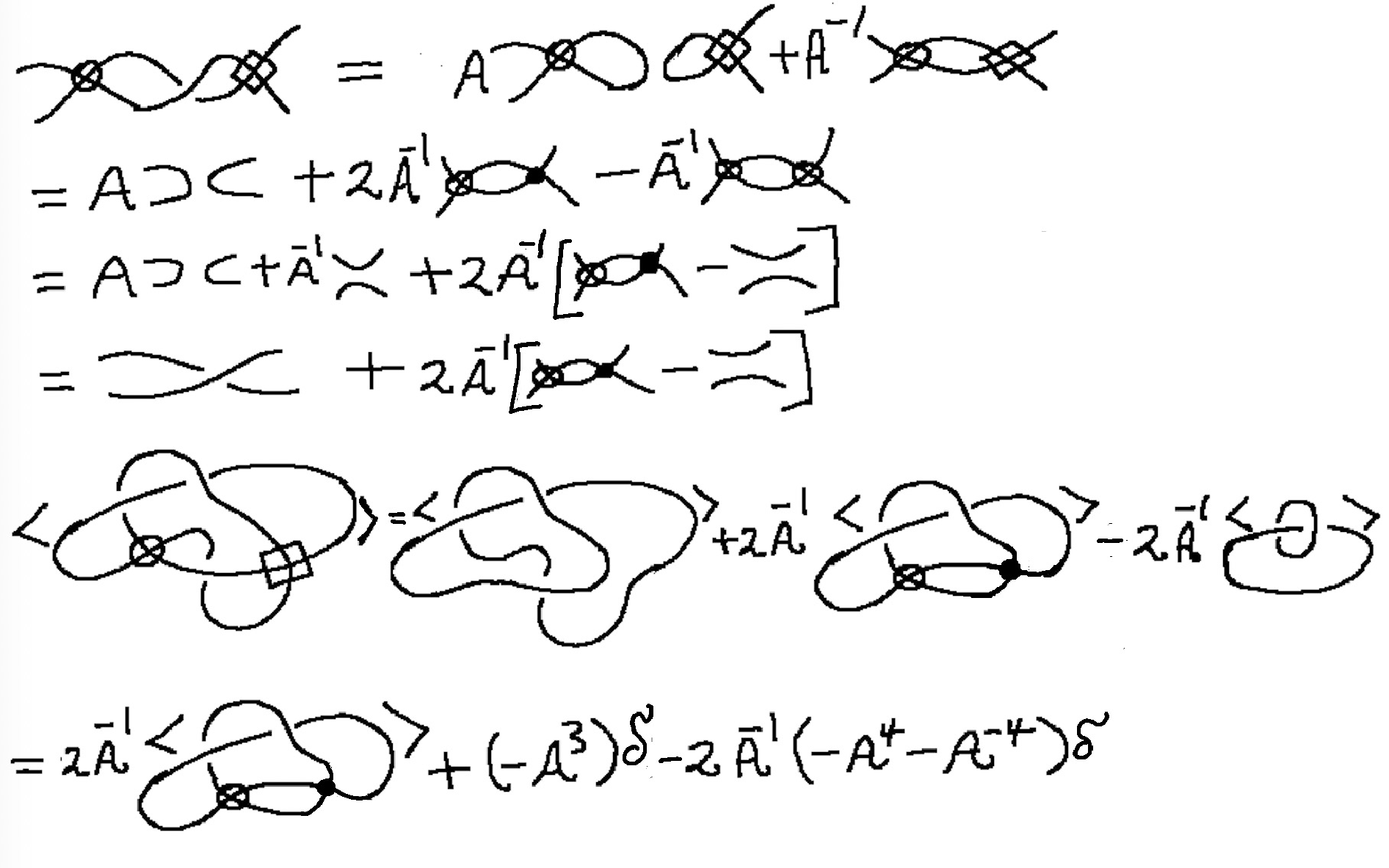}
     \end{tabular}
     \caption{\bf Bracket Expansion, Part 1}
     \label{EFF13}
\end{center}
\end{figure}

\begin{figure}
     \begin{center}
     \begin{tabular}{c}
     \includegraphics[width=8cm]{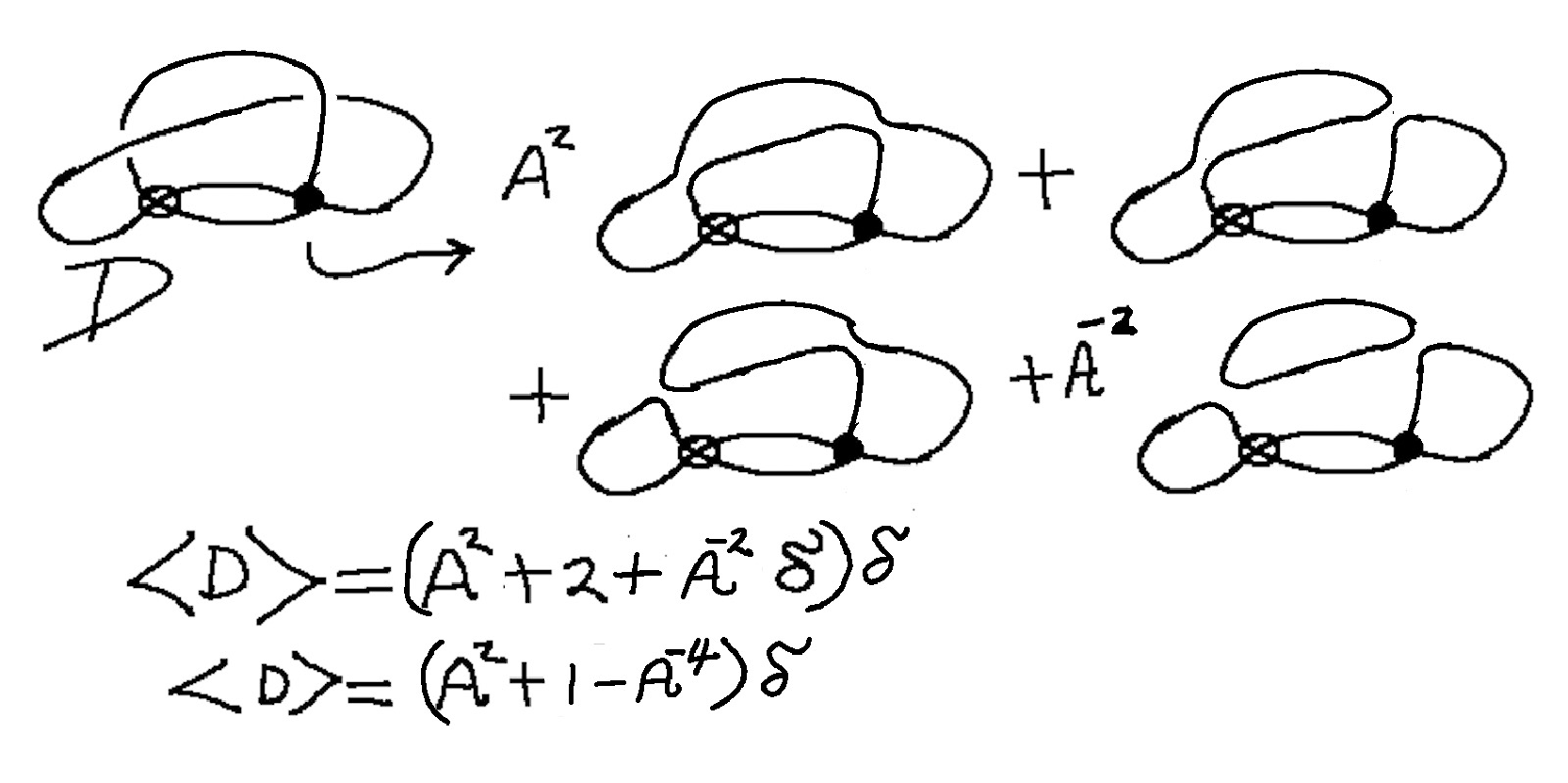}
     \end{tabular}
     \caption{\bf Bracket Expansion, Part 2}
     \label{EFF14}
\end{center}
\end{figure}

\begin{figure}
     \begin{center}
     \begin{tabular}{c}
     \includegraphics[width=8cm]{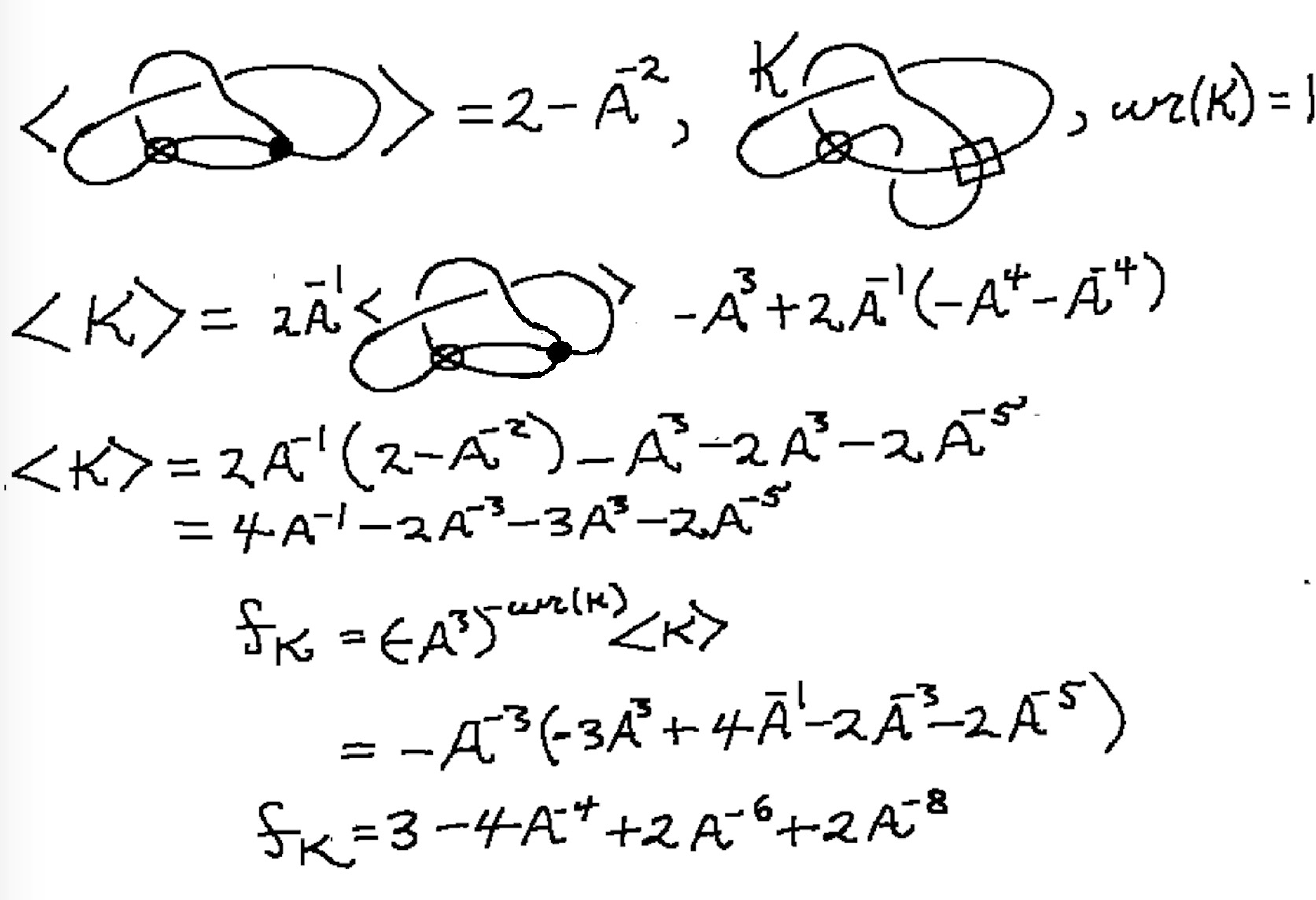}
     \end{tabular}
     \caption{\bf Bracket Expansion, Part 3}
     \label{EFF15}
\end{center}
\end{figure}

Figure~\ref{EFF12} shows the expansion of the generalized bracket for double virtuals on a trefoil diagram with one classical crossing, one box virtual crossing and one round virtual crossing.\\

Figures ~\ref{EFF13}~\ref{EFF14}~\ref{EFF15} illustrate the calculation of the generalized bracket for a modified trefoil diagram $K$ with one box and one round virtual crossing.
In this case the bracket polynomial is non-trivial, while the case with two round or two box virtual crossings gives a trivial bracket polynomial with value $1$ after normalization.
In Figure~\ref{EFF13} we first show the result of expanding a virtual tangle consisting of round and box virtual crossings flanking a classical crossing. We expand the box crossing into rigid node and round virtual crossing as described above. If both virtual crossings had been of the same type, this expansion would be the same as expanding just an unflanked classical crossing. Here we see that the fact that different types of virtual crossing do not cancel at a Reidemeister two configuration creates more structure. In the lower half of the figure we use this expansion to evaluate the diagram $K$ up to a single diagram with one rigid node and one round virtual crossing. In Figure~\ref{EFF14}  and 
Figure~\ref{EFF15} we evaluate this noded diagram. These two figures complete the details of the calculation.\\

\begin{figure}
     \begin{center}
     \begin{tabular}{c}
     \includegraphics[width=8cm]{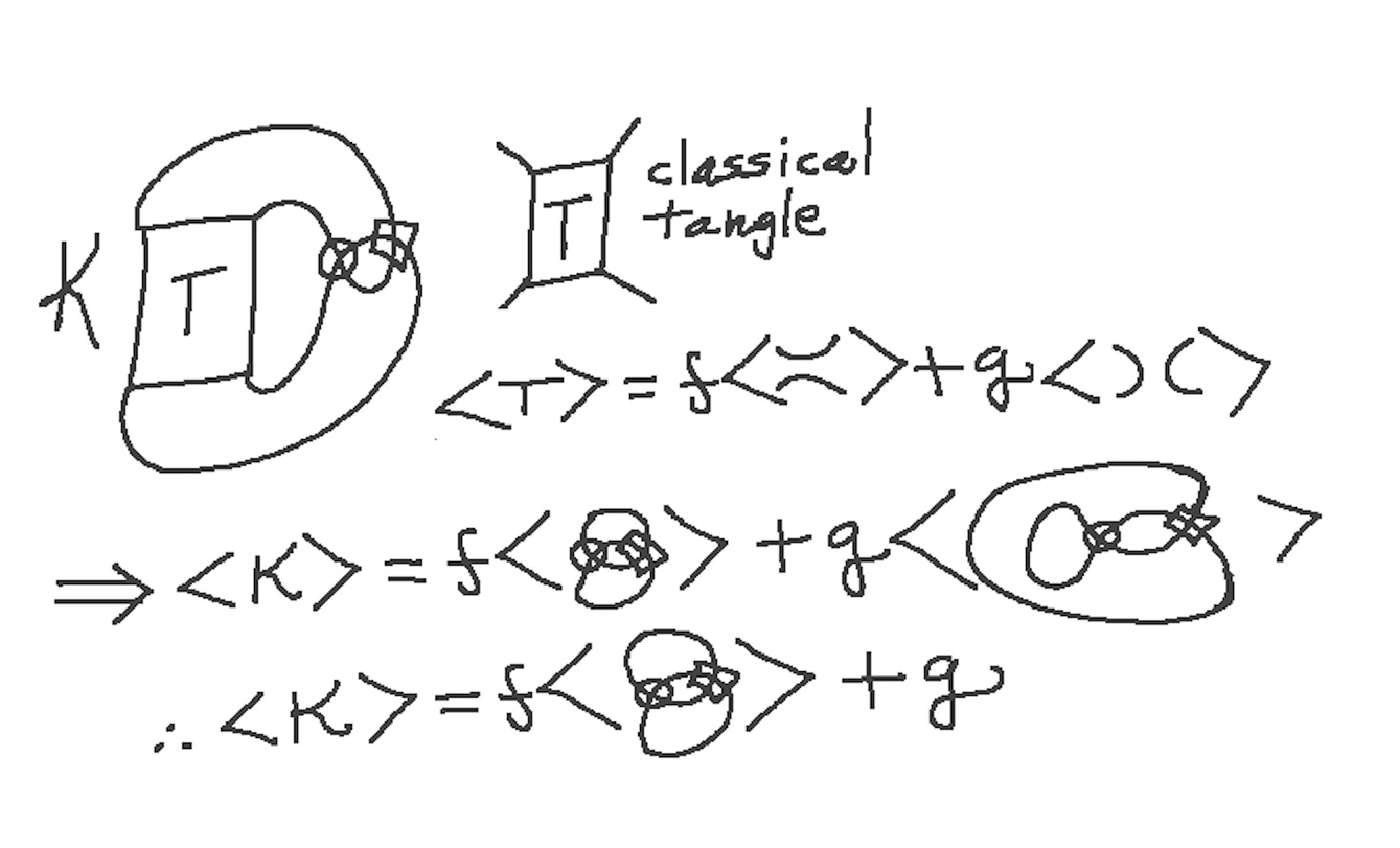}
     \end{tabular}
     \caption{\bf Tangle Expansion}
     \label{EFF27}
\end{center}
\end{figure}

\begin{figure}
     \begin{center}
     \begin{tabular}{c}
     \includegraphics[width=8cm]{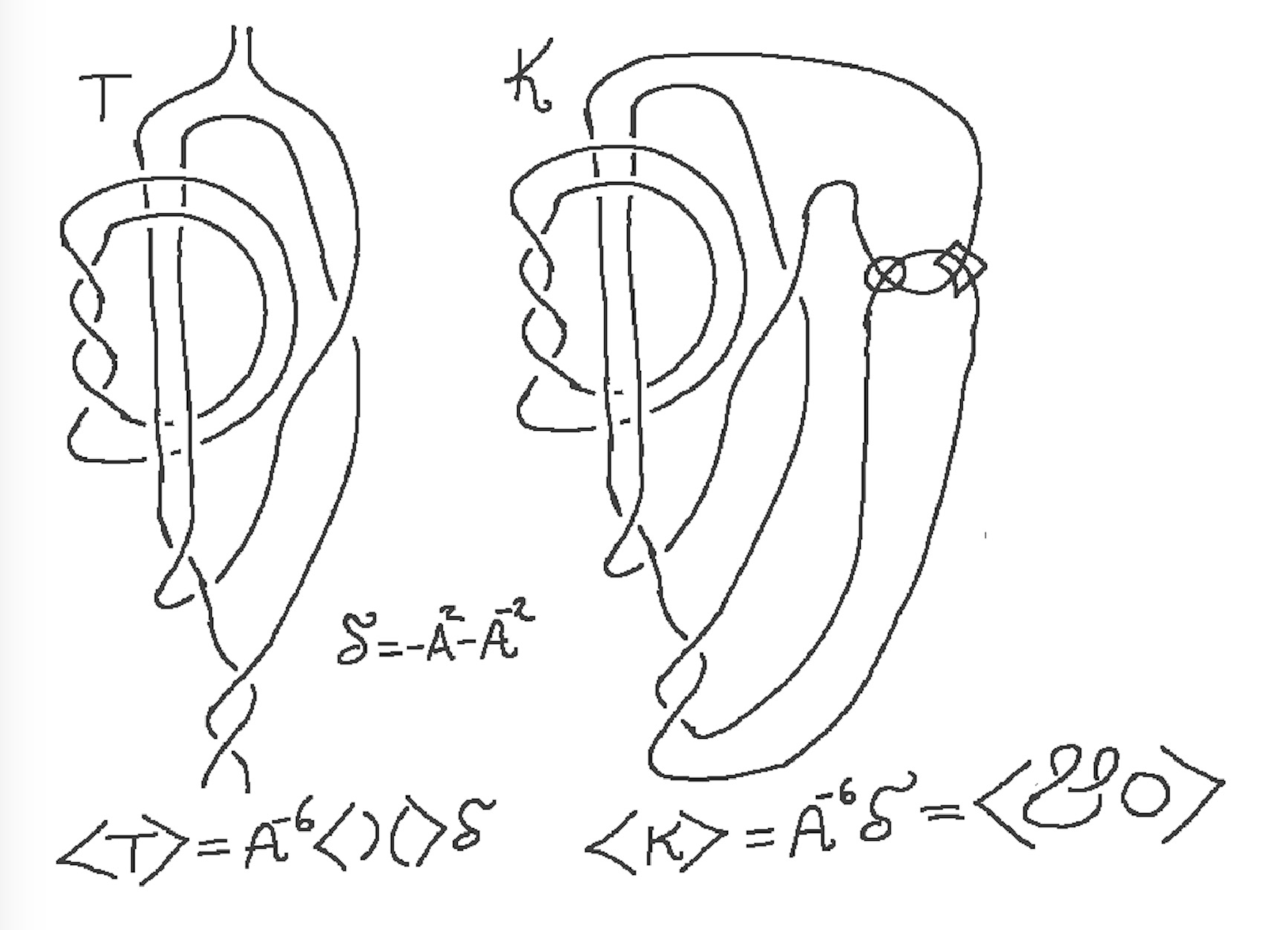}
     \end{tabular}
     \caption{\bf Example of Non-Trivial Double Virtual Knot whose Virtuality is Invisible to Generalized Bracket}
     \label{EFF28}
\end{center}
\end{figure}

In Figure~\ref{EFF27} we illustrate a multiple virtual link $K$ that is composed as a closure of a classical tangle
$T$ with a ``virtual clasp" as illustrated. The direct closure of the virtual clasp consists in two circles that meet in two virtual 
crossings of distinct type. Call these two circles $\Lambda$ We have already seen in Figure~\ref{EFF12} that the graphical evaluation of  $\Lambda$ is 
$2 \delta^2 - \delta$ where $\delta = -A^2 -A^{-2}.$ In the Figure~\ref{EFF27} we see that the bracket evaluation of $K$ is given by the formula
$\langle K \rangle = f \langle \Horiz  \rangle  + g \langle \Virt \rangle $ where $f$ and $g$ are some coefficients dependent on the structure of $T.$
Then we conclude that $$\langle K \rangle = f \langle \Lambda \rangle + g \delta ,$$ from which we conclude that if the tangle $T$ has $f(T) = 0,$ then 
the generalized multibracket polynomial will be unable to distinguish $K$ from the unknot or unlink. In Figure~\ref{EFF28} we give exactly such an example.
The tangle $T$ in that figure is derived in \cite{DKT} and is based on the paper \cite{EKT}.
The link $K$ in Figure~\ref{EFF28} is indetectable to the generalized bracket. That it is a non-trivial multi-virtual link  will be shown in a subsequent paper.\\

\begin{figure}
     \begin{center}
     \begin{tabular}{c}
     \includegraphics[width=12cm]{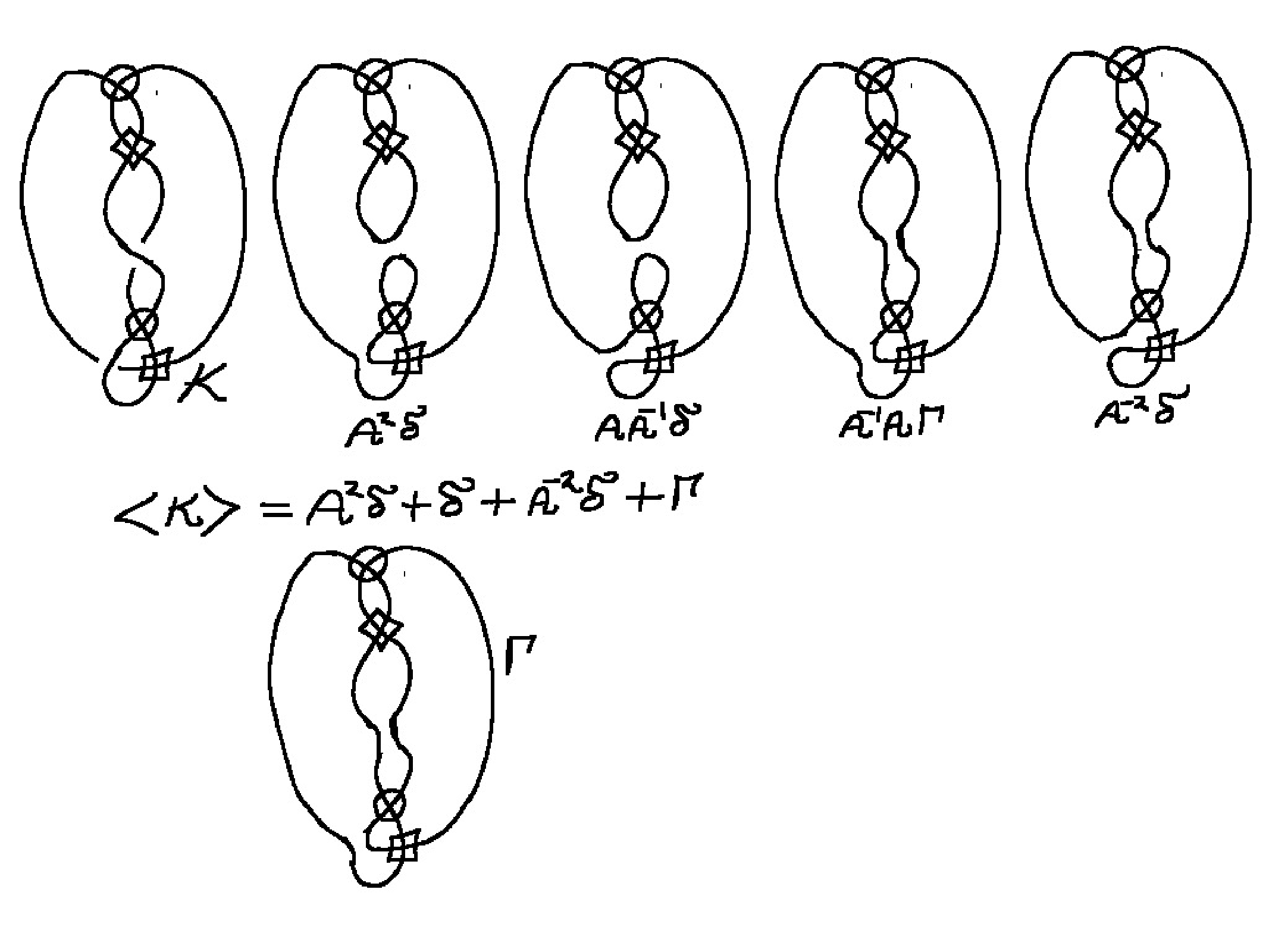}
     \end{tabular}
     \caption{\bf Generalized bracket of a multiple virtual slice knot K.}
     \label{EFF30}
\end{center}
\end{figure}

The next example is shown in Figure~\ref{EFF30}. We call this knot $K$ and find that $\langle K \rangle = A^2 \delta + \delta + A^{-2} \delta + \Gamma$ where $\Gamma$ denotes the 
detour equivalence class of the diagram labeled $\Gamma$ in the figure. It is easy to see the the chromatic evaluation of this diagram does not distinguish $\Gamma$ from two disjoint loops.
Thus we must use a different technique to show that $\Gamma$ is non-trivial and so prove that $K$ is a non-trivial knot. In the next section, we introduce quandle and biquandle invariants of 
multi-virtuals, and we shall apply these invariants to the state $\Gamma.$ We will also point out that the knot $K$ in this figure is an example of a slice knot in a generalized context of 
multi-virtual knot concordance.\\

\section{Quandles and Biquandles for Multiple Virtual Knots and Links}

In this section we introduce generalized quandles and biquandles for multiple virtual knots and links. \\

Recall the definition of a {\it quandle}. A quandle is a set $Q$ together with binary operations $\star, \sharp$ satisfying the following axioms:

\begin{enumerate}
\item $a \star a = a$ and $ a \sharp a = a$ for all $a \in Q.$
\item $ (a \star b) \sharp b = a$ and $ (a \sharp b) \star b = a$ for $a, b \in Q.$
\item $(a \star b) \star c = (a \star c) \star (b \star c)$ and 
$(a \sharp b) \sharp c = (a \sharp c) \sharp (b \sharp c)$ for $a, b, c \in Q.$
\end{enumerate}

\begin{figure}
     \begin{center}
     \begin{tabular}{c}
     \includegraphics[width=10cm]{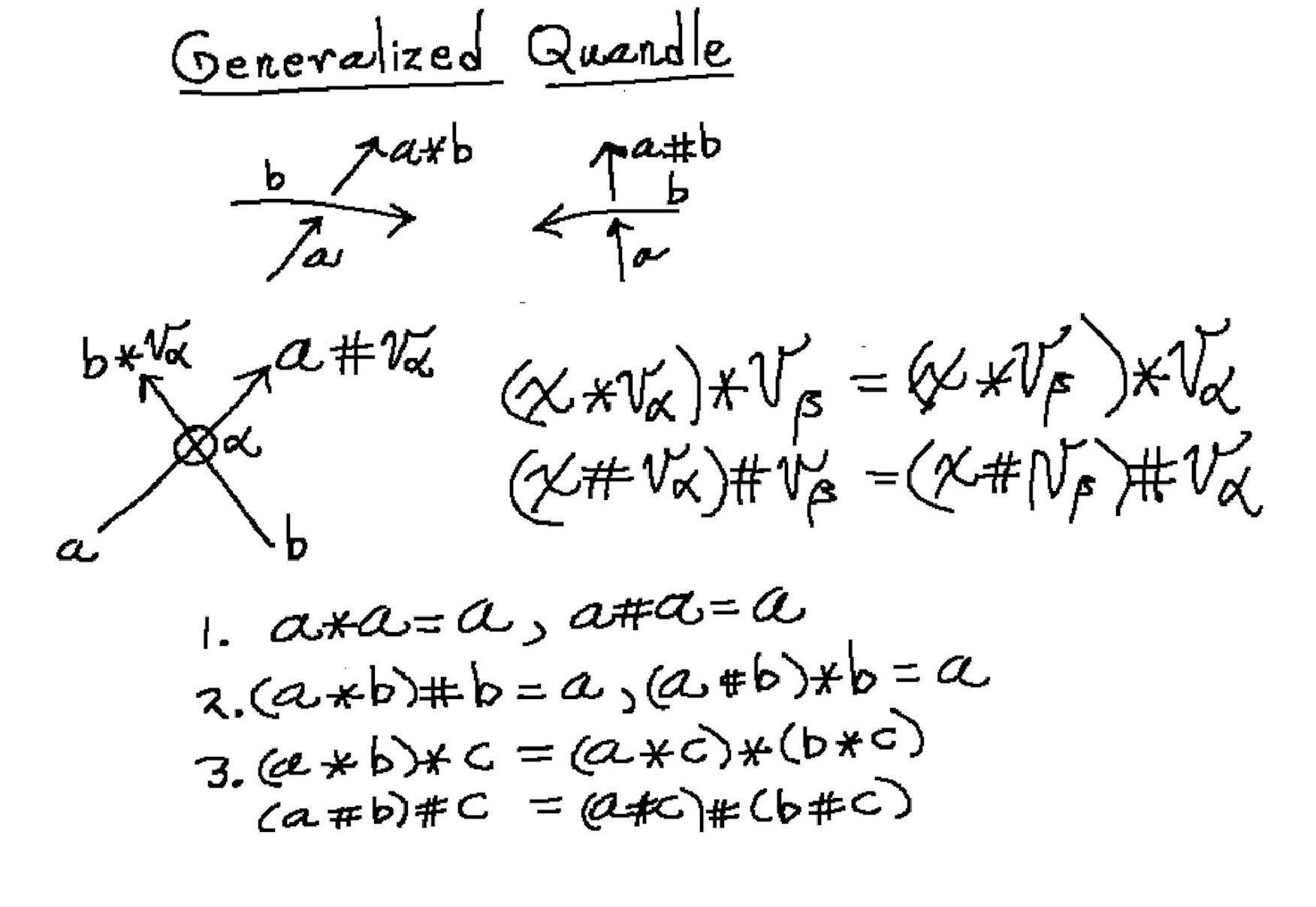}
     \end{tabular}
     \caption{\bf Generalized Quandle}
     \label{EFF22}
\end{center}
\end{figure}

These axioms correspond directly to the three Reidemeister moves when we have an oriented link diagram whose arcs are labeled with elements of a quandle, and so that 
quandle operations $\star$ and $\sharp$ correspond to positive and negative crossing types as shown in Figure~\ref{EFF22}.\\

\noindent {\bf Remark}. Sometimes it is convenient to denote quandle operations by exponent notation. Thus we can write
$$a \star b = a^{b}$$ and $$a \sharp b = a^{\bar{b}}.$$ It is understood that the ``bar" in $\bar{b}$ is strictly notational and does not designate a separate operation on $b.$ This notation is used in Figure~\ref{EFF31} and in the next section on biquandles.\\

\noindent {\bf Definition.} Consider a multiple virtual knot theory with virtual crossings corresponding to indices $\alpha$ in an index set $\cal {A}.$ To each index $\alpha$ associate a free quandle variable
$v_{\alpha}.$ Now associate quandle relations at each classical crossing of the form $c = a \star b$ for positive crossings and $c = a \sharp b$ for negative crossings, as shown
in Figure~\ref{EFF22}. To each virtual crossing associate relations of the form $c = a \star v_{\alpha}$ and $d = b \sharp v_{\alpha}$ as illustrated in Figure~\ref{EFF22}.
Another way to say this is that we can regard $v_{\alpha}$ as a quandle automorphism. Note that $(a\star b) \star v_{\alpha} = (a \star v_{\alpha} )\star (b \star v_{\alpha})$ so that operating on elements of the quandle by  $v_{\alpha}$  preserves the quandle structure.\\

Given a multiple virtual knot or link diagram $K$, let there be one symbol for each arc of the diagram $K$ where an {\it arc} of the diagram extends
from an undercrossing to an undercrossing when meeting no virtual crossings, or from an undercrossing to a virtual crossing, or from a virtual crossing to a virtual crossing.
Thus an arc of the diagram has a constant symbol across any overcrossing. See Figure~\ref{EFF22} for examples of arc labelings. 
Let $Q(K)$ denote the quandle generated by the symbols associated with its arcs and the free variables associated with its virtual crossings, with the relations described above for 
each of its classical and virtual crossings. We call $Q(K)$ the {\it Multi-Virtual Quandle} associated to the link $K.$\\

\noindent {\bf Proposition.} The Multi-virtual Quandle $Q(K)$ associated with a multi-virtual diagram $K$ is invariant under the classical Reidemeister moves and all detour moves.
Hence $Q(K)$ is an invariant of multi-virtual links.\\

\begin{figure}
     \begin{center}
     \begin{tabular}{c}
     \includegraphics[width=10cm]{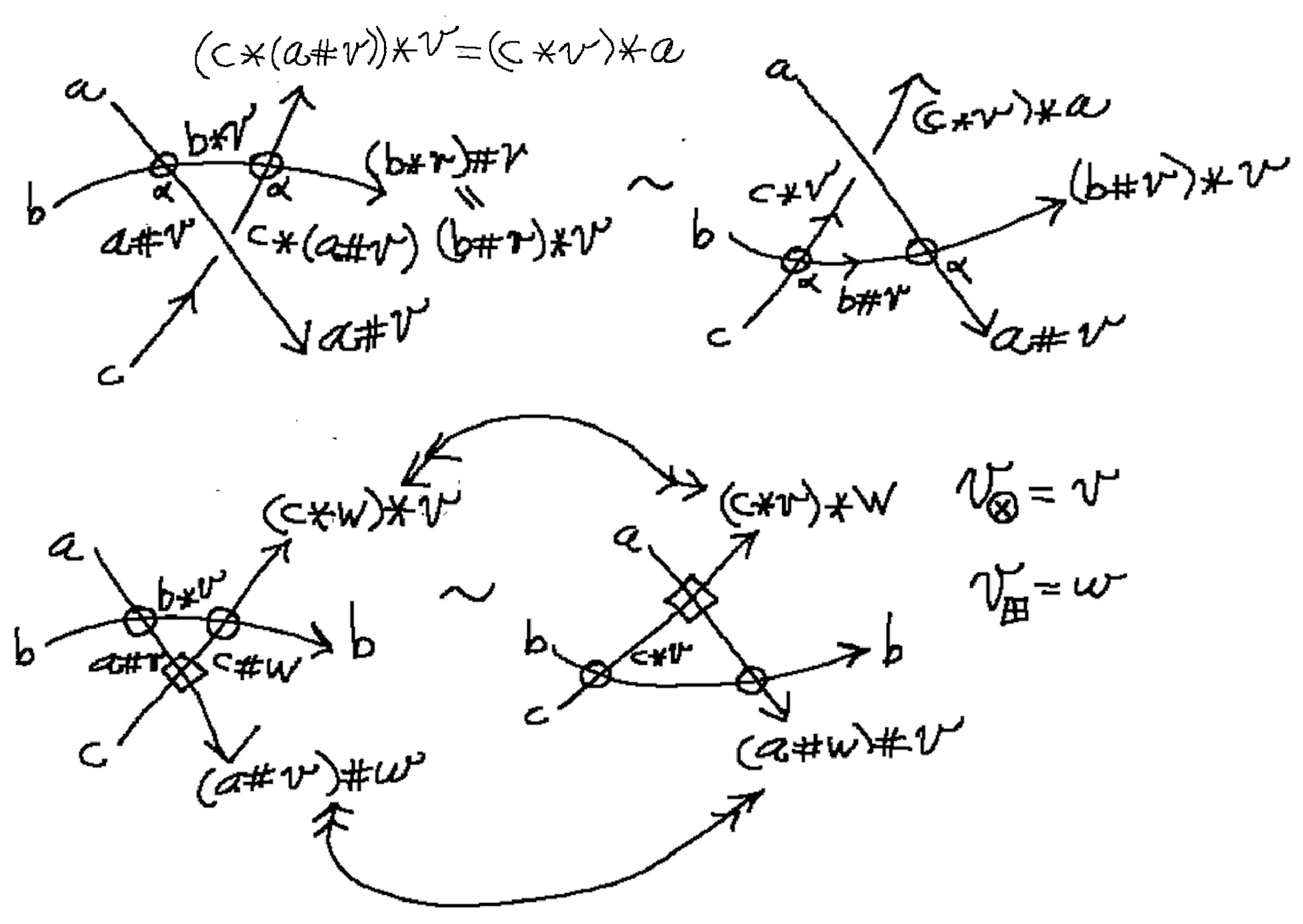}
     \end{tabular}
     \caption{\bf Generalized Quandle Invariance I}
     \label{EFF23}
\end{center}
\end{figure}

\begin{figure}
     \begin{center}
     \begin{tabular}{c}
     \includegraphics[width=10cm]{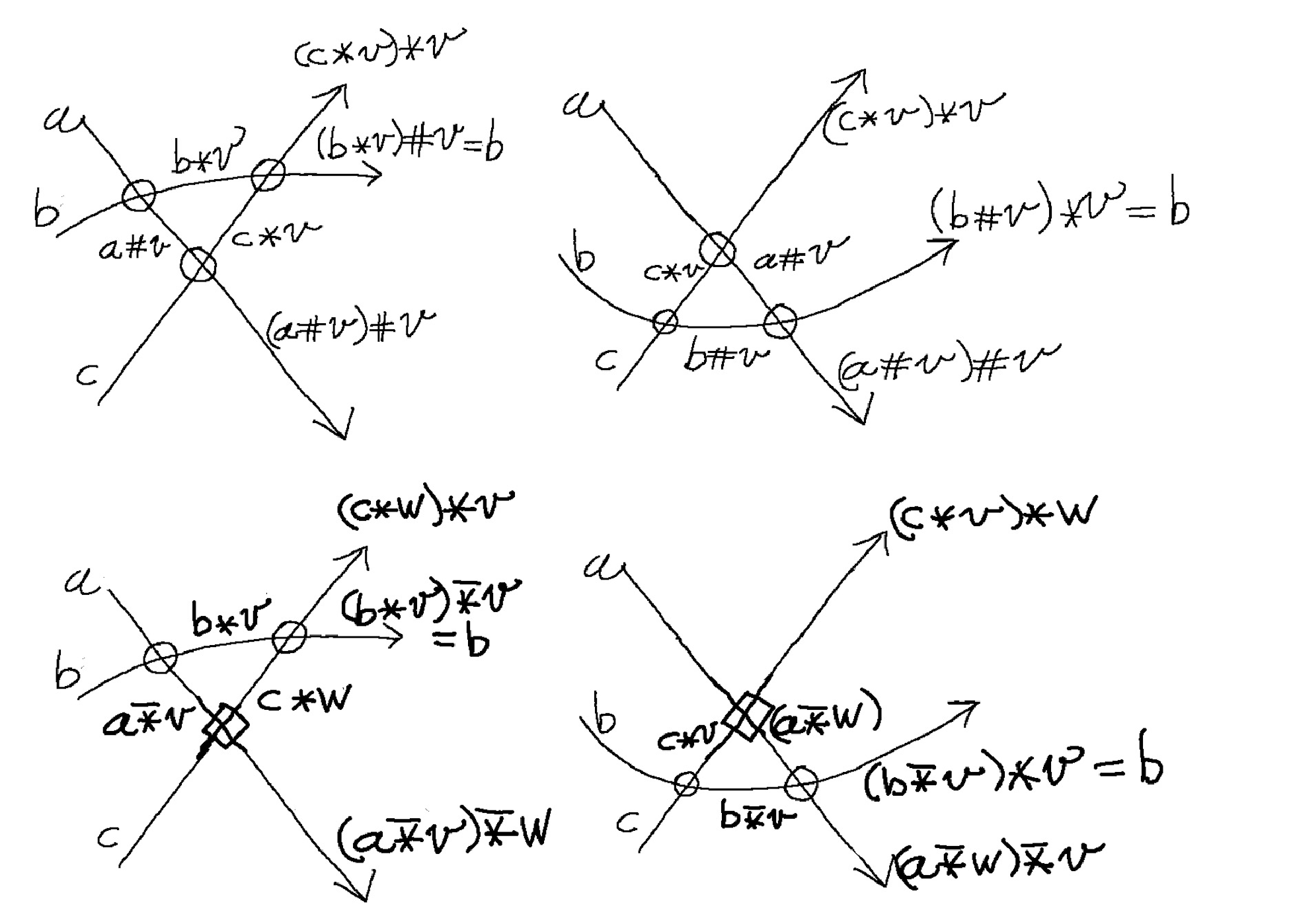}
     \end{tabular}
     \caption{\bf Generalized Quandle Invariance II}
     \label{EFF23.1}
\end{center}
\end{figure}

\begin{figure}
     \begin{center}
     \begin{tabular}{c}
     \includegraphics[width=10cm]{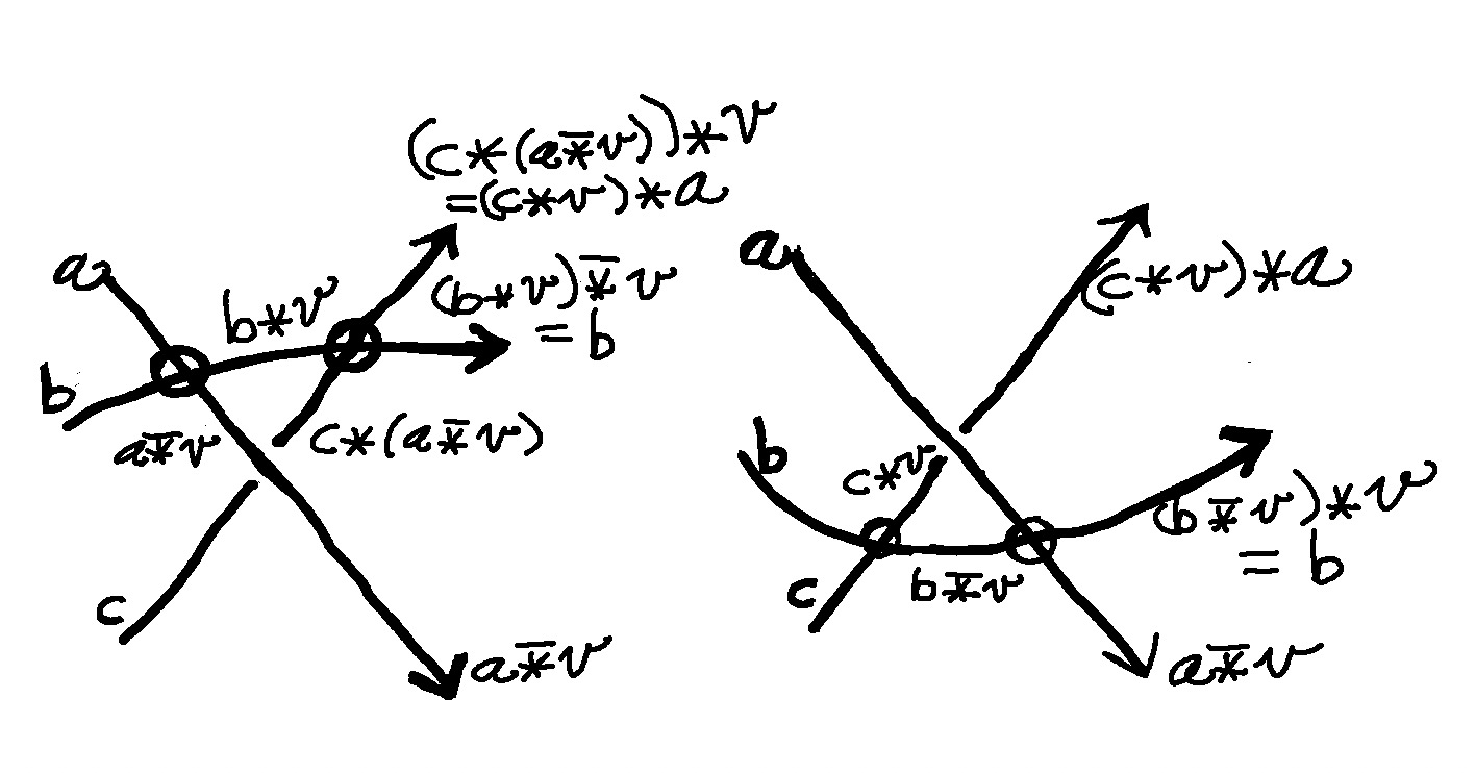}
     \end{tabular}
     \caption{\bf Generalized Quandle Invariance III}
     \label{EFF23.11}
\end{center}
\end{figure}

\noindent {\bf Proof.} The invariance under the classical Reidemeister moves is guaranteed by the three basic quandle axioms.  
The invariance under detours is verified in Figure~\ref{EFF23}, Figure~\ref{EFF23.1} and Figure~\ref{EFF23.11}. Variations on the illustrated cases are left to the reader. This completes the proof of the proposition.$\hfill\Box$ \\

\begin{figure}
     \begin{center}
     \begin{tabular}{c}
     \includegraphics[width=10cm]{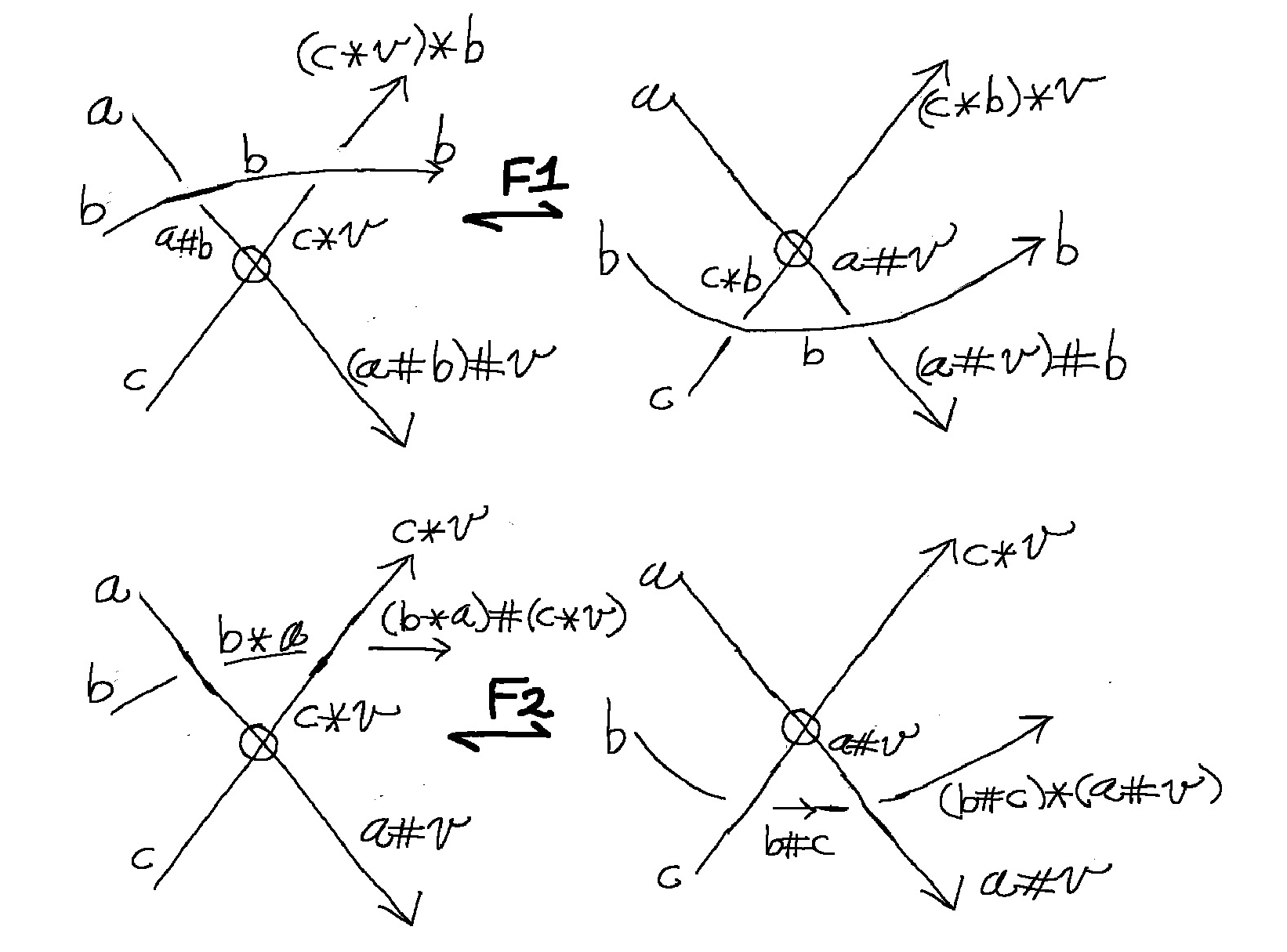}
     \end{tabular}
     \caption{\bf Forbidden Moves in the Quandle }
     \label{EFF23.2}
\end{center}
\end{figure}

\begin{figure}
     \begin{center}
     \begin{tabular}{c}
     \includegraphics[width=8cm]{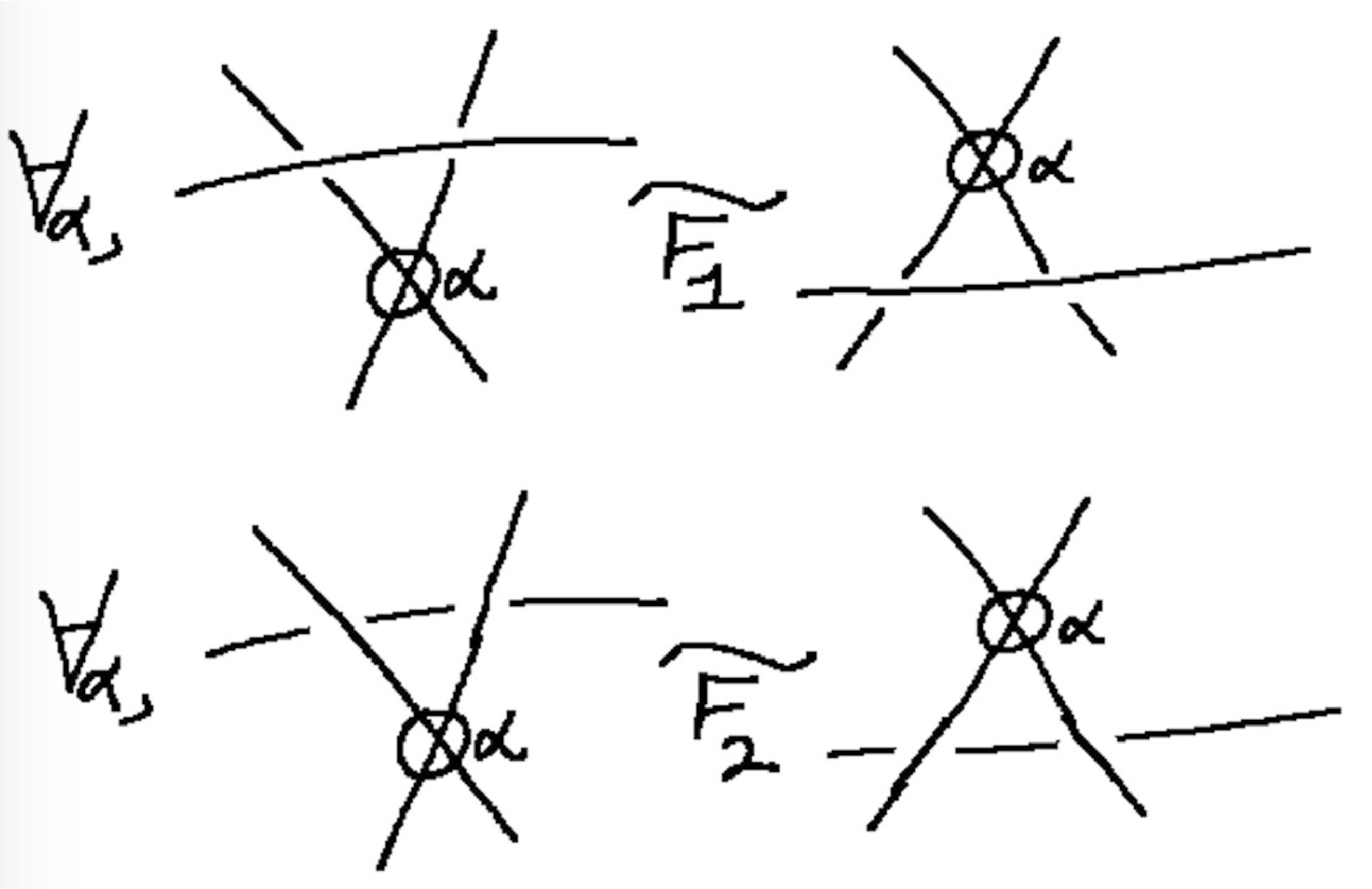}
     \end{tabular}
     \caption{\bf Generalized Forbidden Moves}
     \label{EFF21}
\end{center}
\end{figure}

\noindent{\bf Remark.} Note that one may wish to examine how these generalized quandles behave under the forbidden moves  of Figure~\ref{Figure 3} and Figure~\ref{EFF21}. Note that in Figure~\ref{EFF21} we illustrate the multiplicity of forbidden moves that arise in the presence of many virtual crossing types. When we say that moves are ``forbidden" in a given theory, this does not mean that an invariant for that theory might not
have the same values or be invariant under these moves. The maxim that corresponds to axioms for moves is that ``what is not allowed is forbidden". We emphasize certain forbidden moves to remind ourselves and the reader that these moves are not consequences of the axioms and so invariants are not required to respect these moves. In the case of the forbidden moves under consideration, the bracket polynomials we have already discussed show that these moves are not consequences of the given axiomatic moves. It is often the case that quandles we define are also not invariant under these moves. In detail, see Figure~\ref{EFF23.2} where we see that for a forbidden move of type $1$ to occur at the quandle level, the virtual automorphism representative $v$ would have to commute with an arbitrary element of the quandle. Since we only require the automorphisms for different virtual crossings to commute with one another, this extra relation will generically not occur.
Similarly, the requirement for the forbidden move of type $2$ is also an extra algebraic relation.\\

\begin{figure}
     \begin{center}
     \begin{tabular}{c}
     \includegraphics[width=8cm]{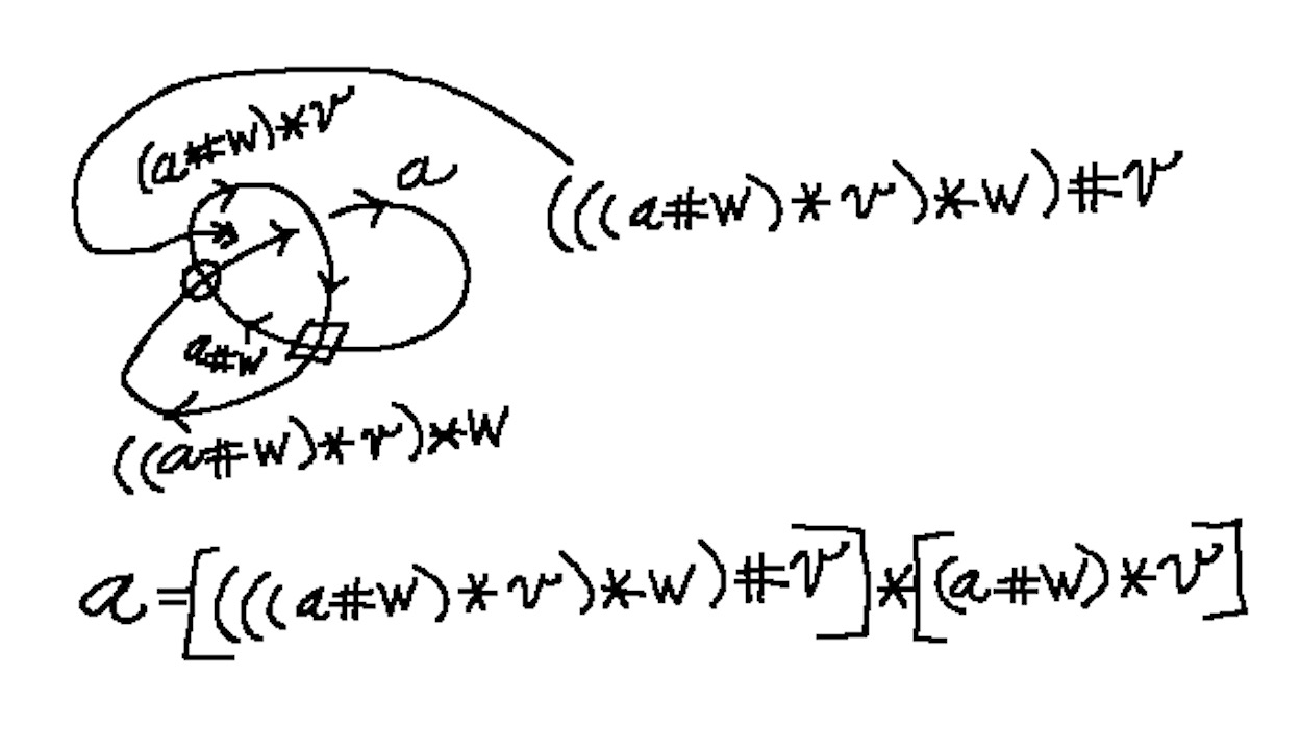}
     \end{tabular}
     \caption{\bf Quandle Coding of Double Virtual Trefoil}
     \label{EFF24}
\end{center}
\end{figure}

In Figure~\ref{EFF24} we illustrate the quandle coding of the double virtual trefoil $K$, a multi-virtual knot with one classical crossing, one round virtual crossing and one box virtual crossing. With the labelling as shown in the figure, and virtual free generators $v$ and $w$ we have that $Q(K)$ has generators $a,v,w$ and one relation
$a = [(((a\sharp w)\star v)\star w)\sharp v]\star[(a \sharp w) \star v].$ We wish to show that this quandle is non-trivial.\\

\begin{figure}
     \begin{center}
     \begin{tabular}{c}
     \includegraphics[width=12cm]{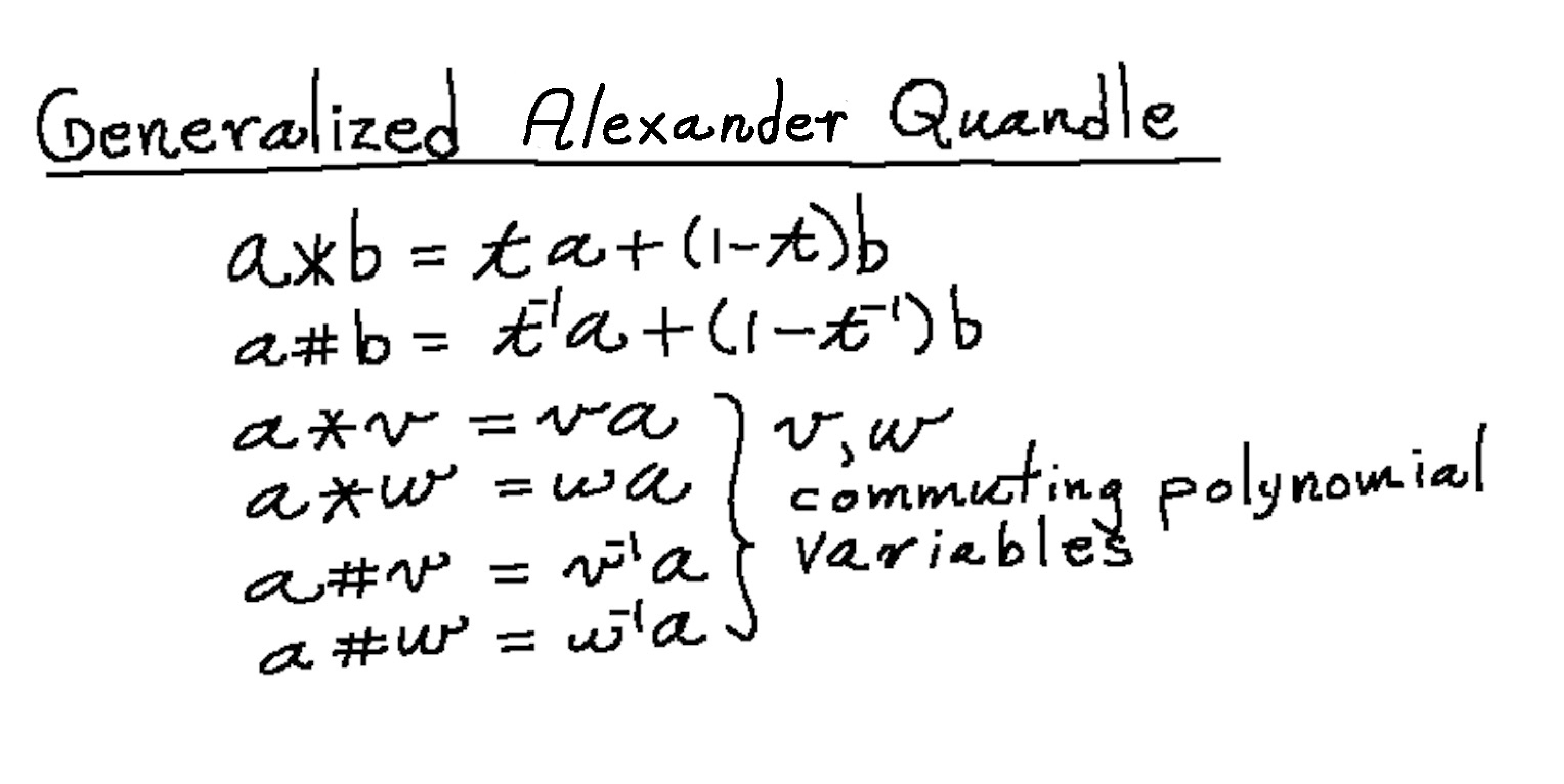}
     \end{tabular}
     \caption{\bf Generalized Alexander Quandle}
     \label{EFF25}
\end{center}
\end{figure}

To facilitate calculation Figure~\ref{EFF25} shows the defining relations for the {\it Generalized Alexander Quandle} where the quandle is a module over the 
ring $Z[t,t^{-1}, \{ v_{\alpha}, v_{\alpha}^{-1} \}]$ of integer-coefficient Laurent polynomials in the variable $t$ and in a set of commuting Laurent polynomial variables corresponding to the 
multiple virtual crossing types. In the figure we illustrate this quandle for just the two virtual types (round and box) with variables $v$ and $w$ respectively. Thus the ring is taken to be  a ring $Z[t,t^{-1}, v, v^{-1}, w, w^{-1}]$ of Laurent polynomials in three variables. The variables $t, v, w$ all commute with one another. Since elements of the quandle are modules over this ring, we can define quandle operations through linear multiplication by the elements of the ring. Then we have 

\begin{enumerate}
\item $a \star b = ta + (1-t)b$ and $a \sharp b = t^{-1} a + (1-t^{-1})b.$
\item $a\star v = va, a\star w = wa, a \sharp v = v^{-1} a, a \sharp w = w^{-1}a$
\item Here $v$ and $w$ are commuting polynomial variables.
\end{enumerate}

\begin{figure}
     \begin{center}
     \begin{tabular}{c}
     \includegraphics[width=10cm]{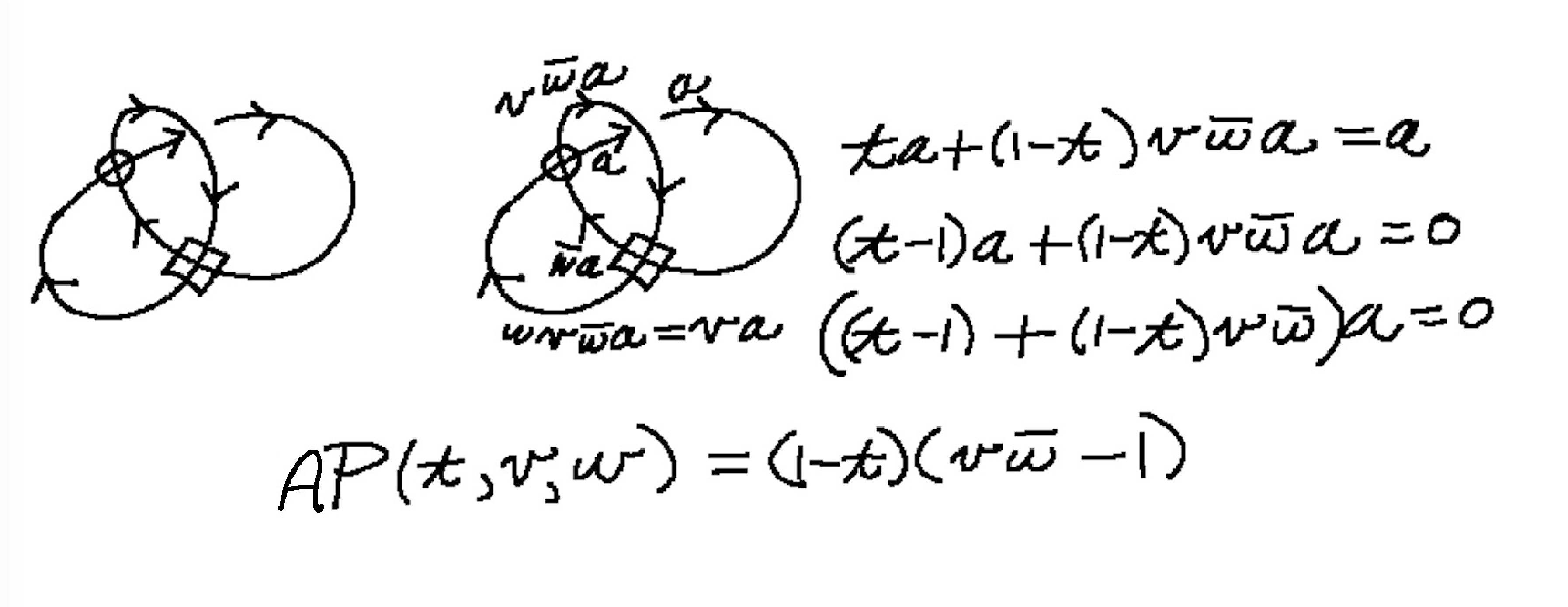}
     \end{tabular}
     \caption{\bf Calculation of Generalized Alexander Polynomial}
     \label{EFF26}
\end{center}
\end{figure}

\begin{figure}
     \begin{center}
     \begin{tabular}{c}
     \includegraphics[width=8cm]{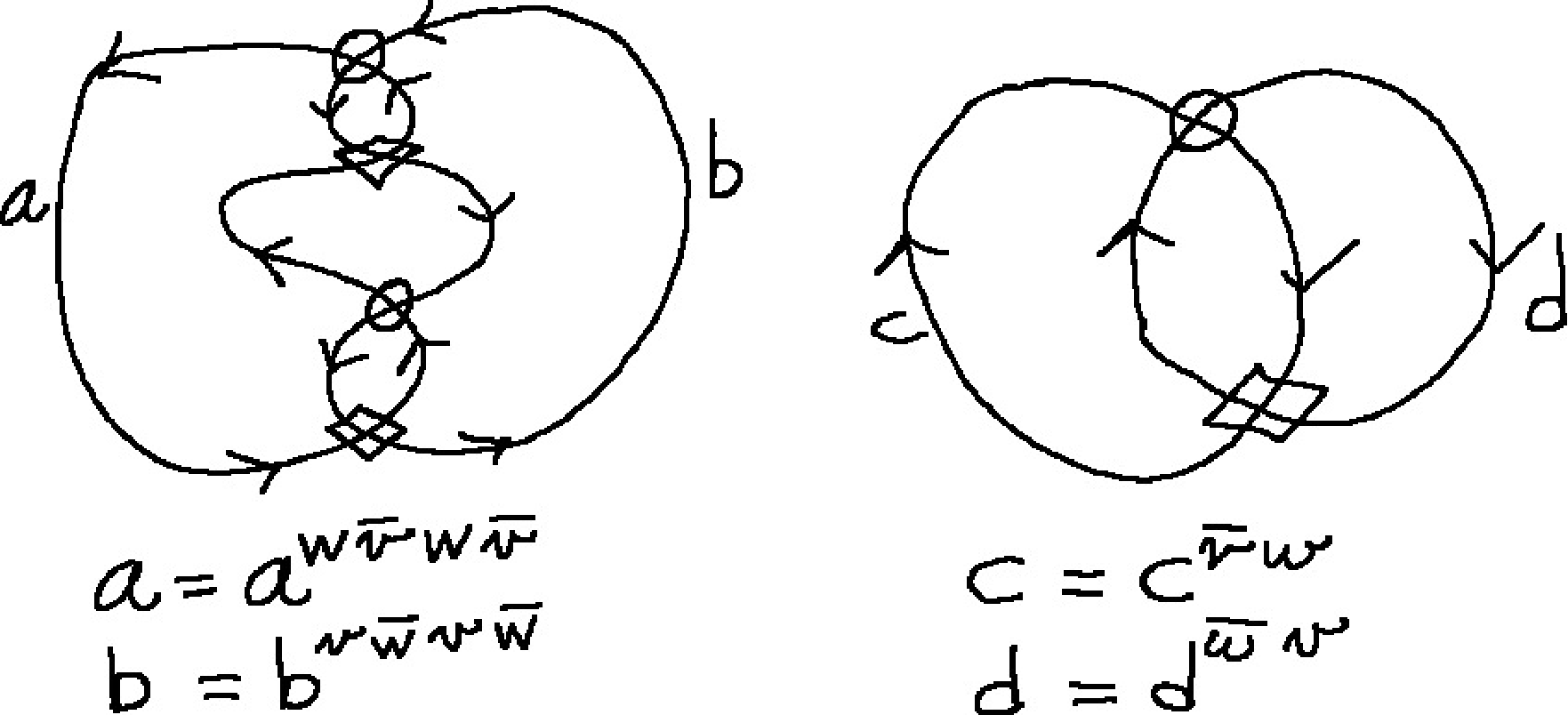}
     \end{tabular}
     \caption{\bf Using multi-quandle to detect non-triviality of virtual states.}
     \label{EFF31}
\end{center}
\end{figure}

It is not hard to verify that the generalized Alexander quandle satisfies all the quandle properties. 
In Figure~\ref{EFF26} we illustrate the labeling and production of the basic relation for the generalized Alexander quandle of the double trefoil.
From that figure the relation is $$(ta + (1-t)v w^{-1}) a = a.$$
The reader should note that this relation also follows algebraically from the general relation we showed for Figure~\ref{EFF24}. That is, we have
$$[(((a\sharp w)\star v)\star w)\sharp v]\star[(a \sharp w) \star v] = $$
$$[v^{-1}wvw^{-1}a]\star[vw^{-1}a] =$$
$$[a]\star[vw^{-1}a] =$$
$$ta + (1-t) v w^{-1}a.$$
Thus the relation is
$$ta + (1-t) v w^{-1}a =  a,$$ as we show directly in Figure~\ref{EFF26}.
Now simplifying this relation we see that it becomes
$$(t-1)a + (1-t) v w^{-1}a = 0$$ 
whence
$$[(t-1) + (1-t) v w^{-1}]a = 0.$$ 
or
$$[(t-1)(1-v w^{-1})]a = 0.$$ 
It is not hard to show that
$$AP(K) = (t-1)(1-v w^{-1})$$
is a well defined invariant of $K$ with value in  the ring $$Z[t,t^{-1}, v, v^{-1}, w, w^{-1}]$$
taken up to multiples of monomials in $t,v,w$ of the form $\pm t^a v^b w^c.$ This is the {Generalized Alexander Polynomial} of the 
multiple virtual knot $K.$ Thus the knot $K$ is non-trivial since we have shown that its polynomial $AP(K)$ is non-trivial. The non-triviality of the polynomial implies that the quandle from which it was computed is non-trivial and that the module that represents this quandle is non-trivial.\\

\noindent {\bf Remark.} The generalized Alexander polynomial that we have here indicated by example is a generalization of the approach to such generalizations found in 
\cite{BIQ}. When we discuss biquandles (below) we will indicate further generalizations along similar lines.\\

In Figure~\ref{EFF31} we illustrate the quandle relations for two virtual states that can appear in the expansion of knots and links by the generalized bracket. The first example is the state
$\Gamma$ that first appeared in Figure~\ref{EFF30}. We pointed out that $\Gamma$ was not visible to the chromatic calculation. Here we see that  the quandle for $\Gamma$ has generators
$a,b$ with relations $a = a^{w \bar{v} w \bar{v}}, b = b^{v \bar{w} v \bar{w}},$ from which it follows that $\Gamma$ is not detour equivalent to disjoint loops. This in turn shows that the knot $K$ in Figure~\ref{EFF30} is non-trivial and non-classical, requiring two virtual crossing types. Also in Figure~\ref{EFF31} we show the quandle for the state $\Lambda$ consisting in two loops intersecting in one round virtual and one box virtual. This quandle is also non-trivial and we have previously detected $\Lambda$ via chromatic evaluation.\\

\clearpage

\subsection{Generalized Biquandles}
Biquandles are a generalization of quandles where two operations occur at each crossing. We will use exponential notation to describe these operations. We will write
$$a \star b = a^{b}$$ and $$a \sharp b = a^{\bar{b}}.$$ It is understood that the ``bar" in $\bar{b}$ is strictly notational and does not designate a separate operation on $b.$ \\

\noindent {\bf Definition.}  A {\em biquandle} $B$ is a set with four binary operations indicated by the conventions we have explained
above:  $a^{b} \,\mbox{,} \, a^{\o{b}} \, \mbox{,} \,  a_{b} \,\mbox{,} \, a_{\o{b}}.$ We shall refer to the operations with barred 
variables as the {\em left} operations and the operations without barred variables as the {\em right} operations. The biquandle is 
closed under these operations and the following axioms are satisfied:

\begin{enumerate}
\item   For any elements $a$ and $b$ in $B$ we have 

$$a = a^{b \o{b_{a}}}  \quad \mbox{and} \quad  b= b_{a \o{a^{b}}} \quad \mbox{and}$$

$$a = a^{\o{b}b_{\o{a}}}  \quad \mbox{and} \quad  b= b_{\o{a} a^{\o{b}}}.$$

\item  Given elements $a$ and $b$
in $B$, then there exists an element $x$ such that 

$$x = a^{b_{\o{x}}} \mbox{,} \quad a = x^{\o{b}} \quad \mbox{and} \quad b= b_{\o{x}a}.$$

\noindent Given elements $a$ and $b$
in $B$, then there exists an element $x$ such that 

$$x = a^{\o{b_{x}}} \mbox{,} \quad a = x^{b} \quad \mbox{and} \quad b= b_{x\o{a}}.$$

\item For any $a$ , $b$ , $c$ in $B$ the following equations hold and the same equations hold when all right operations are 
replaced in these equations by left operations.

$$a^{b c} = a^{c_{b} b^{c}} \mbox{,} \quad c_{b a} = c_{a^{b} b_{a}} \mbox{,} \quad (b_{a})^{c_{a^{b}}} = (b^{c})_{a^{c_{b}}}.$$

\item Given an element $a$ in $B$, then there exists an $x$ in the biquandle such that $x=a_{x}$ and  
$a = x^{a}.$ Given an element $a$ in $B$, then there exists an $x$ in the biquandle such that $x=a^{\o{x}}$ and  
$a = x_{\o{a}}.$
\end{enumerate}

\begin{figure}
     \begin{center}
     \begin{tabular}{c}
     \includegraphics[width=10cm]{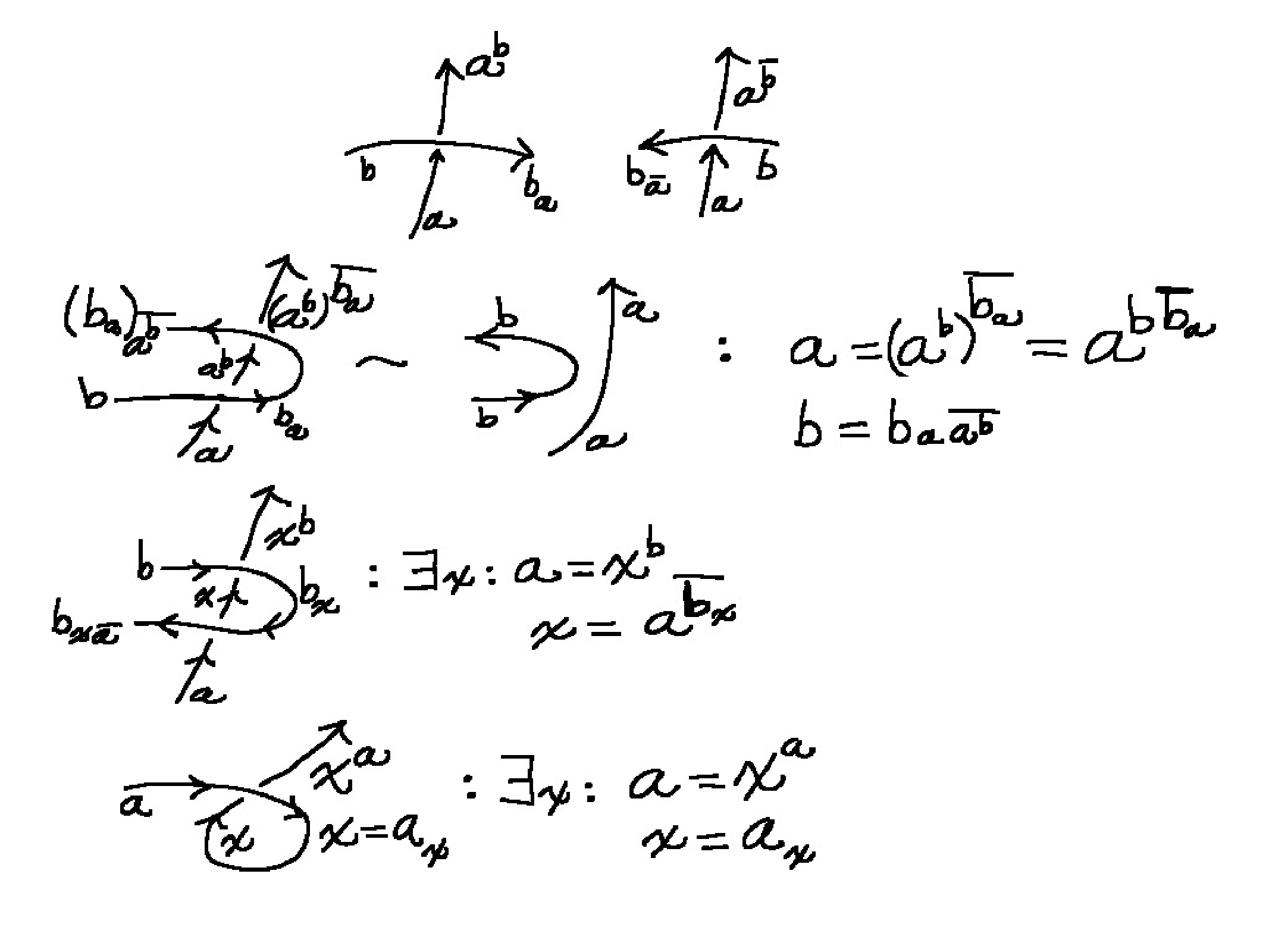}
     \end{tabular}
     \caption{\bf Biquandle Formalism and Reidemeister Moves Two and One.}
     \label{BQ1}
\end{center}
\end{figure}

\begin{figure}
     \begin{center}
     \begin{tabular}{c}
     \includegraphics[width=10cm]{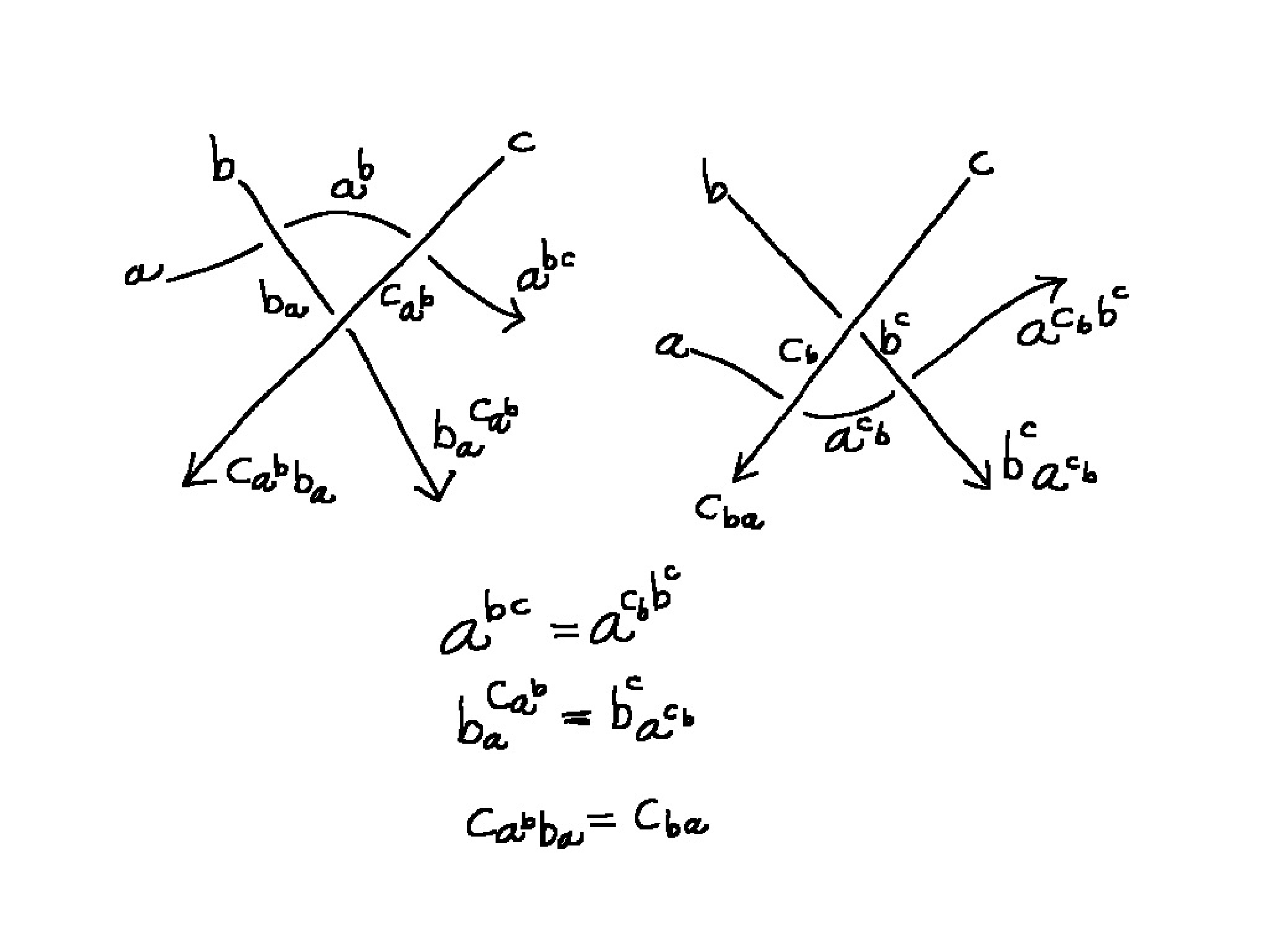}
     \end{tabular}
     \caption{\bf Biquandle Formalism and Reidemeister Move Three.}
     \label{BQ2}
\end{center}
\end{figure}

\begin{figure}
     \begin{center}
     \begin{tabular}{c}
     \includegraphics[width=10cm]{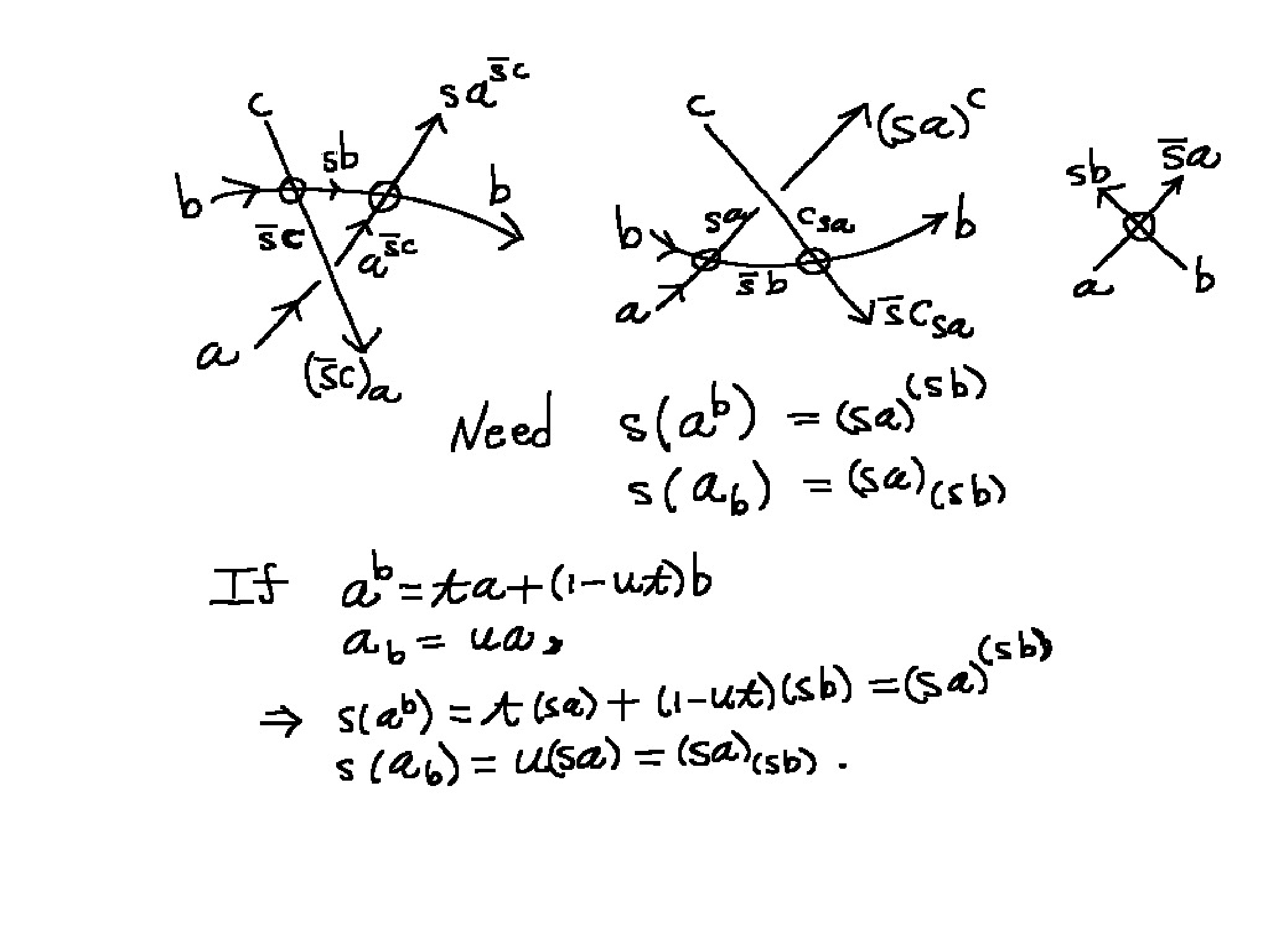}
     \end{tabular}
     \caption{\bf Generalized Biquandle Formalism with respect to Automorphisms at Virtual Crossings.}
     \label{BQ3}
\end{center}
\end{figure}

These axioms are transcriptions of the Reidemeister moves.  
In Figure~\ref{BQ1} we illustrate the link diagrammatic formalism for the biquandle. The two basic operations $a^{b}$ and $b_{a}$ correspond to arrow outputs at a crossing where the 
inputs $a$ and $b$ are at the base of the arrow. This figure illustrates the identities for Reidemeister Two and Reidemeister One moves. In Figure~\ref{BQ2} the labeling leading to the idenities for the Reidemeister Three move are shown. Note the use of the unary operation $\bar{s}$ to denote inverse operations and the use of the formalism $(a^{b})^c = a^{bc}.$ These aspects are used in quandle formalism as well. In Figure~\ref{BQ3} we show the formalism for a biquandle automorphism at a virtual crossing. If $B$ is the biquandle, then we are given
$s:B \longrightarrow B$ such that $s(a^{b}) = (sa)^{sb}$, $s(a_{b}) = (sa)_{sb}$, $s(\bar{a}) = \bar{(sa)}.$ Such an automorphism applied at a virtual crossing, as we have discussed for quandles,
can be regarded as extending the biquandle structure to virtual diagrams. This figure illustrates that when the $s$ is a quandle automorphism, then the virtual crossing will have a detour move with respect to classical crossings. We obtain generalized biquandles by choosing automorphisms $s_{\alpha}$ for each virtual crossing type $\alpha.$ As in the case of quandles, it is then sufficient that $s_{\alpha} s_{\beta} = s_{\beta} s_{\alpha} $ for all $\alpha, \beta$ to have detour moves satisfied among the different virtual crossings. This suffices to define the generalized
biquandle for multiple virtual knots and links.\\

In Figure~\ref{BQ3} we also indicate a specific {\it generalized Alexander biquandle}.  Let $B$ denote this biquandle. As a set, $B$ is a module over the commutative ring 
$$R = Z[t, t^{-1}, u, u^{-1}, s_{\alpha}, s_{\alpha}^{-1}]$$ where $\alpha$ runs over the labels for the virtual crossings intended for use. The operations in $B$ are defined as follows:
$$a^{b} = t a + (1-ut) b,$$
$$a^{\bar{b}} = t^{-1} a + (1 - u^{-1} t^{-1})b,$$
$$a_{b} = ua,$$
$$a_{\bar{b}} = u^{-1} b.$$
The automorphisms are multiplication by $s_{\alpha}$ or $s_{\alpha}^{-1}$ accordingly.\\

Taking the this biquandle for a multi-vitual link $K$ one has a generalized Alexander polynomial $GA_{K} = GA_{K}[t, t^{-1}, u, u^{-1}, s_{\alpha}, s_{\alpha}^{-1}]$ that is the generator of the ideal of polynomials that annihilate all elements of the module $B.$ This is a zero-th order Alexander polynomial. The structure of the biquandle module for $K$ can yield other topological information and higher order polynomials. See \cite{FennBarth,FennGen,KFJ,vkt,Hren} for more information about classical biquandles and other directions for generalization.\\

\section{\bf Virtual Link Cobordism.}  
\begin{figure}
     \begin{center}
     \begin{tabular}{c}
     \includegraphics[width=8cm]{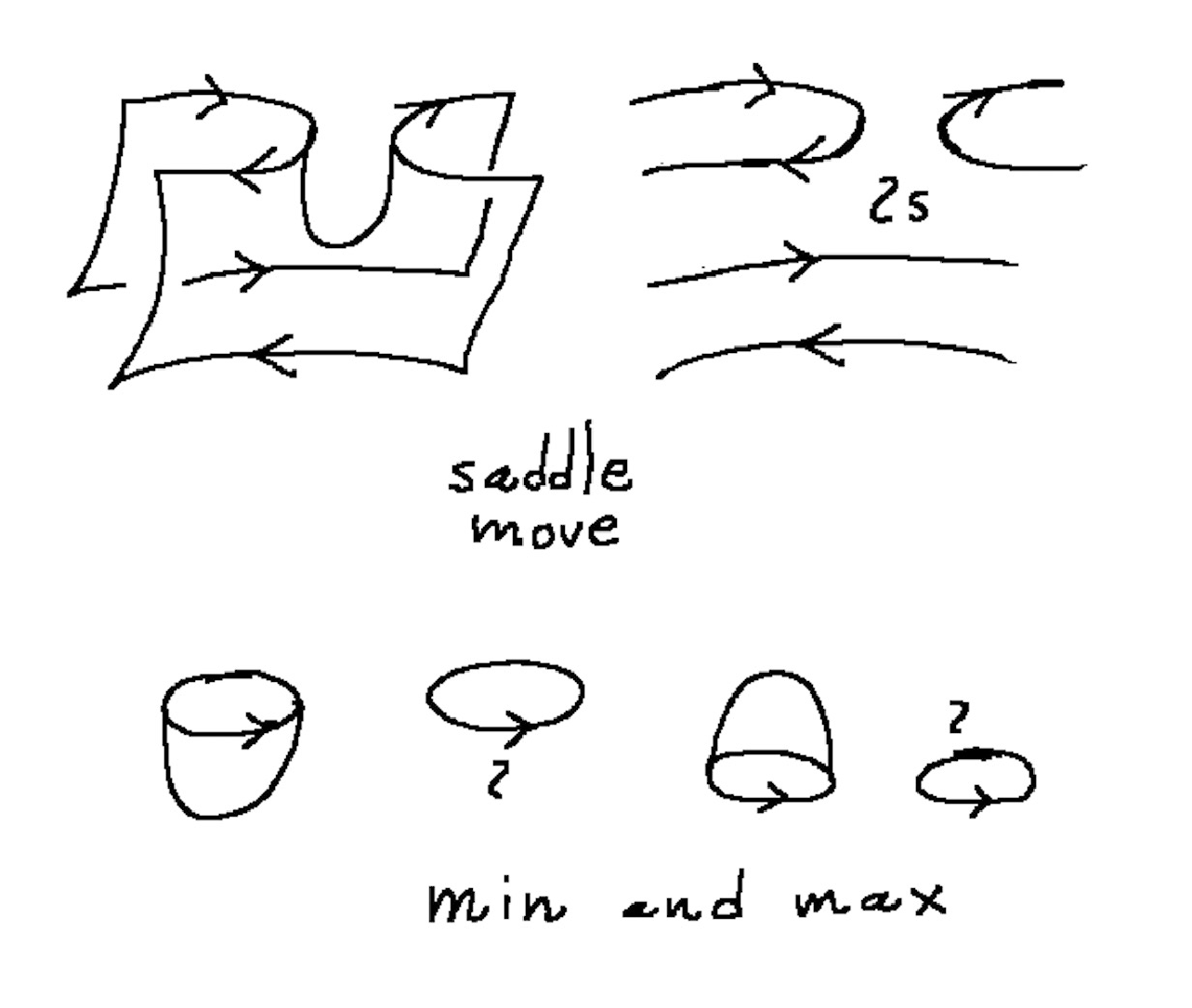}
     \end{tabular}
     \caption{\bf Moves for Cobordisms.}
     \label{Saddle}
\end{center}
\end{figure}

In Figure~\ref{Saddle} we indicate oriented saddle moves and birth and death moves on link diagrams.
We call two (multiple) virtual link diagrams {\it cobordant} if one can be obtained from the other by a combination of Reidemeister moves, detour moves, saddle moves and births and
deaths. This is a generalization of the notion of knot and link cobordism in classical knot theory \cite{FoxMilnor,VKC,VKC1}. Given a cobordism between diagrams $K$ and $K'$ there is 
a schematic surface associated with the cobordism. Individual bits of surface for the saddle and for the birth and death are indicated in Figure~\ref{Saddle}. To conceptualize this 
{\it schematic surface} associated with a cobordism one thinks of each move applied to the diagram as corresponding to the construction of a bit of surface in the cross product of the diagram plane with the unit interval. Virtual crossings are resolved in the schematic surface and we consider its genus. If the surface between $K$ and $k'$ has genus $0$ we call the cobordism a {\it concordance}. If a virtual knot diagram $K$ is concordant to an unknot diagram, we say that $K$ is a {\it slice knot}. \\

\begin{figure}
     \begin{center}
     \begin{tabular}{c}
     \includegraphics[width=8cm]{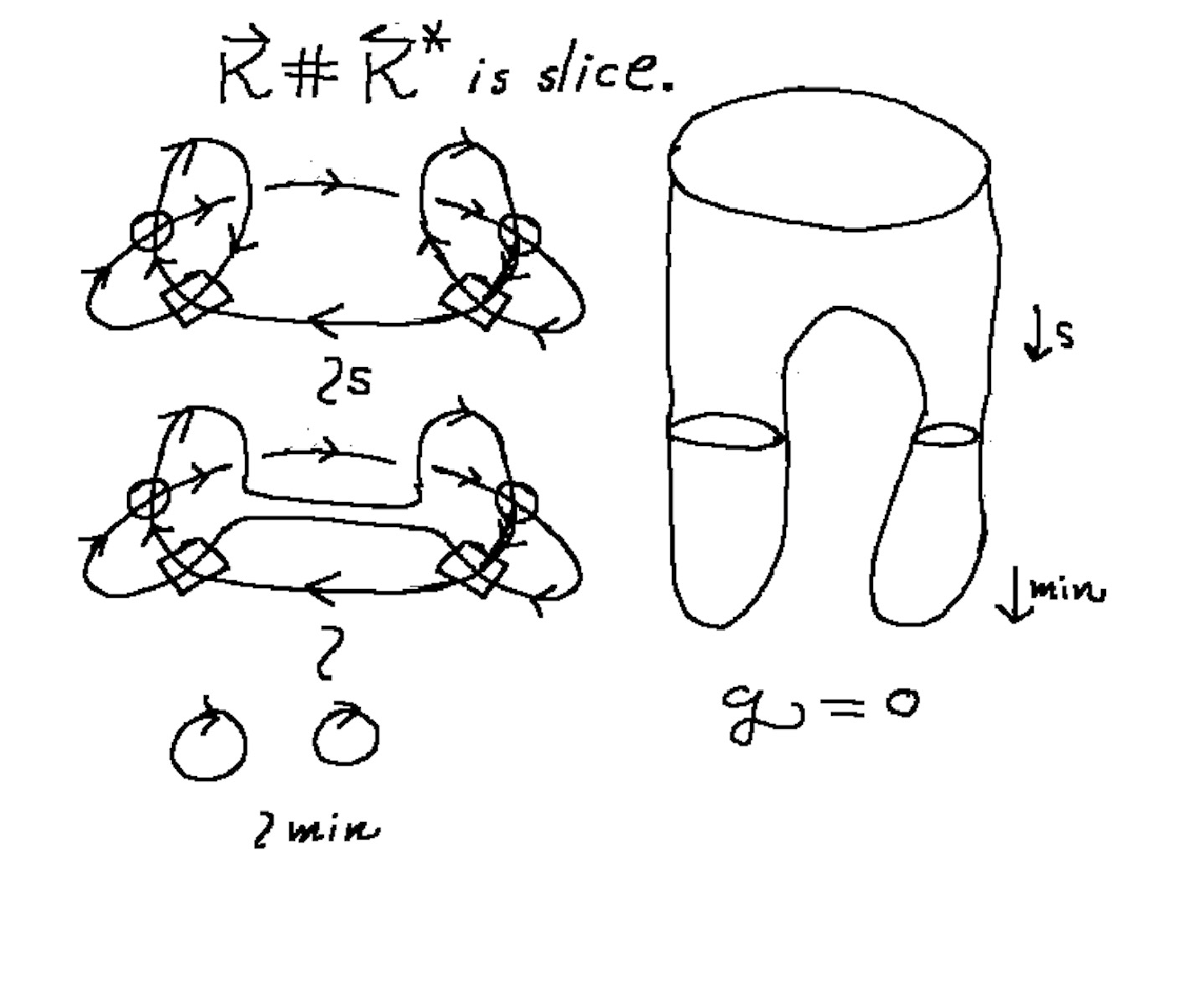}
     \end{tabular}
     \caption{\bf Connected Sum of Knot and Reverse Mirror Image is Slice.}
     \label{ConnSum}
\end{center}
\end{figure}

Figure~\ref{ConnSum} illustrates the general fact that if $\overrightarrow{K}$ is any multiple virtual knot and $\overleftarrow{K}^{\star}$ is its reversed mirror image, then
the connected sum $\overrightarrow{K} \sharp \overleftarrow{K}^{\star}$ is a slice knot. We leave the proof of this fact as an exercise for the reader.\\

\begin{figure}
     \begin{center}
     \begin{tabular}{c}
     \includegraphics[width=8cm]{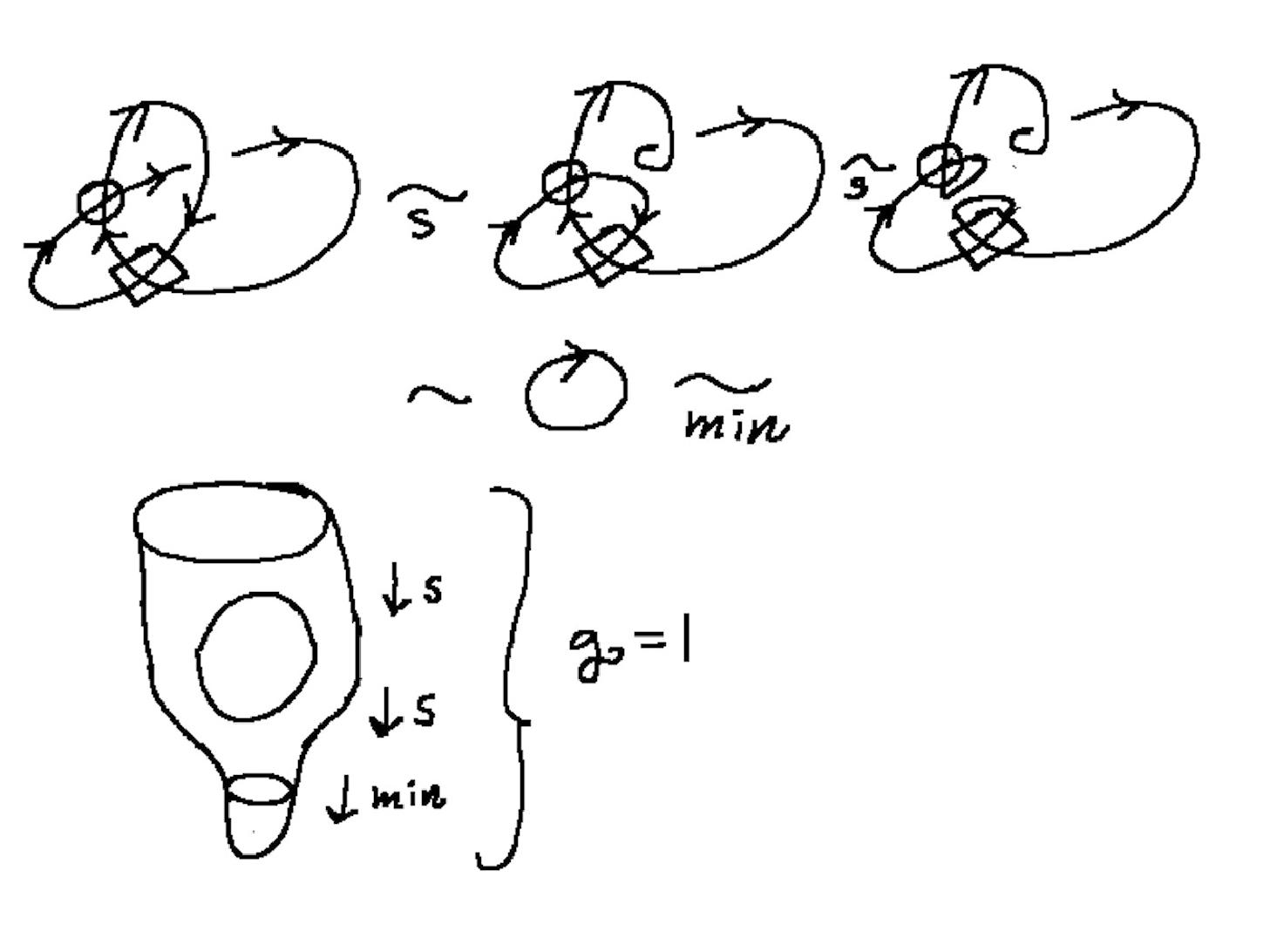}
     \end{tabular}
     \caption{\bf Cobordism for the Traefoil Knot.}
     \label{Traefoil}
\end{center}
\end{figure}

Figure~\ref{Traefoil} illustrates a multiple virtual knot $Traefoil$ with two distinct virtual crossings and one classical crossing  and shows a cobordism of $Traefoil$ to the unknot of genus
one.\\

\begin{figure}
     \begin{center}
     \begin{tabular}{c}
     \includegraphics[width=10cm]{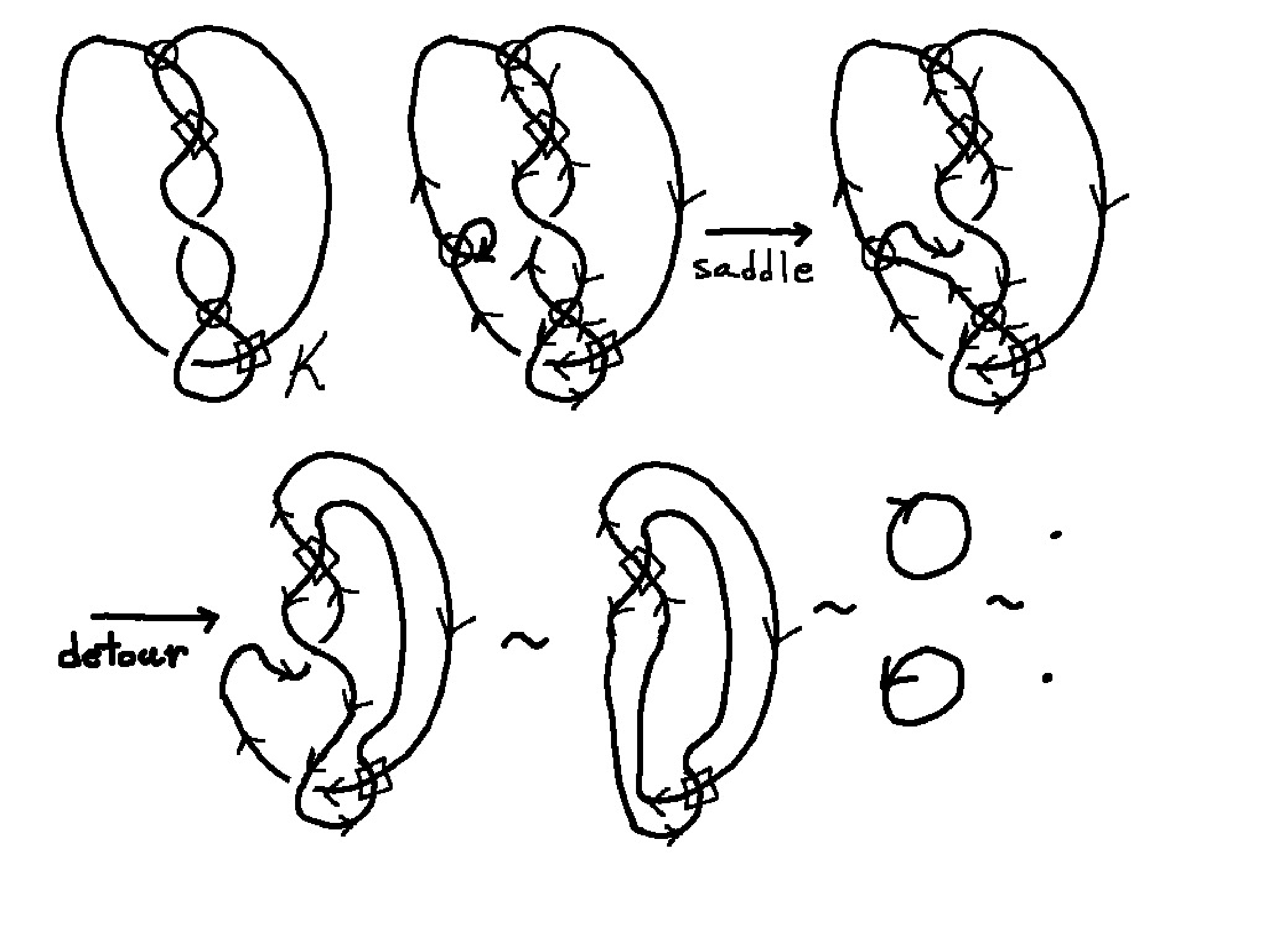}
     \end{tabular}
     \caption{\bf Multiple virtual slice knot K.}
     \label{EFF29}
\end{center}
\end{figure}

Figure~\ref{EFF29} Shows a multiple virtual knot $K$ with two types of virtual crossing. The Figure shows that $K$ is a slice knot.
Figure~\ref{EFF30} shows the bracket expansion for the knot $K$ in Figure~\ref{EFF29}.  As the Figure shows, we calculate that
$\langle K \rangle = A^2 \delta + \delta + B^2 \delta + \Gamma$ where $\Gamma$ is the diagram with four virtual crossings shown in the Figure.
Note that $\delta = -A^2 - A^{-2}$ and $B= A^{-1}.$ We have already proved  that the state $\Gamma$ is non trivial under virtual equivalence.
Figure~\ref{EFF31} accomplishes this task by showing that $\Gamma$ has a non-trivial multiple virtual quandle. Thus the knot $K$ is slice and it is a non-trivial
multiple virtual knot that necessarily carries at least two distinct virtual crossings on all of its diagrams.\\

\section{Parity Bracket Polynomials for Multiple Virtual Knot Theory}
In this section we introduce and generalize the parity bracket polynomial of Vassily Manturov \cite{IP,vkt}.  
Using the parity bracket involves making 4-valent graphical nodes. In this paper we have used such nodes to indicate a fusion evaluation that corresponds to a graphical contraction.
In this section we will not use this convention for the graphical node. Graphs with such nodes will be taken according to equivalence relations as described below.\\

The parity bracket is a generalization of the bracket polynomial to virtual knots and links that uses the {\it parity} of the crossings. Crossing parity is understood as follows: Let $c$ be a classical crossing in a (multi-virtual) diagram. Traverse the diagram from $c$ until one returns to $c$ for the first time. Count the number of classical crossings encountered in this traverse (some crossings will be encountered once, others twice). The parity of this count is the parity of the crossing. Thus if one encounters an even number of crossings in the traverse from $c$ to $c$, then the crossing is {\it even}. If one encounters an odd number of crossings, then $c$ is {\it odd}.\\

\begin{figure}
     \begin{center}
     \begin{tabular}{c}
     \includegraphics[width=8cm]{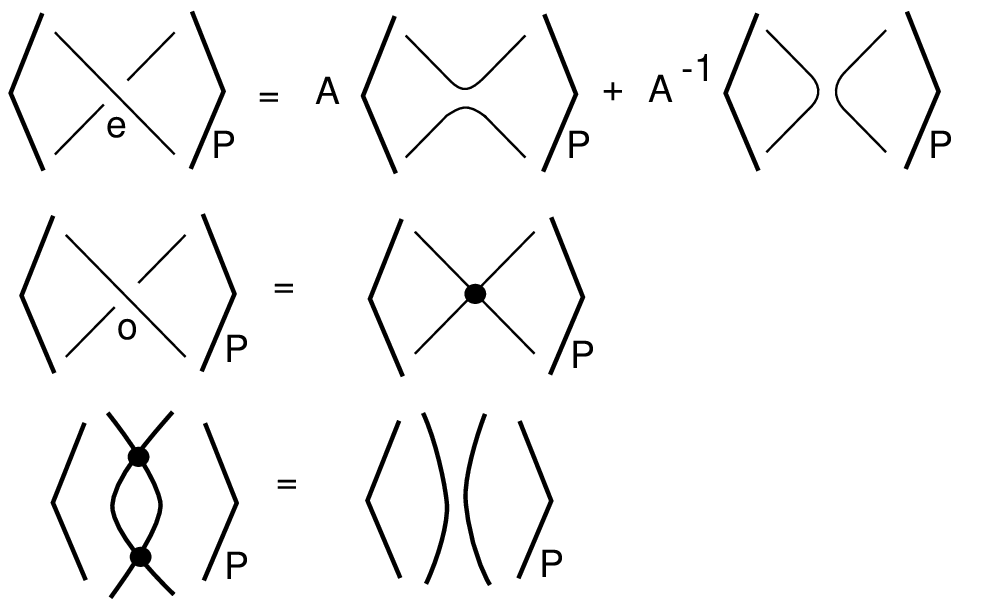}
     \end{tabular}
     \caption{\bf Parity Bracket Expansion}
     \label{Figure 10}
\end{center}
\end{figure}

We define a {\em parity state} of a virtual diagram $K$ to be a labeled virtual
graph obtained from $K$ as follows: For each odd crossing in $K$ replace the crossing by a graphical node. For each even crossing in $K$ replace the crossing by one of its two possible smoothings, and label the smoothing site by $A$ or $A^{-1}$ in the usual way. Then we define the parity bracket by the 
state expansion formula
$$\langle K \rangle _{P} = \sum_{S}A^{n(S)}[S]$$
where $n(S)$ denotes the number of $A$-smoothings minus the number of $A^{-1}$ smoothings and 
$[S]$ denotes the equivalence class of the nodal state $S$ with respect to detour moves and Reidemeister two moves on nodes as shown in Figure~\ref{Figure 10}.  
In the case of single virtual knot theory, a diamond lemma argument \cite{IP,VSA} shows that one can obtain a unique representative for the equivalence class of $S$ by reducing $S$ by Reidemeister two moves on the nodes. In the the multi-virtual case, we only say that we take the equivalence class $[S]$, but in some cases we can determine information about this class and thereby use the invariant.\\

With this we write $$[S] = (- A^{2} - A^{-2})^{l(S)} [G(S)]$$ where $l(S)$ is the number of standard loops in the reduction of the state $S$ and $[G(S)]$ is a  disjoint union of (reduced) graphs that contain nodes.
In this way, we obtain a sum of Laurent polynomials in $A$ multiplying graph classes as the parity bracket. The parity bracket is invariant under regular isotopy and detour moves
and behaves just like the usual bracket under the first Reidemeister move. The definition we have given is a direct generalization of the usual parity bracket for single virtual knot theory.\\

\begin{figure}
     \begin{center}
     \begin{tabular}{c}
     \includegraphics[width=8cm]{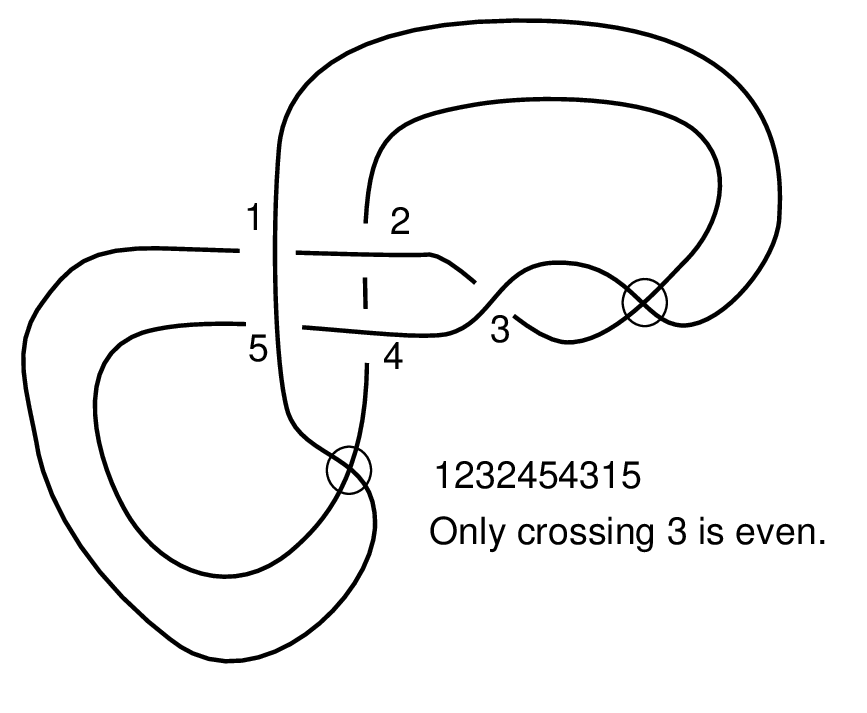}
     \end{tabular}
     \caption{\bf A Knot KS With Unit Jones Polynomial}
     \label{Figure 12}
\end{center}
\end{figure}

\begin{figure}
     \begin{center}
     \begin{tabular}{c}
     \includegraphics[width=8cm]{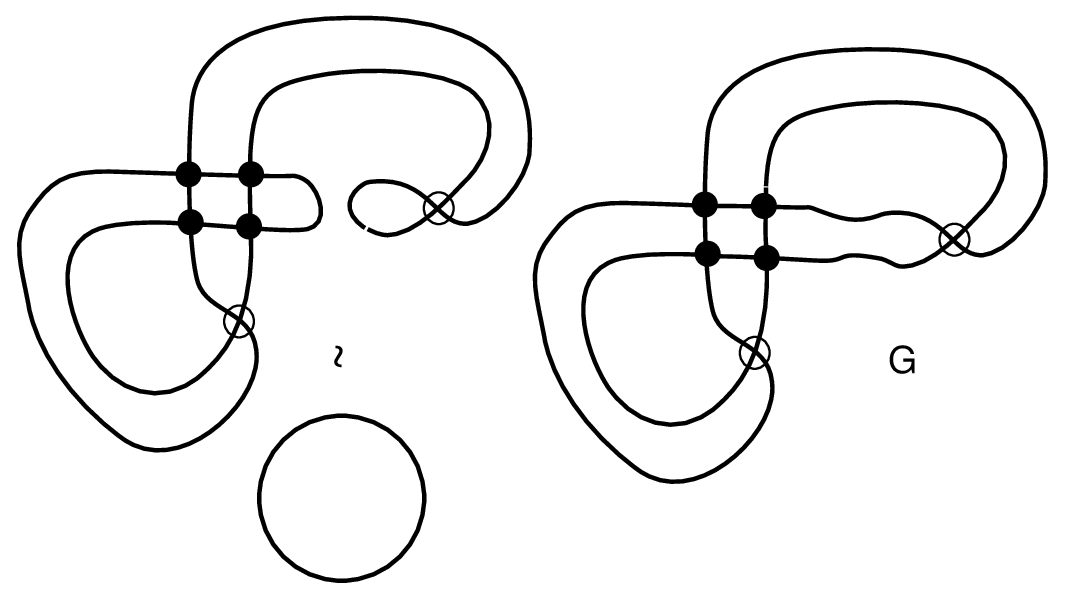}
     \end{tabular}
     \caption{\bf Parity Bracket States for the Knot KS}
     \label{Figure 13}
\end{center}
\end{figure}

\begin{figure}
     \begin{center}
     \begin{tabular}{c}
     \includegraphics[width=8cm]{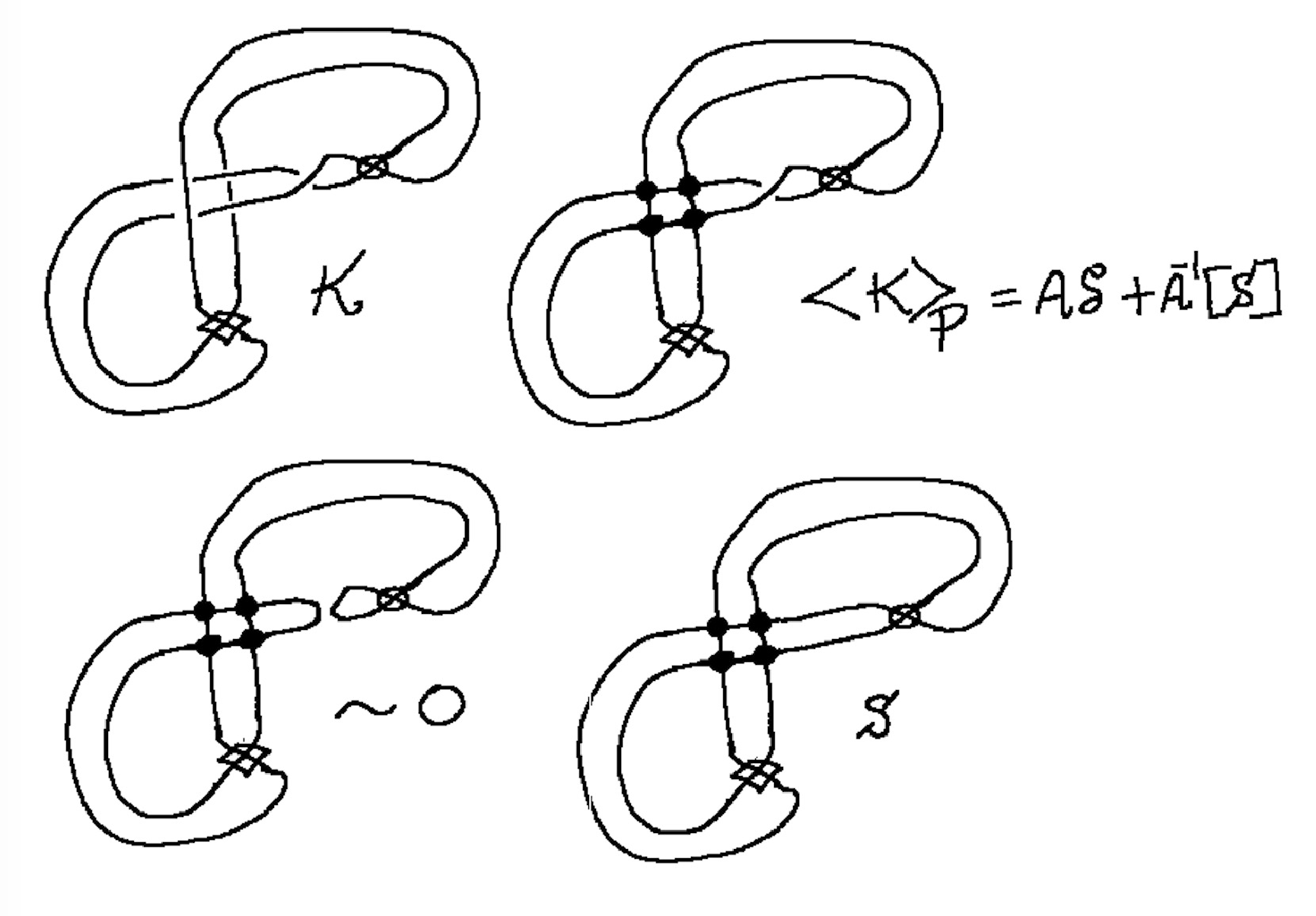}
     \end{tabular}
     \caption{\bf Parity Bracket States for Knot with Two Virtual Types}
     \label{ParityA}
\end{center}
\end{figure}

\begin{figure}
     \begin{center}
     \begin{tabular}{c}
     \includegraphics[width=8cm]{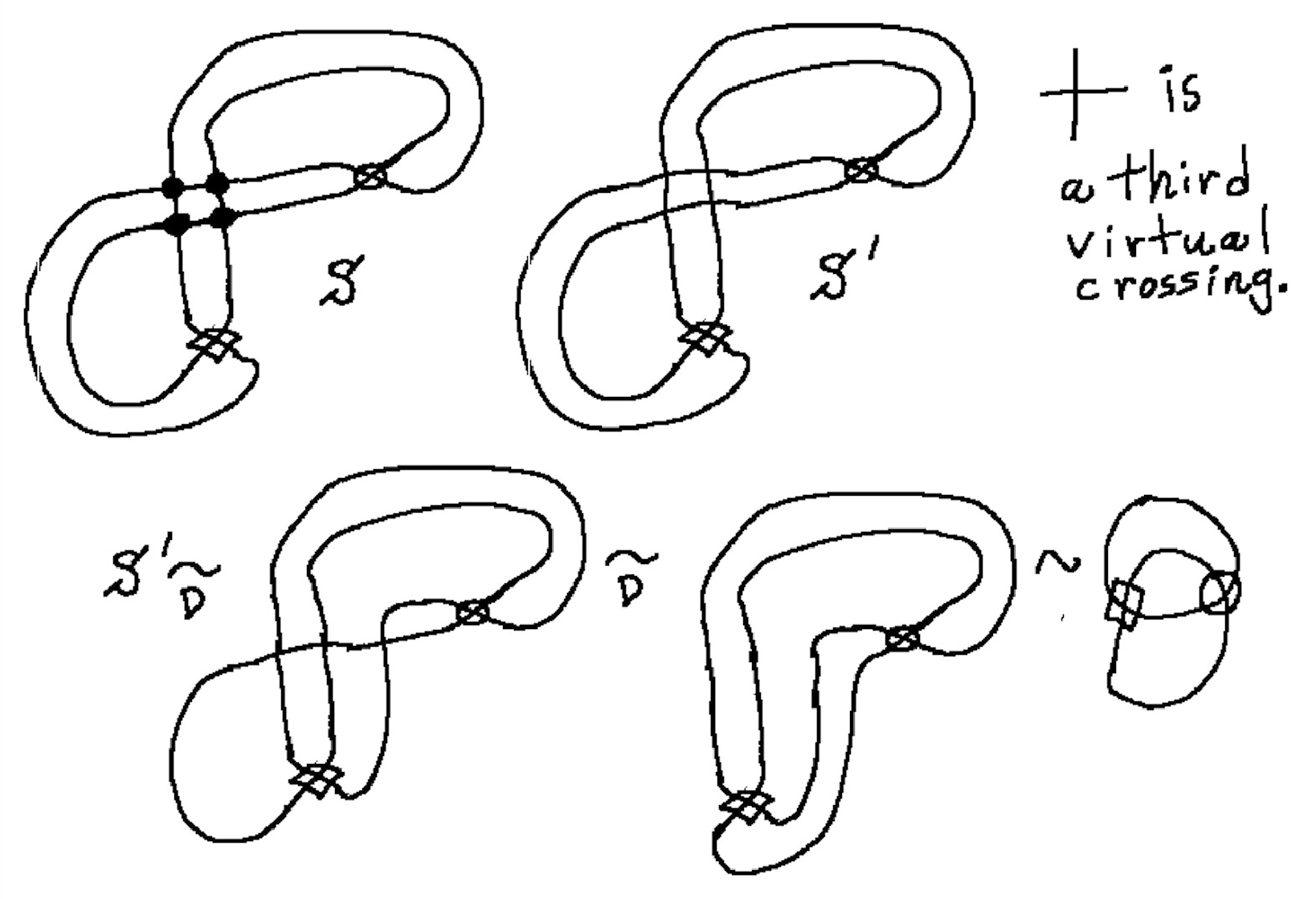}
     \end{tabular}
     \caption{\bf Irreducibility of State S}
     \label{ParityB}
\end{center}
\end{figure}

The use of
parity to make this bracket expand to graphical states gives it considerable extra power in some situations.  
Consider the parity bracket of the knot $KS$ shown in Figure~\ref{Figure 12} and Figure~\ref{Figure 13}. This knot has one even classical crossing and four odd crossings. One smoothing of the even crossing yields a state that reduces to a loop with no graphical nodes, while the other smoothing yields a state that is irreducible. The upshot is that this knot KS is a non-trivial virtual knot. One can verify that $KS$ has unit Jones polynomial. It is a counterexample to a conjecture of Fenn, Kauffman and Maturov \cite{KFM,KFJ} that suggested that a knot with unit Jones polynomial should be $Z$-equivalent to a classical knot. (For the definition of $Z$-equivalence, the reader can consult \cite{KFM,vkt}.)\\

Now examine the knot with two virtual crossings in Figure~\ref{ParityA} and Figure~\ref{ParityB}. This knot is obtained from the previous example by changing one of the round virtual crossings to a box virtual crossing. There is one nodal state to consider. Figure ~\ref{ParityB} shows that this nodal state is irreducible. We replace the nodes by a {\it third} virtual crossing type, and then we show that the replaced state is equivalent by detour moves to a two circle state $\Lambda$ with one round and one box virtual crossing. Since we know that the virtual detour class of $\Lambda$ is non-trivial, this proves that the state $S$ is irreducible, and this proves that the knot in Figure~\ref{ParityA} is a knotted curve in the multi-virtual theory.\\

This last example shows also a new way to use a parity bracket formulation for multi-virtuals. Instead of replacing the odd crossing by graphical nodes, replace them by a new virtual crossing outside the set of virtual crossings being used. Once again define the parity state sum as before
$$\langle K \rangle _{P} = \sum_{S}A^{n(S)}[S]$$
where $[S]$ denotes the detour class of the state with the new virtual crossing in place of the nodes. We have just constructed an example showing how this invariant can be successfully used.
There is much to explore in multi-virtual theory in relation to parity.\\

Parity is an important theme in virtual knot theory and will figure in future investigations of multi-virtuality.  Furthermore the parity notions can be generalized. For more on this theme the reader should consult \cite{MP} and \cite{SL}.\\

\noindent {\bf Remark on Free Knots}
Manturov (See \cite{IP}) studies {\it free knots}. A free knot is a Gauss diagram
(or Gauss code) without any orientations or signs, taken up to abstract Reidemeister moves. Furthermore,
the parity bracket polynomial evaluated at $A = 1$ or $A = -1$ is an invariant of free knots. By using it
on examples where all the crossings are odd, one obtains infinitely many examples of non-trivial free knots.  Any free knot that is shown to be non-trivial has
a number of non-trivial virtual knots overlying it. The free knots can be studied for their own sake. In order to consider a multi-virtual analog for free knots, one can add chords to chord diagrams that correspond to different virtual crossing types, or one can use diagrams as we do in this paper with classical crossings replaced by immersion crossings, and add the $Z$-move to the equivalence relation where the $Z$-move allows the clasical and any virtual crossing that are adjacent to one another to undergo and interchange. This generalization will be pursued separately from the present paper.\\

\section{Generalized Arrow Polynomial}

\begin{figure}
     \begin{center}
     \begin{tabular}{c}
     \includegraphics[width=10cm]{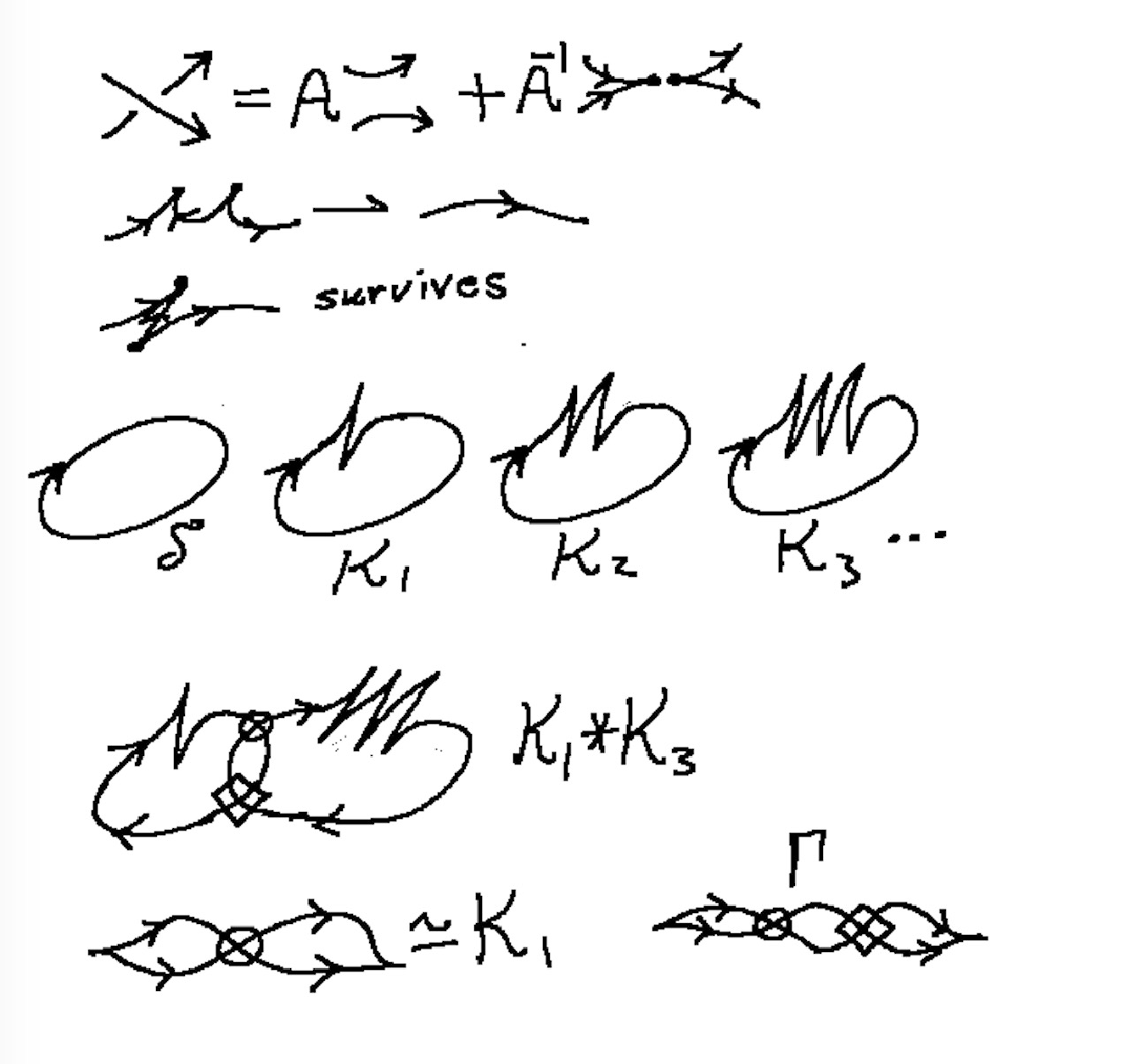}
     \end{tabular}
     \caption{\bf Generalized Arrow Polynomial}
     \label{EFF16}
\end{center}
\end{figure}

\begin{figure}
     \begin{center}
     \begin{tabular}{c}
     \includegraphics[width=8cm]{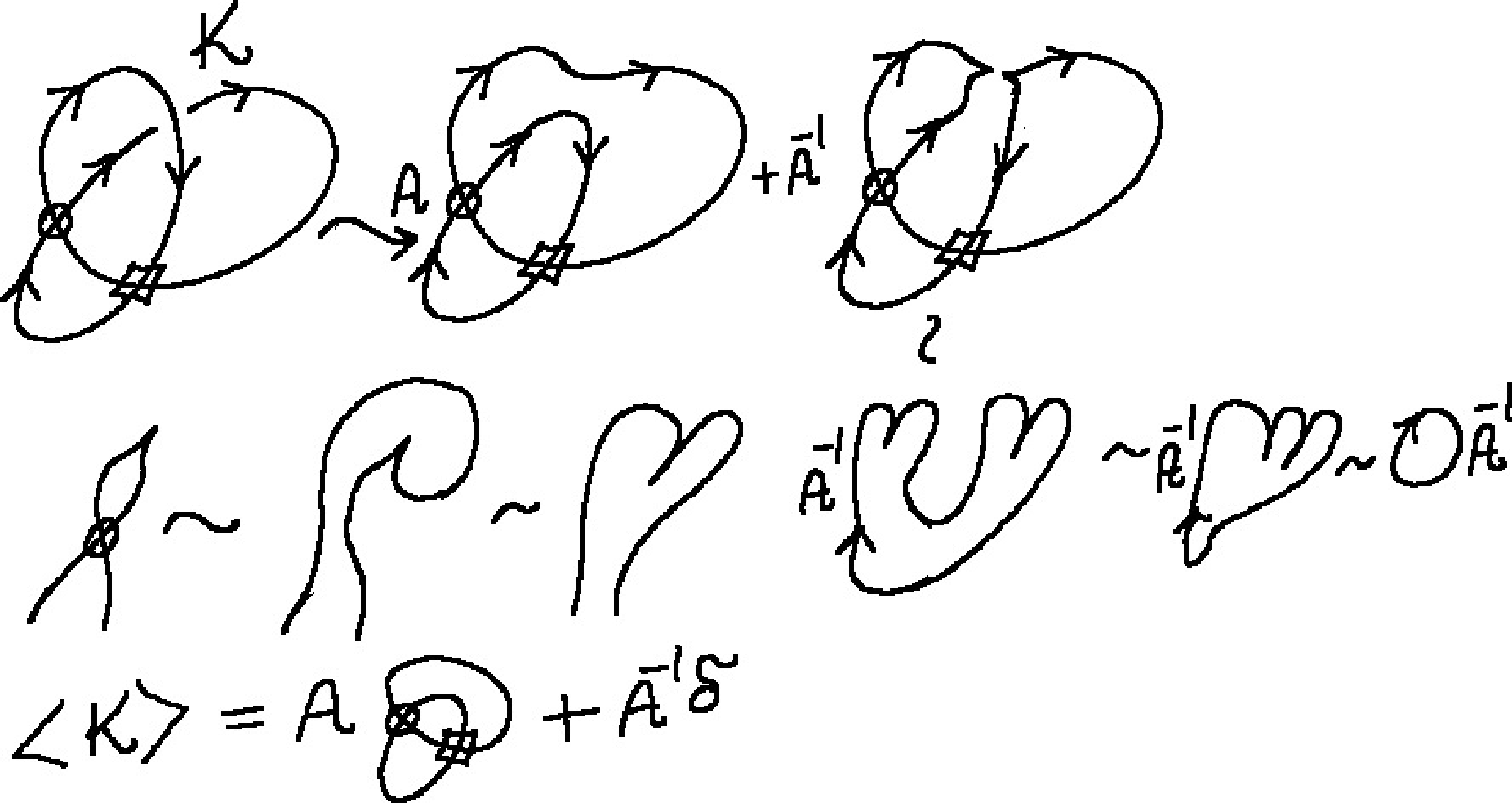}
     \end{tabular}
     \caption{\bf Arrow Expansion}
     \label{EFF17}
\end{center}
\end{figure}

We examine a generalization of the arrow polynomial \cite{vkt}. 
The arrow polynomial uses the bracket expansion on an oriented diagram. This means that there is, at each crossing, an oriented smoothing and a disoriented smoothing, as
shown in Figure~\ref{EFF16}. In that figure we show the bracket expansion at the top of the figure. Note that in the disoriented smoothing we see two cusps where by a {\it cusp}
I mean two opposite orientations meeting at a point and such that the local drawing shows the point of cusp on one side of the given state curve. Thus the cusp indicates a locus
of orientations meeting and it indicates a side of the curve. Analysis of the behaviour of this expansion shows that consecutive same-side cusps need to cancel in order to insure the invariance under the second Reidemeister move. This cancellation is indicated in the figure, and we also indicate how a ``zig-zag" of opposite-side cusps does not cancel.
As a result, we have infinitely many single loop reduced states $K_i$ where $i$ denotes the number of zig-zags in the given loop. These become new variables in the invariant.
In working with multiple virtuals there are many more reduced states as shown at the bottom of the figure. These consist in loops with zig-zags that are bound to one another due to virtual linkage of the multiple virtual crossings. Thus, the generalized arrow polynomial is a summation over the equivalence classes of its states, with polynomial coefficients.
See Figure~\ref{EFF17} for an example calculation. In this example we show how linked virtual states can occur in the detection of a multiple virtual knot.\\

\clearpage

\begin{figure}
     \begin{center}
     \begin{tabular}{c}
     \includegraphics[width=10cm]{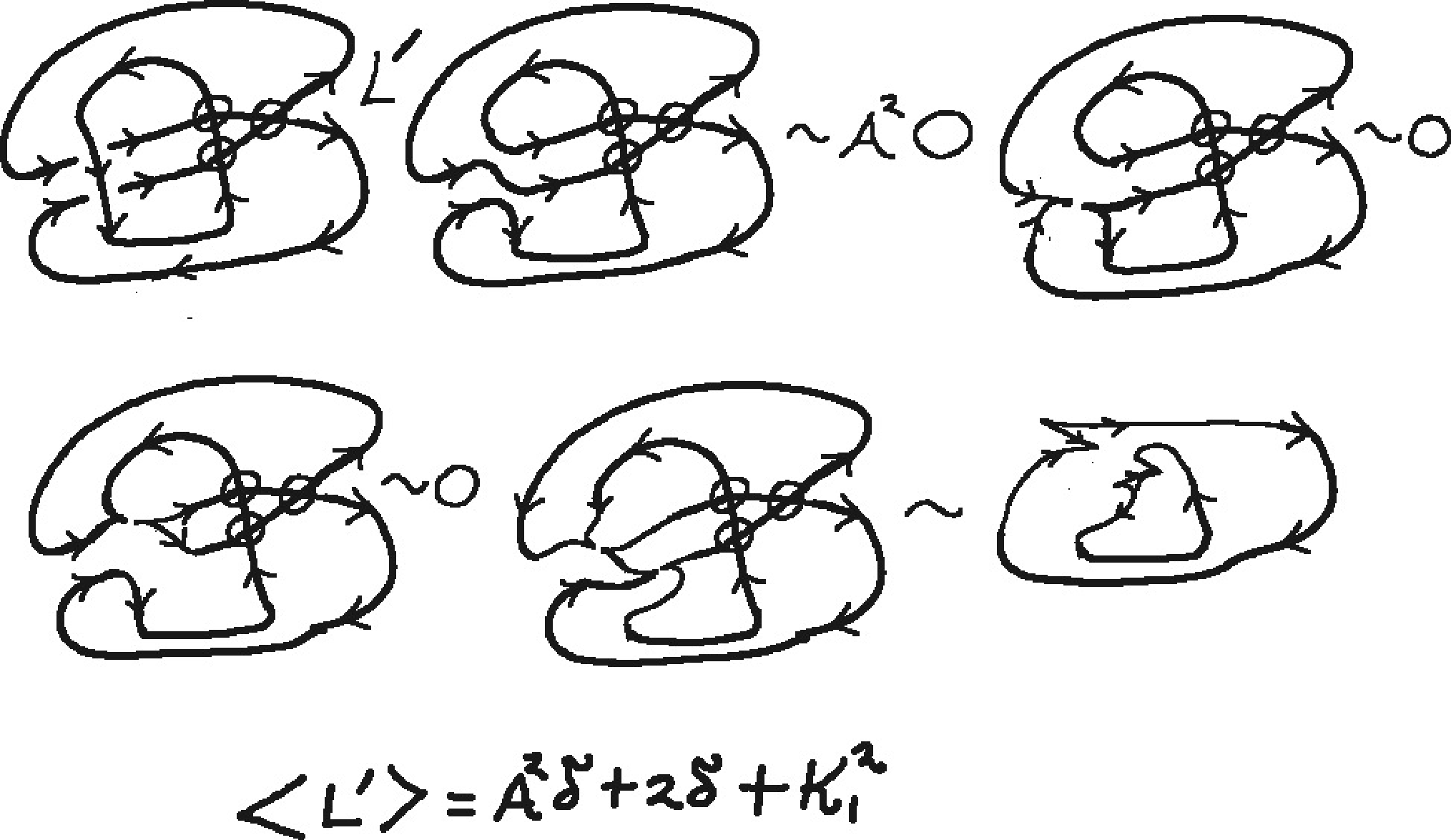}
     \end{tabular}
     \caption{\bf Simple Link}
     \label{MVLink0}
\end{center}
\end{figure}

\begin{figure}
     \begin{center}
     \begin{tabular}{c}
     \includegraphics[width=10cm]{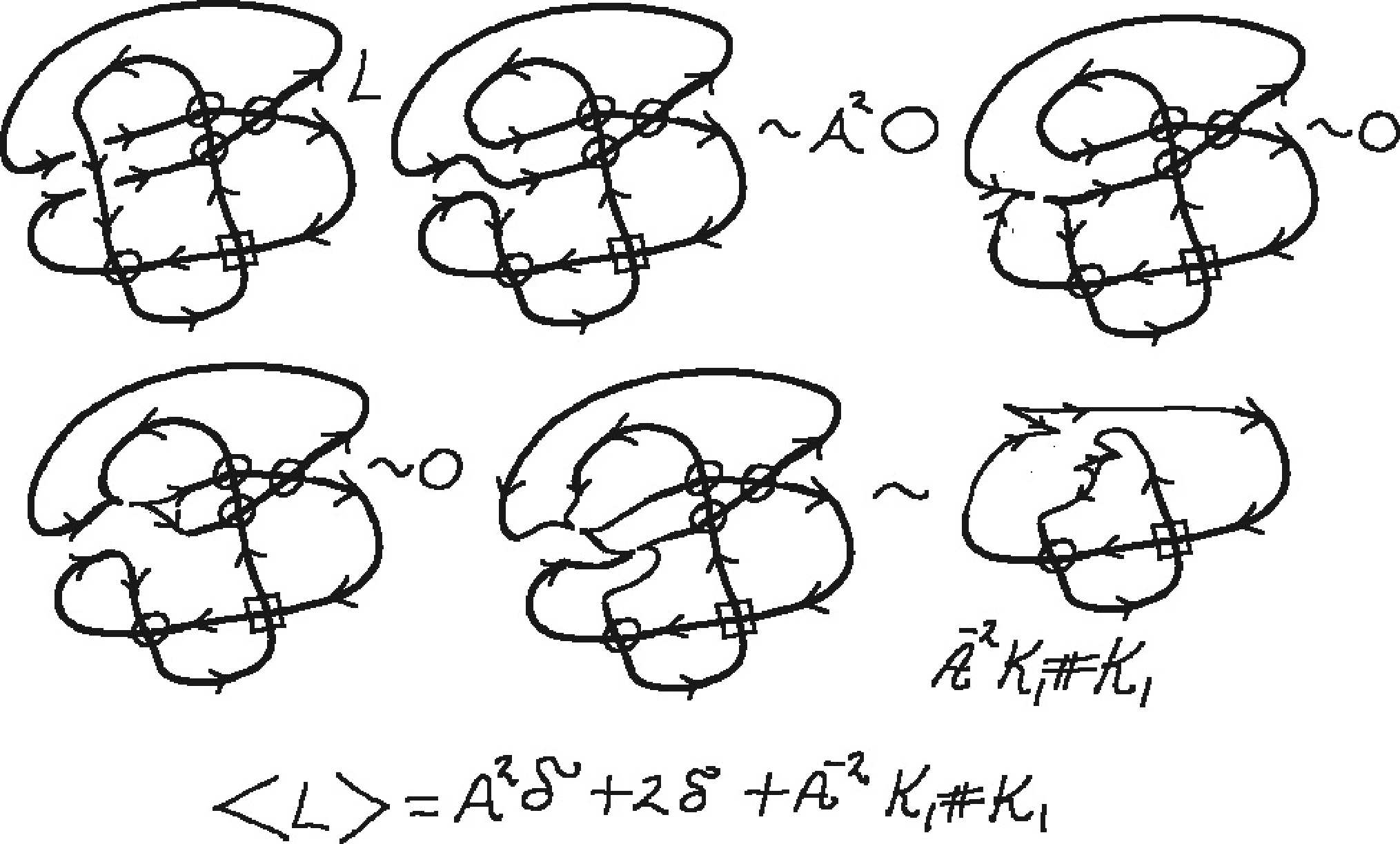}
     \end{tabular}
     \caption{\bf Complex Link}
     \label{MVLink1}
\end{center}
\end{figure}

In Figures~\ref{MVLink0} and ~\ref{MVLink1} we illustrate two arrow polynomial calculations. In the first figure we have a link with one type of virtual crossing and find that it has 
arrow polynomial equal to $A^2 \delta + 2 \delta + K_1^{2}.$ In Figure~\ref{MVLink1} we have a similar link with a clasp consisting of a round virtual crossing and a box virtual crossing. The figure ilustrates how the resulting polynomial has an irreducible state that we call $K_1 \sharp K_1.$ The polynomial is $A^2 \delta + 2 \delta + K_1 \sharp K_1.$ 
Here we evaluate the arrow polynomial in terms of reduced representatives for its states. We do not expand the box virtual as we did for the bracket. This expansion can be used if there are only two types of virtual crossing present, but we wish to illustrate here the pattern of the general case. Note that $\delta = -A^2 - A^{-2}$ as in the bracket. It is
a problem of interest to classify the irreducible states that can occur in the arrow polynomial for general multiple virtual theory.\\

\section{Multiple Welded Knot Theory}
\begin{figure}
     \begin{center}
     \begin{tabular}{c}
     \includegraphics[width=8cm]{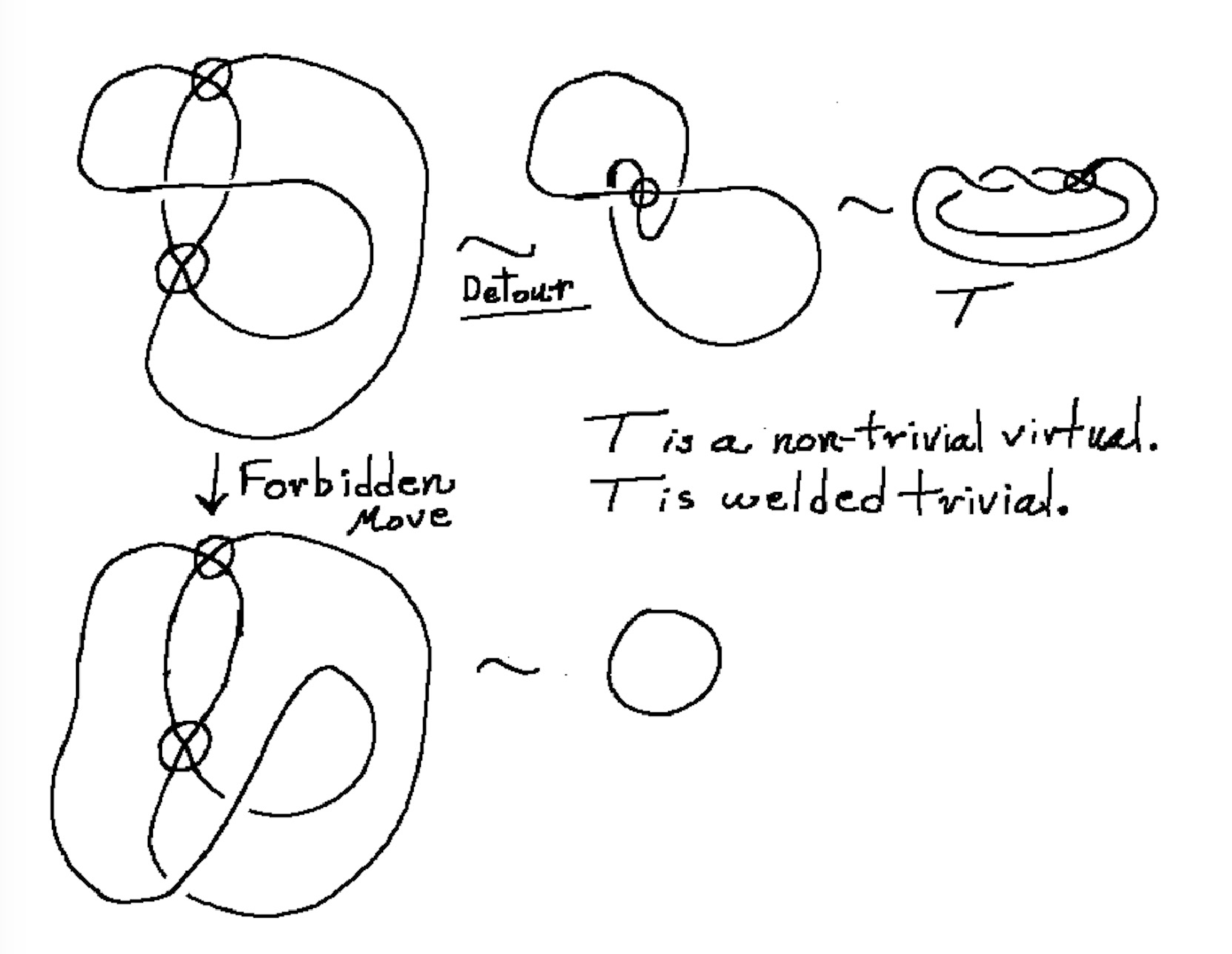}
     \end{tabular}
     \caption{\bf Virtual Non Triviality and Welded Triviality}
     \label{EFF18}
\end{center}
\end{figure}

\begin{figure}
     \begin{center}
     \begin{tabular}{c}
     \includegraphics[width=8cm]{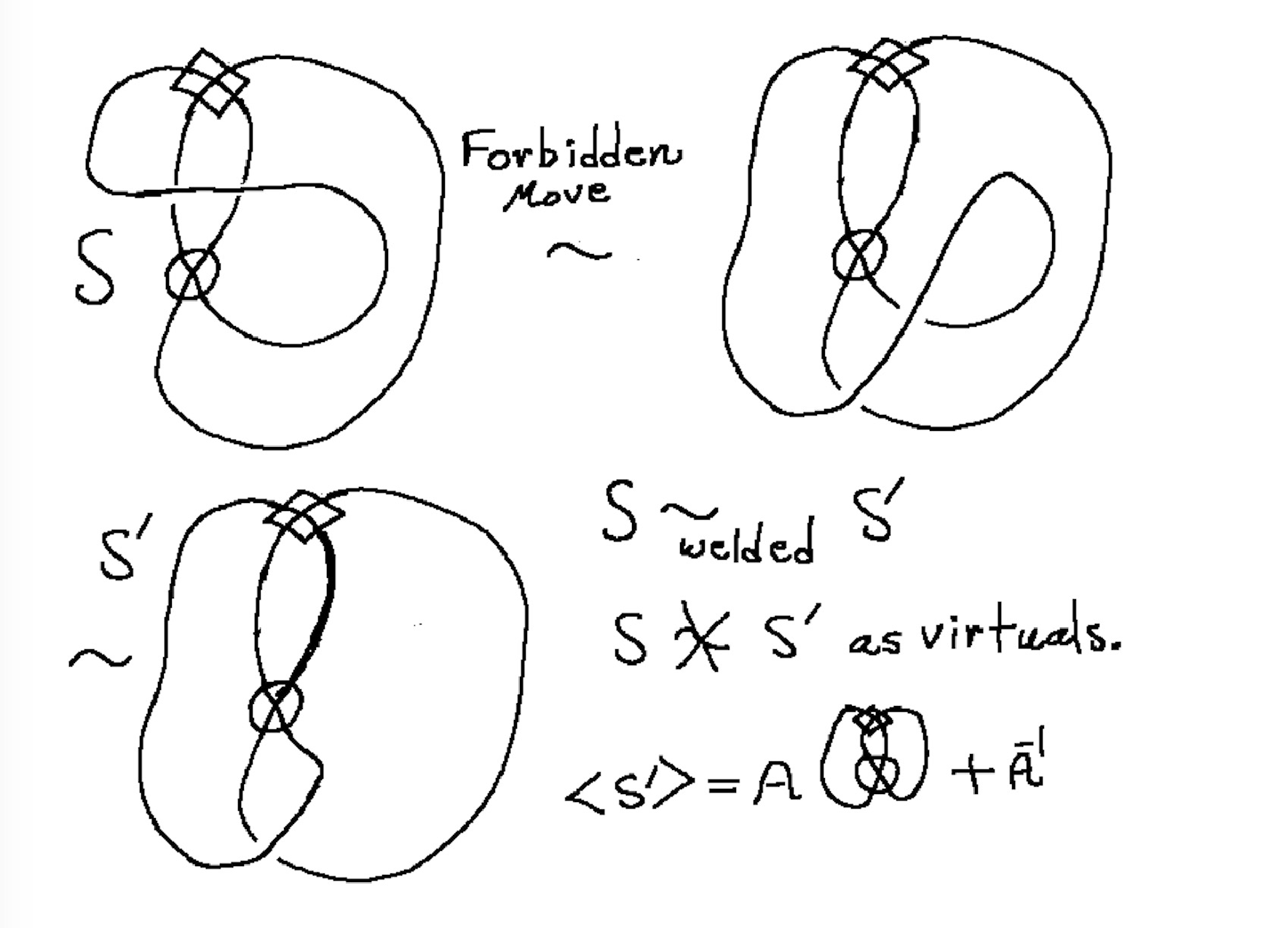}
     \end{tabular}
     \caption{\bf Generalized Virtual and Generalized Welded}
     \label{EFF19}
\end{center}
\end{figure}

In Figure~\ref{EFF18} we remind the reader of the definition of {\it welded knot theory} where we add the move $F1$ (the over-forbidden move) to the already given moves of 
virtual knot theory. The figure shows how the non-trivial virtual trefoil knot becomes a trivial knot in welded knot theory.  Figure~\ref{EFF19} gives an example of two multiple virtual knots $S$ and $S'$ such that a generalized bracket calculation shows that $S$ and $S'$ are distinct virtual knots, but equivalent welded knots. \\

\begin{figure}
     \begin{center}
     \begin{tabular}{c}
     \includegraphics[width=8cm]{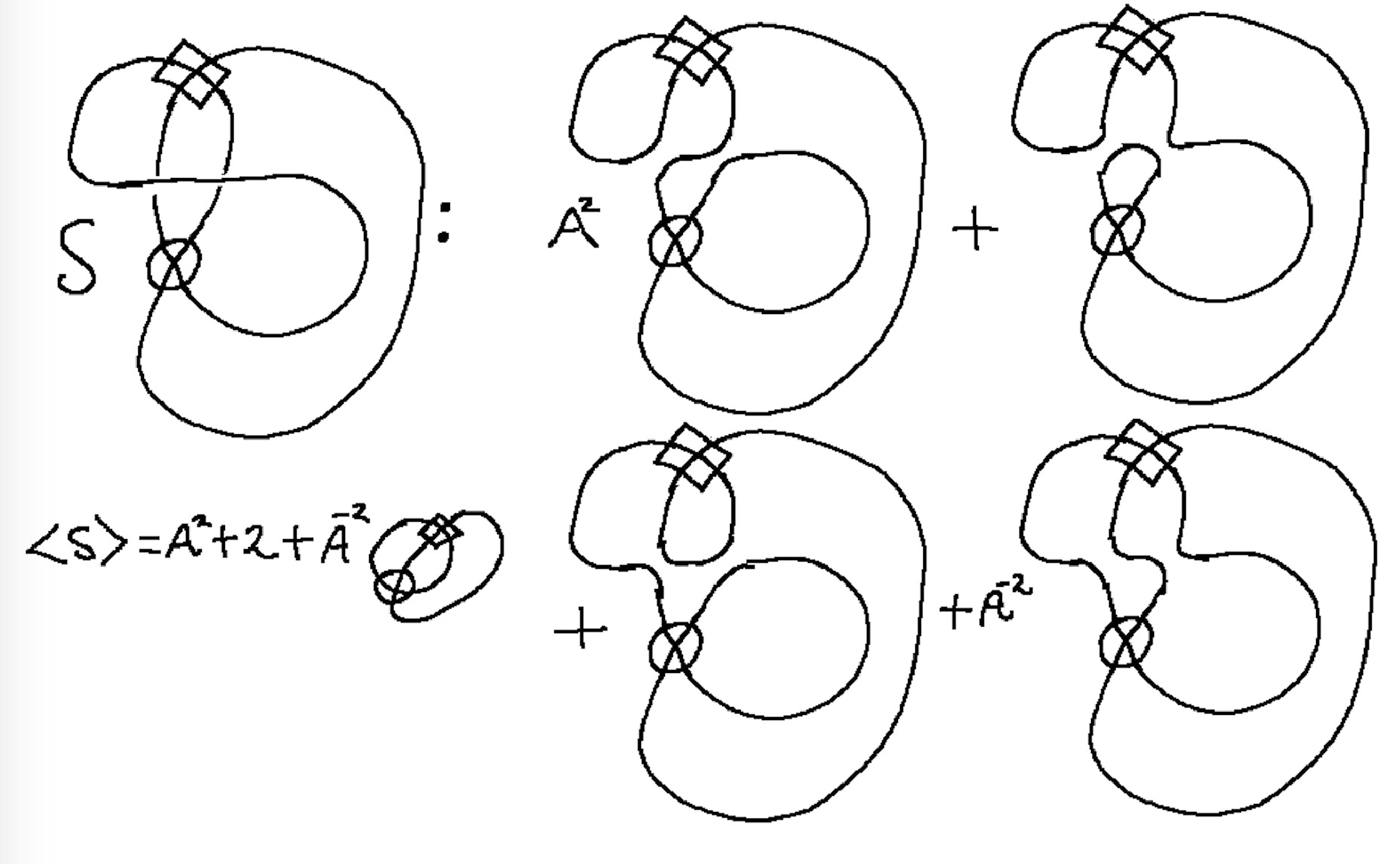}
     \end{tabular}
     \caption{\bf A Bracket Expansion}
     \label{EFF20}
\end{center}
\end{figure}

Just as we have opened the domain of multiple virtual knot theory, there is a correspondingly opened domain of multiple welded knot theory with many virtual crossings.
Welded knot theory is based on allowing the forbidden move $F1$ but not allowing the move $F2.$ See Figure~\ref{EFF19} and \ref{EFF20}.\\

We can also consider allowing both $F1$ and $F2$ in a {\it Free Multiple Welded 
Theory}. See Figure~\ref{EFF21} for an illustration of the generalized forbidden moves of type $F1$ and $F2.$
We conjecture that the knot shown in Figure~\ref{EFF19} is non-trivial in free welded theory.\\

Standard welded knot theory with a single type of virtual crossing has a (partial) interpretation in terms of ribbon torus embeddings in four dimensional space
\cite{Satoh,Rourke,Ogasa}. We can indicate the question: {\it Is there a way to intrepret multiple welded virtual links in terms of surfaces in four-space?} \\

\clearpage

\section{Rotational Multiple Virtuals}
In \cite{vkt,rotvkt} we introduce Rotational Virtual Knot Theory, where the detour moves are restricted to be regular isotopies in the plane or on the two dimensional sphere.
A feature of this theory is that essentially all quantum link invariants for classical knots and links generalize to the rotational virtuals. In \cite{rotvkt} we show limitations on quantum link invariants for detecting rotational virtuals. The rotational theory generalizes to multi-virtuals, and there are many problems to explore, \\ 

\begin{figure}
     \begin{center}
     \begin{tabular}{c}
     \includegraphics[width=12cm]{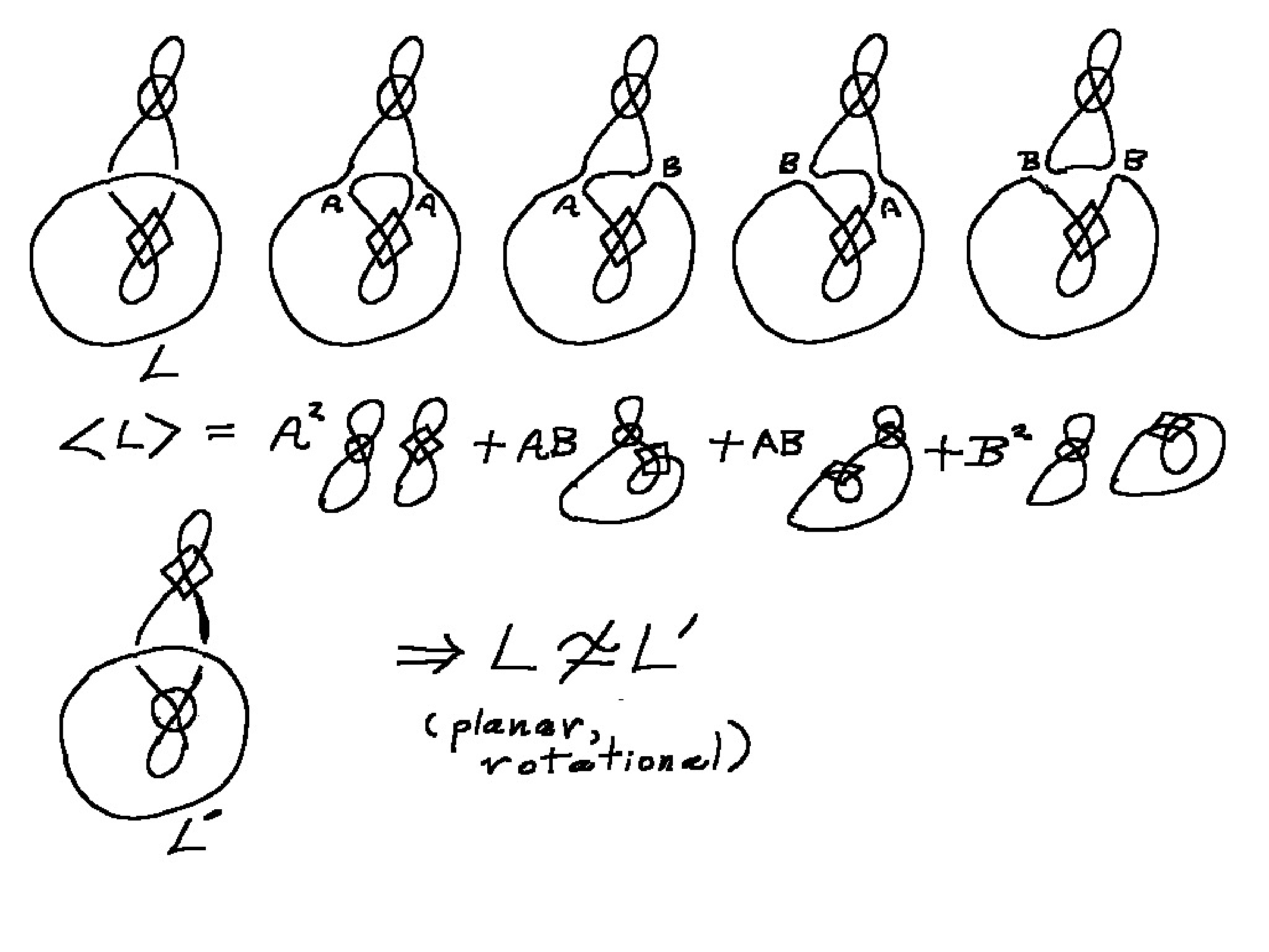}
     \end{tabular}
     \caption{\bf Planar Rotational Example}
     \label{PR}
\end{center}
\end{figure}

In this section we give one example of a calculation of the generalized rotational bracket polynomial. This invariant is defined for multi-virtuals just as the generalized bracket except that the states are taken up to rotational detour
equivalence. This means that detours involve regular homotopy and the Whitney degree for plane curves is an invariant when we are using the planar version of the bracket.
In Figure~\ref{PR} we show examples of multiple virtual links $L$ and $L'$ that are planar rotationally distinct  as shown by the computation of the multiple rotational virtual bracket.\\

Further exploration of rotational multiple virtuals will occur in subsequent work.\\

 \section{The Multiple Virtual Braid Group.}
\begin{figure}
     \begin{center}
     \begin{tabular}{c}
     \includegraphics[width=12cm]{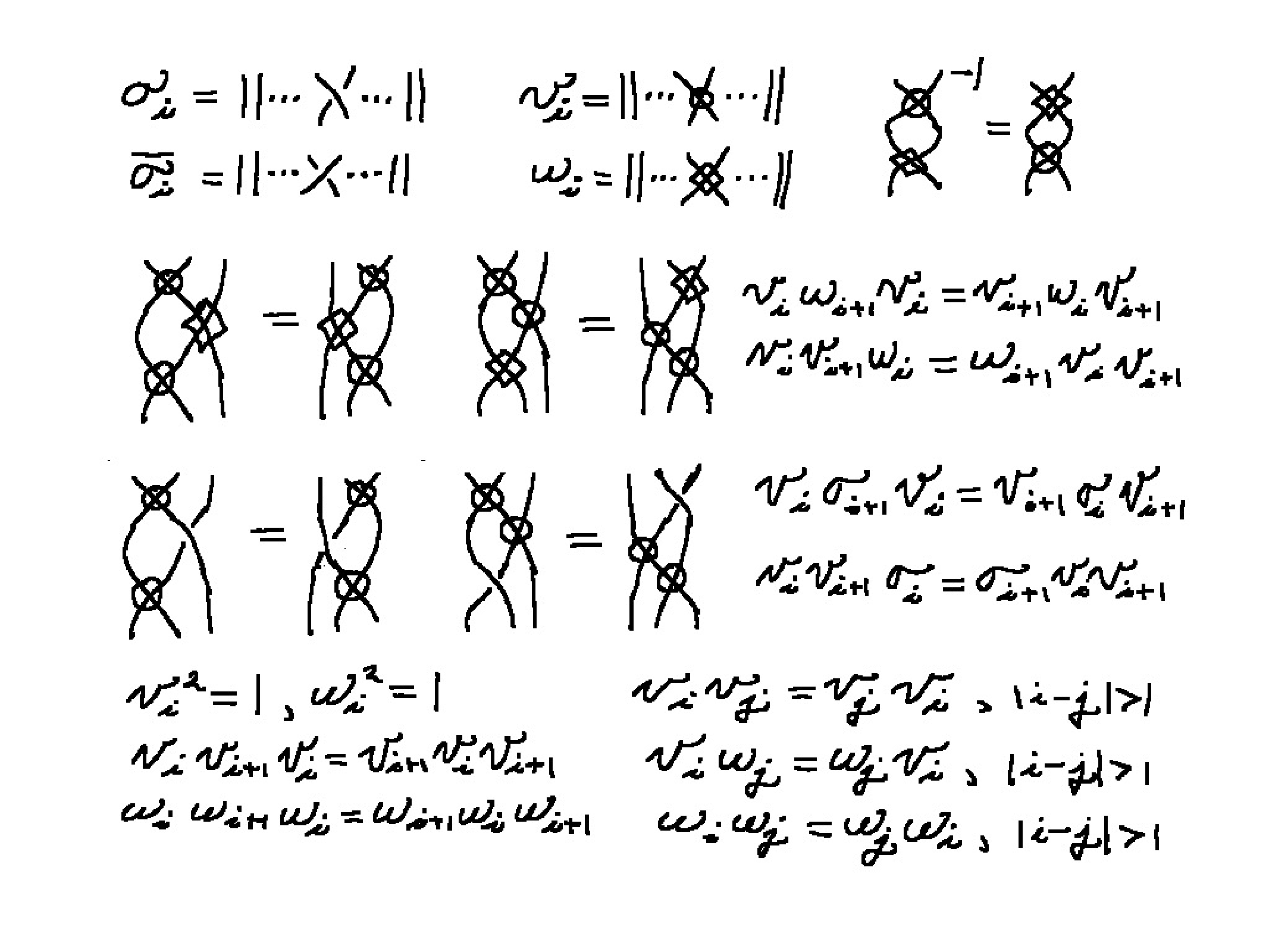}
     \end{tabular}
     \caption{\bf Multiple Virtual Braid Group}
     \label{braid}
\end{center}
\end{figure}

The multiple virtual braid group is defined in direct analogy with the virtual braid group.
See Figure~\ref{braid}.\\

We have virtual braiding generators for each type of virtual crossing and relations corresponding to the detour moves.
We have braiding generators for the classical crossings, and braiding relations as in the classical braid group.
This means that we have braiding generators $\sigma_{i}$ and virtual generators $v[\alpha]_{i}$ where $\alpha$ runs over some index set for the virtual generators.
Let $v_{i}$ and $w_{i}$ denote two distinct virtual generators ($v = v[\alpha], w = w[\beta]$). Then the relations for the multiple virtual braid group are:
$$\sigma_{i} \sigma_{i+1} \sigma_{i } =  \sigma_{i+1} \sigma_{i} \sigma_{i +1},$$
$$\sigma_{i} \sigma_{j} = \sigma_{j} \sigma_{i}, |i-j|>1,$$
$$v_{i}^2 = w_{i}^2 = 1,$$
$$v_{i} v_{j} = v_{j} v_{i},  |i-j|>1,$$
$$v_{i} w_{j} = w_{j} v_{i},  |i-j|>1,$$
$$v_{i} v_{i+1} v_{i } =  v_{i+1} v_{i} v_{i +1},$$
$$v_{i} v_{i+1} w_{i } =  w_{i+1} v_{i} v_{i +1},$$
$$v_{i} w_{i+1} v_{i } =  v_{i+1} w_{i} v_{i +1},$$
$$v_{i} \sigma_{i+1} v_{i} = v_{i+1} \sigma_{i} v_{i+1},$$
$$v_{i} v_{i+1} \sigma_{i} = \sigma_{i+1} v_{i} v_{i+1},$$
$$v_{i} w_{i+1} v_{i} = v_{i+1} w_{i} v_{i+1},$$
$$v_{i} v_{i+1} w_{i} = w_{i+1} v_{i} v_{i+1}.$$

 Note that for a finite number of strands $n$ the group $MVB_{n}$ involves the classical braid group $B_{n}$ and as many copies of the symmetric group $S_{n}$ as there are 
 virtual crossing types. Thus there is a rich algebraic structure in the multiple virtual braid group and even in the subgroups of these braid groups that have no braiding, just the intertwining of the multiple copies of the symmetric group.\\
 
\subsection {Algebraic Markov Equivalence for Multiple Virtual Braids}
 We will prove Alexander and Markov Theorems for the Multiple Virtual Braid Group $MVB$ in a subsequent paper.
In this section we state the Markov Theorem. The proof of this theorem can be accomplished by generalizing the work of 
Kauffman and Lambropoulou in \cite{L,LB}.\\

Let  $MVB_{n}$ denote a multiple virtual braid group on
$n$ strands and let $\sigma_i, (v_i)^{\alpha}$ be its generating classical and virtual crossings. The superscript $\alpha$ indicates the virtual crossing type. The $\sigma_i$'s
satisfy the relations of the classical braid group and the $(v_i)^{\alpha}$'s satisfy the relations of the permutation
group and the detour relations for virtual crossings. The characteristic relation in $MVB_{n}$  is the {\it special detour move} relating both:
 
$$ \begin{array}{cccl} 
(v_{i})^{\alpha} \sigma_{i+1} (v_{i})^{\alpha} & = & (v_{i+1})^{\alpha} \sigma_i (v_{i+1})^{\alpha}. &   \\ 
\end{array}$$

\noindent The group $MVB_{n}$ embedds naturally into $MVB_{n+1}$ by adding one identity strand at the right of
the braid. So, it makes sense to define $MVB_{\infty} :=
\bigcup_{n=1}^{\infty} MVB_{n}$,  the disjoint union of all multiple virtual braid groups (for a given choice of virtual crossing types). We can now state our result.\\

\begin{figure}
     \begin{center}
    $$\includegraphics[width=10cm]{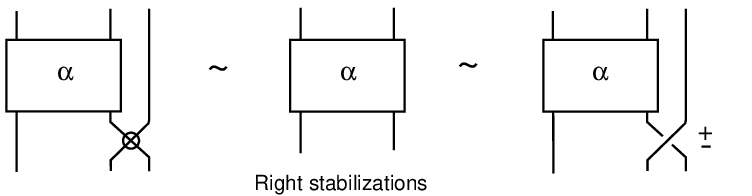}$$
     $$\includegraphics[width=10cm]{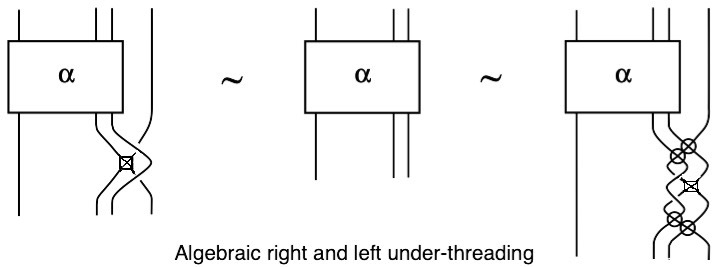}$$
     \caption{\bf The Moves (ii), (iii) and (iv) of the Algebraic Markov Theorem.}
     \label{L}
\end{center}
\end{figure}

\noindent {\bf Theorem. (Algebraic Markov Theorem for multi-virtuals).} Two oriented  multi-virtual links are isotopic if and
only if any two corresponding virtual braids differ by a finite sequence of braid relations in $MVB_{\infty}$ and the
following moves or their inverses. In the statement below and in Figure~\ref{L}, $v_{n}$ stands for any given virtual crossing type.

\begin{itemize}
\item[(i)]Virtual and real conjugation:    \ \ \ \ \ \ \ \ \ \ \ \  $ v_i \alpha v_i \sim \alpha \sim 
{\sigma_i}^{-1}\alpha \sigma_i $
\item[(ii)]Right virtual and real stabilization:  \ \ \ \ \  $\alpha v_n \sim \alpha
\sim \alpha \sigma_n^{\pm 1}$ 
\item[(iii)]Algebraic right under--threading:  \ \  $\alpha \sim \alpha \sigma_n^{-1} v_{n-1} \sigma_n^{+1} $
\item[(iv)]Algebraic left under--threading:  \ \ \ \ $\alpha \sim  \alpha v_n v_{n-1} \sigma_{n-1}^{+1} (v_n)'
\sigma_{n-1}^{-1} v_{n-1} v_n $,
\end{itemize} 

\noindent where $\alpha,  v_i, \sigma_i \in VB_n$ and  $v_n, \sigma_n \in VB_{n+1}$ (see Figure~\ref{L}) and $(v_n)'$ denotes a possibly different virtual crossing type from $v_n$ . Note that in Figure~\ref{L} this possible difference in virtual crossing type is indicated by a box at the crossing rather than a circle.\\

The proof of this theorem is a generalization of the work in \cite{L}, and this theorem and its consequences will be the subject of a separate paper by Kauffman and Lambropoulou. Having a Markov Theorem is a first step toward 
constructing new invariants of multiple virtual knots and links via the properties of the braid category.\\

 \section{Epilogue}
 This paper has been about a generalization of virtual knot theory with many types of virtual crossings, all allowed to detour with respect to one another. We have endeavored to show how this kind of virtual knot theory is related to coloring problems in graph theory and to many interesting topological and combinatorial questions related to the generalizations of knot and link invariants. We have raised many questions and the reader may feel satisfied in working on some of the 
 problems we have raised. While we have discussed a number of topics, some have not been developed - such as the tangle category for multiple virtuals and quantum link invariants for rotational multiple virtuals. Needless to say, there are many
 questions related to link homology for multiple virtual knots.\\
 
Nevertheless, the generalizations of virtual knot theory that we have discussed herein are not the only possibliities. The reader may also be interested in trying other avenues to explore. For example, one can consider knot diagrams with more than one type of classical (topological) crossing. Each such crossing type can be allowed Reidemeister moves with respect to other crossings of its kind, but no interactions with crossings not of its kind. Virtual crossings need not detour over each other,.
There can be more specific rules for allowing detours. Relationships with knots in thickened surfaces can be explored, and there remains the possibility of representations in terms of surfaces in four-space. We encourage the reader to 
explore this territory and to enjoy the freedom of construction that is suggested by these actions.\\

\end{document}